\documentclass[leqno]{amsart}
\usepackage{amsmath,amssymb}
\usepackage{booktabs,mathrsfs,enumerate,mathtools}
\usepackage{etoolbox}
\usepackage[svgnames]{xcolor}

\usepackage{hyperref,framed}
\hypersetup{pdftitle={Quiver Hecke algebras for alternating groups},
  pdfauthor={Clinton Boys and Andrew Mathas},
  colorlinks = true,
  linkcolor =Green,
  anchorcolor = red,
  citecolor = Green,
  urlcolor = Green
}

\DeclareMathOperator\qdim{\mathrm{dim}_q}
\DeclareMathOperator\defect{\mathrm{def}}

\usepackage{todonotes}

\usepackage{aliascnt}
\def\NewTheorem#1{%
    \newaliascnt{#1}{equation}
    \newtheorem{#1}[#1]{#1}
    \aliascntresetthe{#1}
    \expandafter\def\csname #1autorefname\endcsname{#1}
}
\def\equationautorefname~#1\null{(#1)\null}

\theoremstyle{plain}
\newcounter{mainTheorem}

\newtheorem{MainTheorem}[mainTheorem]{Theorem}

\newcounter{mainCorollary}
\numberwithin{mainCorollary}{mainTheorem}

\newtheorem{MainCorollary}[mainCorollary]{Corollary}

\swapnumbers
\numberwithin{equation}{section}
\NewTheorem{Assumption}
\NewTheorem{Proposition}
\NewTheorem{Theorem}
\NewTheorem{Corollary}
\NewTheorem{Conjecture}
\NewTheorem{Lemma}
\theoremstyle{definition}
\NewTheorem{Definition}
\theoremstyle{remark}
\NewTheorem{Remark}
\newtheorem*{Remark*}{Remark}
\NewTheorem{Remarks}

\theoremstyle{definition}
\NewTheorem{Example}
\AtEndEnvironment{Example}{\null\hfill$\Diamond$}%

\def\Sum{\displaystyle\sum}

\usepackage{xparse,mathtools}
\newcommand\SetBox[2][34]{\Big\{\vcenter{\hsize#1mm\centering#2}\Big\}}
\DeclarePairedDelimiterX{\set}[1]{\{}{\}}{\setargs{#1}}
\NewDocumentCommand{\setargs}{>{\SplitArgument{1}{|}}m}
{\setargsaux#1}
\NewDocumentCommand{\setargsaux}{mm}
{\IfNoValueTF{#2}{#1} {#1\,\delimsize|\,\mathopen{}#2}}

\def\map#1#2{\,{:}\,#1\!\longrightarrow\!#2}
\def\bijection{\overset{\simeq}{\longrightarrow}}

\def\({\big(}
\def\){\big)}

\let\<=\langle
\let\>=\rangle

\def\pmod#1{\text{ }(\text{mod } #1)\,}

\usepackage{tikz}
\usetikzlibrary{matrix}

\newcount\tableauRow\newcount\tableauCol
\newcommand\Tableau[2][-2]{
\begin{tikzpicture}[scale=0.4,draw/.append style={thick,black},baseline=#1mm]
  \tableauRow=0
  \foreach \Row in {#2} {
  \tableauCol=1
  \foreach\k in \Row {
  \draw(\the\tableauCol,\the\tableauRow)+(-.5,-.5)rectangle++(.5,.5);
  \draw(\the\tableauCol,\the\tableauRow)node{\k};
  \global\advance\tableauCol by 1
  }
  \global\advance\tableauRow by -1
  }
\end{tikzpicture}
}

\def\i{\hat\imath}
\def\j{\hat\jmath}

\def\bi{\mathbf{i}}
\def\bj{\mathbf{j}}
\def\bk{\mathbf{k}}

\let\eps=\varepsilon

\let\gedom=\unrhd
\let\notgedom\ntrianglerighteq
\let\gdom=\rhd

\newcommand{\J}{\mathcal{J}}
\newcommand{\K}{\mathcal{K}}
\renewcommand{\O}{\mathcal{O}}

\newcommand\RootMr[1][\pm]{\sqrt{\vrule height 2mm width 0pt\smash{M_r'}}}
\newcommand\RootLr[1][\pm]{\sqrt{\vrule height 2mm width 0pt\smash{L_r'}}}

\newcommand\dhash{\dot\#}

\newcommand{\Sn}[1][n]{\mathfrak{S}_{#1}}
\newcommand{\An}[1][n]{\mathfrak{A}_{#1}}

\renewcommand{\H}{\mathscr{H}}

\newcommand{\HO}[1][t]{\H^\O_{#1}}

\newcommand{\KHAn}{\H_t^{\K}(\An)}
\newcommand{\KHSn}{\H_t^{\K}(\Sn)}
\newcommand{\FHSn}[1][F]{\H_\xi^{#1}(\Sn)}
\newcommand{\FHAn}[1][F]{\H_\xi^{#1}(\An)}
\newcommand{\HSn}[1][n]{\H_\xi(\Sn[#1])}
\newcommand{\HAn}[1][n]{\H_\xi(\An[#1])}
\newcommand{\RSn}[1][n]{\mathscr{R}_e(\Sn[#1])}
\newcommand{\RpSn}[1][n]{\mathscr{R}^\varepsilon_{#1}}
\newcommand{\RpSnp}[1][+]{\mathscr{R}^{\varepsilon#1}_\gamma}
\newcommand{\Ralpha}[1][\alpha]{\RSn_{#1}}
\newcommand{\FRSn}[1][F]{\mathscr{R}^{#1}_e(\Sn)}
\NewDocumentCommand\RAn{D<>{e} O{n}}{\mathscr{R}_{#1}(\An[#2])}
\newcommand{\FRAn}[1][F]{\mathscr{R}^{#1}_e(\An)}

\newcommand{\N}{\mathbb{N}}
\newcommand{\Klesh}[1][n]{\mathcal{R}_{#1}}
\newcommand{\Parts}[1][n]{\mathcal{P}_{#1}}

\newcommand{\Z}{\mathbb{Z}}
\newcommand{\F}{\mathbb{F}}
\newcommand{\Lscr}{\mathscr{L}}
\newcommand{\Acal}{\mathcal{A}}

\newcommand{\Zcal}{\mathcal{Z}}
\def\Std{\mathop{\rm Std}\nolimits}

\newcommand{\la}{\lambda}

\newcommand\fo[1][\bi]{f^\O_{#1}}
\newcommand\yp{y^+}
\newcommand\ym{y^-}
\newcommand\yo{y^\O}
\newcommand\psio{\psi^\O}
\newcommand\psip{\psi^+}
\newcommand\psim{\psi^-}

\newcommand\dyo[1]{y^{\<#1\>}}

\newcommand{\wpsio}{{\dot\psi}^\O}
\newcommand{\wyo}{\dot{y}^\O}
\newcommand{\wfo}[1][\bi]{\dot{f}^\O_{#1}}
\newcommand\ddyo[1]{\dot{y}^{\<#1\>}}
\newcommand\RO[1][\O]{\dot R^{#1}_\gamma}

\DeclareMathOperator{\row}{row}
\DeclareMathOperator{\col}{col}
\DeclareMathOperator{\rank}{rank}
\DeclareMathOperator{\height}{ht}
\DeclareMathOperator{\Deg}{Deg}
\DeclareMathOperator{\codeg}{codeg}

\newcommand{\m}{\mathfrak{m}}
\newcommand{\mull}{\mathbf{m}}

\DeclareMathOperator\sh{\text{Shape}}

\renewcommand{\a}{\mathsf{a}}
\renewcommand{\b}{\mathsf{b}}
\newcommand{\s}{\mathsf{s}}
\renewcommand{\t}{\mathsf{t}}
\renewcommand{\u}{\mathsf{u}}
\renewcommand{\v}{\mathsf{v}}

\newcommand\balpha{{\boldsymbol{\alpha}}}
\newcommand{\ei}{{\mathbf{i}}}
\newcommand{\ej}{{\mathbf{j}}}

\newcommand{\sgn}{\mathtt{sgn}}
\DeclareMathOperator{\res}{res}

\DeclareMathOperator{\Res}{Res}

\DeclareMathOperator\noedge{\:\rlap{\hspace*{0.25em}/}\text{---}\:}

\newcommand\Qpm{Q^\varepsilon_n}

\makeindex
\def\Email#1{\email{\href{mailto:#1}{#1}}}
\begin{document}
\title{Quiver Hecke algebras for alternating groups}
\author{Clinton Boys}
\address{School of Mathematics and Statistics F07, University of
Sydney, NSW 2006, Australia.}
\Email{clinton.boys@sydney.edu.au}
\author{Andrew Mathas}
\Email{andrew.mathas@sydney.edu.au}

\begin{abstract}
  The main result of this paper shows that, over large enough fields of
  characteristic different from~$2$, the alternating Hecke algebras are
  $\Z$-graded algebras that are isomorphic to fixed-point subalgebras of
  the quiver Hecke algebra of~the symmetric group~$\Sn$. As a special
  case, this shows that the group algebra of the alternating group,
  over large enough fields of characteristic different from~$2$, is a
  $\Z$-graded algebra. We give a homogeneous presentation for these
  algebras, compute their graded dimension and show that the blocks of
  the quiver Hecke algebras of the alternating group are graded
  symmetric algebras.
\end{abstract}

\maketitle

\section*{Introduction}

  In a landmark paper, Brundan and Kleshchev \cite{BK:GradedKL}
  constructed an explicit $\Z$-grading on the cyclotomic Hecke algebras of
  type~$A$. These algebras include, as special cases, the group algebras
  of the symmetric group and the Iwahori-Hecke algebras of type~$A$. This
  paper extends these results to the group algebras of the alternating
  groups and, more generally, to Mitsuhashi's alternating Hecke
  algebras~\cite{Mitsuhashi:A}.

  Let $\HSn$ be the Iwahori-Hecke algebra of the symmetric group with
  parameter $\xi\in F^\times$, where $F$ is a field. Then $\HSn$ is a
  deformation of the group algebra of~$\Sn$. The algebra $\HSn$ has an
  automorphism $\#$ that can be considered as a $\xi$-deformation of the
  sign automorphism of~$F\Sn$. The alternating Hecke algebra $\HAn=\HSn^\#$
  is the fixed-point subalgebra of~$\HSn$ under~$\#$.

  Brundan and Kleshchev showed that $\HSn$ is a $\Z$-graded algebra by
  constructing an explicit family of isomorphisms $\theta:\RSn\bijection\HSn$,
  where $\RSn$ is a quiver Hecke algebra
  of~$\Sn$~\cite{BK:GradedKL,BrundanStroppel:KhovanovI,KhovLaud:diagI,Rouq:2KM}.
  Here, $e$ is \textbf{quantum characteristic} of~$\xi$, so $e>0$ is
  minimal such that $1+\xi+\dots+\xi^{e-1}=0$.

  As observed in \cite[(3.14)]{KMR:UniversalSpecht}, the algebra $\RSn$
  has a homogeneous automorphism $\sgn$ that is a graded analogue of the
  sign automorphism of the symmetric group. Let $\RAn=\RSn^\sgn$ be the
  fixed-point subalgebra of~$\RSn$ under~$\sgn$. Then $\RAn$ is a
  homogeneous subalgebra of~$\RSn$. It is natural to hope that $\theta$
  restricts to an isomorphism $\RAn\bijection\HAn$.  Unfortunately, the
  isomorphisms constructed by Brundan and Kleshchev do not restrict to
  isomorphisms between the alternating subalgebras; see
  \autoref{Ex:BadBKRestriction}.

  Let $F$ be a field and $\xi\in F$ an element of quantum
  characteristic~$e$. By definition, the field~$F$ is \textbf{large
  enough} for~$\xi$ if~$F$ contains square roots $\sqrt\xi$ and
  $\sqrt{1+\xi+\xi^2}$ whenever $e>3$.

  \begin{MainTheorem}\label{T:Main}
    Suppose that $\xi\in F$ is an element of quantum characteristic~$e\ne2$,
    where~$F$ is a large enough field for~$\xi$ of characteristic
    different from~$2$. Then $\HAn\cong\RAn$.
  \end{MainTheorem}

  To prove this result we construct a new isomorphism $\RSn\bijection\HSn$
  that intertwines the involutions,~$\sgn$ and~$\#$, on the two algebras. We
  do this using the framework developed by Hu and the second-named
  author~\cite{HuMathas:SeminormalQuiver}, which shows that the KLR grading
  can be described explicitly in terms of seminormal forms.

  As the algebra $\RAn$ is a graded subalgebra of $\RSn$ we immediately
  obtain the following.

  \begin{MainCorollary}\label{C:HAnZGraded}
    Suppose that $\xi\in F$ is an element of quantum
    characteristic~$e\ne2$, where~$F$ is a large enough field
    for~$\xi$ of characteristic different from~$2$. Then $\HAn$ is a
    $\Z$-graded algebra.
  \end{MainCorollary}

  In particular, over large enough fields of characteristic different
  from~$2$ the group algebra $F\An$ of the alternating group is
  $\Z$-graded, for $n\ge1$. The alternating group corresponds to the
  case when~$\xi=1$, so if~$F$ is a field of characteristic~$p\ne2$ then
  $F\An\cong \RAn<p>$ if~$3$ has a square root in $F$ whenever $p>3$.

  Applying \autoref{T:Main} twice shows that, up to isomorphism, $\HAn$
  depends only on~$e$, the quantum characteristic of~$\xi$, rather than
  on~$\xi$ itself. Hence, the following holds:

  \begin{MainCorollary}
    Let $F$ be a field of characteristic different from~$2$. Suppose
    that $\xi,\xi'\in F$ are elements of quantum characteristic~$e\ne2$,
    where~$F$ is a large enough field for~$\xi$ and for $\xi'$. Then
    $\HAn\cong\H_{\xi'}(\An)$.
  \end{MainCorollary}

  In particular, over a large enough field, the decomposition matrix
  of~$\HAn$ depends only on~$e$, and the field~$F$, and not on the choice
  of~$\xi$.

  The quiver Hecke algebra $\RSn$ has a homogeneous presentation by
  generators and relations that is described in terms of the quiver
  $\Gamma_e$ with vertex set $I=\Z/e\Z$ and edges $i\to i+1$, for $i\in I$.
  In \autoref{S:KLRAlgebras} we associate to~$\Gamma_e$ a set of
  simple roots
  $\set{\alpha_i|i\in I}$ and a positive root lattice
  $Q^+=\bigoplus_{i\in I}\N\alpha_i$. If $\alpha\in Q^+$ let
  $I^\alpha=\set{\bi\in I^n|\alpha=\alpha_{i_1}+\dots+\alpha_{i_n}}$.
  Let~$\sim$ be the equivalence relation on~$I^n$ generated by $\bi\sim\bj$ if $\bj=-\bi$.
  This relation induces an equivalence relation on~$Q^+$ where
  $\alpha\sim\beta$ if there exists $\bi\in I^\alpha$ such that
  $-\bi\in I^\beta$.
  Let $\Qpm=Q^+_n/{\sim}$. If $\gamma\in\Qpm$ let
  $I^\gamma=\bigcup_{\alpha\in\gamma}I^\alpha$.

  Using the KLR presentation of~$\RSn$, and the realisation of~$\RAn$ as
  a fixed-point subalgebra of~$\RSn$, gives the following homogeneous
  presentation for~$\RAn$ with respect to the $\Z$-grading. For
  $\bi,\bj\in I^n$ set $\delta_{\bi\sim\bj}=1$ if $\bi\sim\bj$ and set
  $\delta_{\bi\sim\bj}=0$ otherwise.

  \begin{MainTheorem}\label{T:MainRelations}
    Suppose that $e\ne2$ and that $2$ is invertible in~$\Zcal$. Then
    \[\RAn=\bigoplus_{\gamma\in \Qpm}\RAn_\gamma \]
    where $\RAn_\gamma$ is the unital associative $\Zcal$-algebra generated by
    elements
    \[\set{\Psi_{r}(\bi), Y_{s}(\bi), \eps(\bi) |
    1\le r\le n, 1\le s<n\text{ and }\bi\in I^\gamma}
    \]
    subject to the relations
    {\setlength{\abovedisplayskip}{2pt}
    \setlength{\belowdisplayskip}{1pt}
    \begin{align*}
      Y_{1}(\bi)^{(\Lambda_0,\alpha_{i_1})}&=0,
      &\eps(\bi)\eps(\bj)&= \delta_{\bi\sim\bj}\eps(\bi),
      &\textstyle\sum_{\bi\in I^\gamma}\frac12\eps(\bi)&= 1,\\
      Y_r(-\bi)&=-Y_r(\bi)
      &\Psi_r(-\bi)&=-\Psi_r(\bi)
      &\eps(-\bi)&=\eps(\bi)\\
      \eps(\bi)Y_r(\bj)\eps(\bk)&=\delta_{\bi\bj}\delta_{\bj\bk} Y_r(\bj),
      &\eps(\bi)\Psi_r(\bj)\eps(\bk)&
      =\delta_{s_r\cdot\bi,\bj}\delta_{\bj\bk}\Psi_r(\bj),
      &Y_r(\bi)Y_s(\bj)&=\delta_{\bi\bj}Y_s(\bi)Y_r(\bj),
    \end{align*}
    \begin{align*}
      \Psi_r(\bi)Y_{r+1}(\bi) & =\(Y_r(s_r\cdot\bi)\Psi_r(\bi)
      +\delta_{i_r i_{r+1}}\eps(\bi)\),\\
      Y_{r+1}(s_r\cdot\bi)\Psi_r(\bi) &=\(\Psi_r(\bi)Y_r(\bi)
      +\delta_{i_r i_{r+1}}\eps(\bi)\),
    \end{align*}
    \begin{align*}
      \Psi_r(\bi)Y_s(\bi) &=Y_s(s_r\cdot\bi)\Psi_r(\bi),
      &&\text{if }s\neq r,r+1,\\
      \Psi_r(s_t\cdot\bi)\Psi_t(\bi)
      &=\Psi_t(s_r\cdot\bi)\Psi_r(\bi),&&\text{if }|r-s|>1,
    \end{align*}
    \begin{align*}
      \Psi_r(s_r\cdot\bi)\Psi_r(\bi)&= \begin{cases}
        Y_r(\bi)-Y_{r+1}(\bi),&\text{if }i_r\to i_{r+1},\\
        Y_{r+1}(\bi)-Y_r(\bi),&\text{if }i_r\leftarrow i_{r+1},\\
        0,&\text{if }i_r=i_{r+1},\\
        \eps(\bi),&\text{otherwise}
      \end{cases}
    \end{align*}
    and $\Psi_r(s_{r+1}s_r\cdot\bi)\Psi_{r+1}(s_r\cdot\bi)\Psi_r(\bi)
    -\Psi_{r+1}(s_rs_{r+1}\cdot\bi)\Psi_r(s_{r+1}\cdot\bi)\Psi_{r+1}(\bi)$
    is equal to
    \begin{align*}
      \begin{cases}
        \eps(\bi), &\text{if }i_r=i_{r+2}\leftarrow i_{r+1},\\
        -\eps(\bi), &\text{if }i_r=i_{r+2}\to i_{r+1},\\
        0,&\text{otherwise,}
      \end{cases}
    \end{align*}
    }%
    for all $\bi,\bj,\bk\in I^\gamma$ and all admissible $r,s$ and $t$.
  \end{MainTheorem}

  The fraction~$\frac12$ appears in the relation $\sum_{\bi\in
  I^\gamma}\frac12\eps(\bi)=1$ because $\eps(\bi)=\eps(-\bi)$.
  (It follows from the relations in \autoref{T:MainRelations} that
  $\eps(\bi)=0$ if $\bi=-\bi$, for $\bi\in I^n$.)

  The relations in \autoref{T:MainRelations} are homogeneous with respect
  to the degree function
  \[
  \deg \eps(\bi)=0, \quad \deg Y_r(\bi)=2\quad\text{and}\quad
  \deg\Psi_r(\bi)=-(\alpha_{i_r},\alpha_{i_{r+1}}),
  \]
  where $(\ ,\ )$ is the Cartan pairing. Hence, $\RAn$ is a $\Z$-graded
  algebra. Quite surprisingly, this presentation is almost identical to
  the KLR-presentation of~$\RSn$. The main difference being the three
  ``sign relations'' relating the generators indexed by $\bi$ and $-\bi$, for
  $\bi\in I^\gamma$. The key idea behind the proof
  \autoref{T:MainRelations} is to introduce a $(\Z_2\times\Z)$-grading on
  the KLR algebra~$\RSn$. With respect to the $(\Z_2\times\Z)$-grading,
  the algebra $\RAn$ is the even part of~$\RSn$. Using this observation
  we can deduce the relations for~$\RAn$ directly from those for~$\RSn$.

  Finally, using the graded cellular bases of Hu and
  Mathas~\cite{HuMathas:GradedCellular} we construct a homogeneous basis
  for the $\Z$-graded algebra~$\RAn$. As a corollary we obtain the
  graded dimension of~$\RAn$. See \autoref{S:Basis} below for the
  unexplained notation.

  \begin{MainTheorem}\label{T:GDim}
    Suppose that $e\ne2$ and that $2$ is invertible in~$F$. Then the
    alternating quiver Hecke algebra $\RAn$ has graded dimension
    \[
    \sum_{\substack{(\s,\t)\in\Std^2(\Parts)\\\res(\s)\in I^n_+}}
    q^{\deg\s+\deg\t}.
    \]
  \end{MainTheorem}

  In fact, using the work of Li~\cite{GeLi:IntegralKLR} it follows that
  over any ring in which $2$ is invertible the algebra $\RAn$ is free
  with the same graded rank. As a second application of our basis
  theorem we show that the blocks of~$\RAn$ are graded symmetric
  algebras.

  The results in this paper exclude the cases when~$F$ is a field of
  characteristic~$2$ and when $e=2$ or, equivalently, $\xi=-1$.  This is
  because most our arguments fail, and most of our results are false,
  when we drop these assumptions.

  The paper is organised as follows. \autoref{S:HeckeAlgebras} starts by
  defining the quiver Hecke algebra~$\RAn$ of the alternating group as
  the fixed-point subalgebra of~$\RSn$ under the homogeneous sign
  involution~$\sgn$. We then prove \autoref{T:MainRelations} by first
  introducing a $(\Z_2\times\Z)$-grading on~$\RSn$ and showing
  that~$\RAn$ is the \textit{even} part of~$\RSn$, with respect to the
  $\Z_2$-grading.  \autoref{S:SeminormalForm} starts by setting up the
  framework of \textit{seminormal coefficient systems} and showing how
  seminormal bases behave under the ungraded sign involution~$\#$.
  Building on ideas from~\cite{HuMathas:SeminormalQuiver}, we give a new
  presentation of~$\HO$, over a specially chosen ring~$\O$, that we use
  to construct a new isomorphism $\Theta:\RSn\bijection\HSn$ over the
  residue field of~$\O$. Unlike the known isomorphisms in the
  literature, $\Theta$ intertwines the two sign
  involutions,~$\sgn$ and~$\#$, implying \autoref{T:Main}. In
  \autoref{S:Basis} we give a homogeneous basis of~$\RAn$ and hence
  prove \autoref{T:GDim}. As an application we show that the blocks
  of~$\RAn$ are graded symmetric algebras. Finally, using Clifford
  theory, the classification of the blocks and irreducible graded
  modules of~$\RAn$ is given.

  \subsection*{Acknowledgments}
  We thank the referee for their careful reading of our manuscript.
  The first-named author was supported by an Australian Postgraduate Award
  and the second-named author by the Australian Research Council. Some of
  the work in the first-named author's PhD thesis \cite{Boys:PhDThesis}
  was based on an earlier version of this paper.

\section{Iwahori-Hecke algebras and quiver Hecke
algebras}\label{S:HeckeAlgebras}
This chapter defines both the alternating Hecke algebras and the
alternating quiver Hecke algebras of type~$A$. Both algebras are defined as fixed
point subalgebras of the corresponding Hecke algebras. In the
final section we prove that \autoref{T:MainRelations} gives a homogeneous
presentation for the alternating quiver Hecke algebra.

\subsection{Iwahori-Hecke algebras and alternating Hecke algebras}
We start by defining the Iwahori-Hecke algebras of the symmetric groups.
These algebras are well-studied deformations of the group algebras of
the symmetric groups that arise naturally in the representation theory
of the general linear groups.

Fix a (unital) integral domain $\Zcal$ and an invertible element $\xi\in \Zcal ^\times$.

\begin{Definition}\label{D:Hecke}
  The \textbf{Iwahori-Hecke algebra} $\HSn=\H_\xi^\Zcal(\Sn)$ is the (unital)
  associative $\Zcal $-algebra with generators $T_1,\ldots,T_{n-1}$ subject to
  relations
  \begin{align*}
    (T_r-\xi)(T_r+1)&= 0,&\quad\text{for }r=1,\ldots,n-1,\\
    T_rT_s&= T_sT_r,&\quad\text{if }|r-s|>1,\\
    T_rT_{r+1}T_r&= T_{r+1}T_rT_{r+1},&\quad\text{for }r=1,\ldots,n-2.
  \end{align*}
\end{Definition}

For $1\le r<n$ let $s_r=(r,r+1)\in\Sn$. Then $\{s_1,\dots,s_{n-1}\}$ is
the standard set of Coxeter generators for~$\Sn$.  If $w\in\Sn$ then the
\textbf{length} of $w$ is the integer
$\ell(w)=\min\set{l\ge0|w=s_{r_1}\dots s_{r_l}\text{ with }1\le r_j<n}$.  A
\textbf{reduced expression} for~$w$ is any word $w=s_{r_1}\dots
s_{r_\ell}$ with $\ell=\ell(w)$ and $1\le r_j<n$, for $1\le j\le \ell$.

If $w\in\Sn$ define $T_w=T_{r_1}\dots T_{r_\ell}$, where
$w=s_{r_1}\dots s_{r_\ell}$ is any reduced expression.  As is
well-known, because the braid relations hold in~$\HSn$ the element~$T_w$
depends only on~$w$ and not on the choice of reduced expression.
Moreover, the algebra $\HSn$ is free as a $\Zcal$-module with basis
$\set{T_w|w\in\Sn}$. See, for example, \cite[Chapter~1]{M:ULect}.
In particular, if $\xi=1$ then $\HSn\cong\Zcal\Sn$ via the map
$T_w\mapsto w$, for $w\in\Sn$.

Following Goldman~\cite[Theorem~5.4]{Iwahori:Hecke},
let  $\#\map\HSn\HSn$ be the
unique $\Zcal $-linear automorphism of $\H_\xi(\Sn)$ such that
\begin{equation}\label{E:hash}
  T_r^\#=-\xi T_r^{-1}=-T_r+(\xi-1),
\end{equation}
for $1\le r<n$.  It follows directly from the definitions that when
$\xi=1$ the automorphism $\#$ of $\H^{\mathcal{Z}}_1(\Sn)\cong \Zcal
\Sn$ is the usual ``sign involution'' which sends each simple
transposition $s_r$ to  $-s_r$, for $1\le r<n$~\cite[p5]{James}.  Since
the group algebra of the alternating group is the fixed-point subalgebra
of the sign automorphism, the following definition gives a $\xi$-analogue
of the group ring $\Zcal\An$ of the alternating group~$\An$.

\begin{Definition}[Mitsuhashi~\cite{Mitsuhashi:A}]\label{D:AltHecke}
  Suppose that $\xi\ne-1$. Then the \textbf{ alternating Hecke algebra} is the
  fixed-point subalgebra
  \[
  \HAn=\H_\xi^\Zcal(\An)=\set{h\in\HSn\mid h^\#=h}
  \]
  of~$\HSn$ under the hash involution.
\end{Definition}

\begin{Remark}
  Mitsuhashi's~\cite[Definition 4.1]{Mitsuhashi:A} original definition of
  the alternating Hecke algebra was by generators and relations, giving a
  deformation of a well-known presentation of the alternating group.
  \autoref{D:AltHecke} is equivalent to Mitsuhashi's definition by
  \cite[Proposition~1.5]{MathasRatliff}.
\end{Remark}

\subsection{Graded modules and algebras}
The main result of this paper shows that the alternating Hecke algebra
is a graded algebra, so we quickly review this terminology. For the
most part we will work with $\Z$-graded modules and algebras, however,
to prove \autoref{T:MainRelations} we consider more general
gradings.

Recall that $\Zcal$ is a unital integral domain. In this paper all
modules will be assumed to be free and of finite rank as
$\Zcal$-modules.

Let $(G,+)$ be an abelian group. A \textbf{$G$-graded $\Zcal$-module} is
a $\Zcal$-module $M$ that admits a vector space decomposition
$M=\bigoplus_{g\in G} M_g$. If $g\in G$ and $0\ne m\in M_g$ then $m$ is
\textbf{homogeneous of degree}~$g$.  Similarly, a
\textbf{$G$-graded $\Zcal$-algebra} is a $G$-graded $\Zcal$-module
$A=\bigoplus_{g\in G}A_g$ that is a $\Zcal$-algebra such that
$A_fA_g\subseteq A_{f+g}$, for all $f,g\in G$. A  \textbf{$G$-graded $A$-module} is
a $G$-graded $\Zcal$-module $M$ such that $M_fA_g\subset M_{f+g}$, for
$f,g\in G$.

Unless otherwise specified, $G=\Z$ and a \textbf{graded module} will mean a
$\Z$-graded module. Similarly, a \textbf{graded algebra} is a $\Z$-graded
algebra.

\subsection{Cyclotomic quiver Hecke algebras of type~$A$}\label{S:KLRAlgebras}
We now define the second class of algebras that we are interested in:
the cyclotomic quiver Hecke algebras of~$\Sn$. These algebras are
certain quotients of the $\Z$-graded quiver Hecke algebras introduced,
independently, by Khovanov and Lauda~\cite{KhovLaud:diagI} and
Rouquier~\cite{Rouq:2KM}.

Fix  $e\in\set{3,4,5,\ldots}\cup\set{\infty}$ and define $\Gamma_e$ to
be the quiver with vertex set $I=\Z/e\Z$ and edges $i\to i+1$, for $i\in
I$. (By convention, $I=\Z$ if $e=\infty$.) Thus, $\Gamma_e$ is the
infinite quiver of type $A_\infty$ if $e=\infty$ and the finite quiver
of Dynkin type $A_{e-1}^{(1)}$ if $e\geq 3$.  We exclude $e=2$ only because
this corresponds to the case $\xi=-1$ in \autoref{D:AltHecke}, which we
do not consider in this paper.

Following Kac~\cite{Kac}, to the quiver $\Gamma_e$ we attach the usual
Lie theoretic data of the positive roots $\set{\alpha_i|i\in I}$, the
fundamental weights $\set{\Lambda_i|i\in I}$,
the positive weight lattice $P^+=\bigoplus_{i\in I}\N\Lambda_i$,
the positive root lattice $Q^+=\bigoplus_{i\in I}\N\alpha_i$,
the non-degenerate
pairing $(\  ,\ )\map{P^+\times Q^+}\Z$ given by
$(\Lambda_i,\alpha_j)=\delta_{ij}$, for $i,j\in I$, and the
\textbf{Cartan matrix} $C=(c_{ij})_{i,j\in I}$ where
\[
c_{ij}=\begin{cases}
  2,&\text{if }i=j,\\
  -1,&\text{if $i\leftarrow j$ or $i\rightarrow j$},\\
  0,&\text{otherwise.}
\end{cases}
\]
The \textbf{height} of $\alpha\in Q^+$ is the
non-negative integer $\height\alpha=\sum_i(\Lambda_i,\alpha)$. Fix $n\ge0$ and let
$Q^+_n=\set{\alpha\in Q^+|\height\alpha=n}$. For $\alpha\in Q^+_n$, let
\[
  I^\alpha=\set{\bi=(i_1,\dots,i_n)\in I^n|\alpha=\alpha_{i_1}+\dots+\alpha_{i_n}}.
\]

\begin{Definition}[Khovanov and Lauda \cite{KhovLaud:diagI}
  and Rouquier~\cite{Rouq:2KM}]
  \label{D:klr}
  Suppose that  $\alpha\in Q^+$ and
  $e\in\set{3,4,5,\ldots}\cup\set{\infty}$. The \textbf{cyclotomic quiver
  Hecke algebra} $\Ralpha$ is the unital associative
  $\Zcal$-algebra  with generators
  \begin{equation*}
    \set{\psi_1,\ldots,\psi_{n-1}}\cup \set{y_1,\ldots,y_n}
    \cup \set{e(\ei)| \ei\in I^\alpha}
  \end{equation*}
  and relations
  {\setlength{\abovedisplayskip}{2pt}
  \setlength{\belowdisplayskip}{1pt}
  \begin{xalignat*}{3}
    y_1^{(\Lambda_0,\alpha_{i_1})}e(\bi)&= 0,& e(\ei)e(\ej)&= \delta_{\ei\ej}e(\ei),
    &\textstyle\sum_{\ei\in I^\alpha}e(\ei)&= 1,\\
    y_re(\ei)&= e(\ei)y_r,& \psi_re(\ei)&= e(s_r\cdot\ei)\psi_r,& y_ry_s&= y_sy_r,
  \end{xalignat*}
  \begin{xalignat*}{2}
    \psi_ry_{r+1}e(\ei)&= (y_r\psi_r+\delta_{i_r i_{r+1}})e(\ei),
    & y_{r+1}\psi_re(\ei)&= (\psi_ry_r+\delta_{i_r i_{r+1}})e(\ei),
  \end{xalignat*}
  \begin{align*}
    \psi_ry_s&= y_s\psi_r,&\text{if }s\neq r,r+1,\\
    \psi_r\psi_s&= \psi_s\psi_r,&\text{if }|r-s|>1,
  \end{align*}
  \begin{align*}
    \psi_r^2e(\ei)&= \begin{cases}
      0,&\text{if }i_r=i_{r+1},\\
      (y_r-y_{r+1})e(\ei),&\text{if }i_r\to i_{r+1},\\
      (y_{r+1}-y_r)e(\ei),&\text{if }i_r\leftarrow i_{r+1},\\
      e(\ei),&\text{otherwise},
    \end{cases}\\
    \psi_r\psi_{r+1}\psi_re(\ei)&=\begin{cases}
      (\psi_{r+1}\psi_r\psi_{r+1}-1)e(\ei),&\text{if }i_r=i_{r+2}\to i_{r+1},\\
      (\psi_{r+1}\psi_r\psi_{r+1}+1)e(\ei),
      &\text{if }i_r=i_{r+2}\leftarrow i_{r+1},\\
      \psi_{r+1}\psi_r\psi_{r+1}e(\ei),&\text{otherwise,}
    \end{cases}
  \end{align*}
  }%
  for $\ei,\ej\in I^\alpha$ and all admissible $r$ and $s$. If $n\ge0$
  then the \textbf{quiver Hecke algebra of~$\Sn$}
  is the algebra
  \[
  \RSn =\bigoplus_{\alpha\in Q^+_n}\Ralpha.
  \]
\end{Definition}

Note that the algebra $\Ralpha$ depends on $e$, $\Gamma_e$, $\Lambda_0$
and $\alpha\in Q^+$.

We write $\RSn=\FRSn[\Zcal]$ when we want to emphasise that $\RSn$ is a
$\Zcal$-algebra. The main advantage of the relations in \autoref{D:klr} is
that they are homogeneous with respect
to the following $\Z$-valued degree function:
\begin{align}
  \deg e(\ei)&= 0, &\text{for all }\ei\in I^n,\nonumber\\
  \deg y_r&= 2,&\text{for }1\le r\le n,\label{degfunc}\\
  \deg \psi_re(\ei)&= -c_{i_r,i_{r+1}},&\text{for $1\le r<n$ and $\ei\in I^n$}.
  \nonumber
\end{align}
Therefore, $\RSn$ is a $\Z$-graded algebra.

\begin{Remark}
  There are fewer relations appearing in \autoref{D:klr} than in
  \cite[Theorem~1.1]{BK:GradedKL}. This is because we are assuming that
  $e\ne 2$ (and $\Lambda=\Lambda_0$). We have also made a sign change
  compared with \cite{BK:GradedKL}, which is consistent with
  \cite{HuMathas:SeminormalQuiver}.
\end{Remark}

In examples, we write $e(\ei)=e(i_1i_2\dots i_n)$ if
$\ei=(i_1,i_2,\dots,i_n)$.

\begin{Example}\label{Ex:OS3}
  Let $n=3$, $e=3$ and $\Lambda=\Lambda_0$. First, $y_1=0$ because of the
  relations $y_1^{(\Lambda_0,\alpha_{i_1})}e(\ei)=0$ and $\sum_\bi e(\bi)=0$.
  It is not difficult to see that $\psi_1=0=y_2$ and
  that $e(\ei)=0$ unless $\ei=(012)$ or $(021)$; see, for example
  \cite[Proposition~2.4.6]{Mathas:Singapore}.
  Hence, $\RSn[3]$ is generated by $\psi_2,y_3,e(012)$ and $e(021)$. Using
  the quadratic relation,
  \[
  \psi_2^2e(\ei)=\begin{cases}
    -y_3e(\ei),&\text{if }\ei=(012),\\
    \phantom{-}y_3e(\ei),&\text{if }\ei=(021).
  \end{cases}
  \]
  In turn, this implies that $y_3^2e(\ei)=\pm y_3\psi_2e(\ei)=\pm \psi_2y_2e(\ei)=0$,
  so $y_3^2=0$. Therefore, $\RSn[3]$ is spanned by
  \[
  \set{e(012), e(021), \psi_2e(012), \psi_2e(021), y_3e(012), y_3e(021)}.
  \]
  By \autoref{T:BKiso} below, these elements are a basis of~$\RSn$.
\end{Example}

To connect the algebras $\RSn$ and $\HSn$ define the \textbf{quantum
characteristic} of $\xi$ to be the smallest non-negative integer $e$
such that $1+\xi+\dots+\xi^{e-1}=0,$ and set $e=\infty$ if no such
integer exists.  By definition,
$e\in\set{2,3,4,5,6,\dots}\cup\set{\infty}$ and $e=2$ if and only if
$\xi=-1$.

Brundan and Kleshchev proved the following remarkable theorem, which
connects the two algebras $\RSn$ and $\HSn$.

\begin{Theorem}[\protect{Brundan and Kleshchev~\cite{BK:GradedKL},
  Rouquier~\cite[Corollary~3.20]{Rouq:2KM}}]\label{T:BKiso}
  \leavevmode\newline
  Suppose that $F$ is a field and that $\xi\ne-1$ has quantum
  characteristic~$e>2$. Then $\FRSn[F]\cong\FHSn[F]$.
\end{Theorem}

Hence, if $F$ is a field then we can consider $\FHSn$ as a $\Z$-graded
algebra via the isomorphism $\FHSn\cong\FRSn$. We have stated a special
case of Brundan and Kleshchev's result because this is all that we need.
Over a field, Brundan and Kleshchev prove more generally that the cyclotomic Hecke
algebras of type $A$ are isomorphic to cyclotomic quiver Hecke algebras
--- and they also allow $e=2$. To prove \autoref{T:BKiso}
Brundan and Kleshchev construct an explicit isomorphism (in fact, they
construct different isomorphisms for the cases when $\xi=1$ and
$\xi\ne1$). Our proof of \autoref{T:MainRelations} builds from a
variation on their ideas, following~\cite{HuMathas:SeminormalQuiver}.

\subsection{Alternating quiver Hecke algebras of type~$A$}
In this section we introduce a homogeneous analogue of the
$\#$-involution of~$\HSn$ and use it to define the alternating quiver
Hecke algebras of type~$A$.

If $\bi=(i_1,\dots,i_n)\in I^n$ let $-\bi=(-i_1,\dots,-i_n)\in I^n$.
Following~\cite[(3.14)]{KMR:UniversalSpecht}, define $\sgn$ to be the unique
automorphism of $\RSn$ such that
\[
\psi_r^\sgn= -\psi_r,\quad y_s^\sgn= -y_s,\quad\text{and}\quad e(\ei)^\sgn= e(-\ei),
\]
for $1\le r<n$, $1\le s\le n$ and $\bi\in I^n$.

If $\alpha\in Q^+_n$ let $\alpha'$ be the unique element of~$Q^+_n$
such that $(\Lambda_i,\alpha)=(\Lambda_{-i},\alpha')$, for all $i\in I$.
Recall from the introduction that $\sim$ is the equivalence relation
on~$I^n$ generated by $\bi\sim\bj$ if $\bj=-\bi$ and that if
$\alpha,\beta\in Q^+$ then $\alpha\sim\beta$ if there exists $\bi\in
I^\alpha$ such that $-\bi\in I^\beta$. Hence, $\alpha\sim\beta$ if and
only if $\beta\in\set{\alpha,\alpha'}$.

Checking the relations in \autoref{D:klr} reveals the following.

\begin{Proposition}[\protect{\cite[(3.14)]{KMR:UniversalSpecht}}]
  \label{P:SgnHomogeneous}
  The map $\sgn$ restricts to a homogeneous isomorphism of $\Z$-graded
  algebras
  $\Ralpha\bijection\Ralpha[\alpha']$, for $\alpha\in Q^+_n$.
  Hence, $\sgn$ is a homogeneous automorphism of~$\RSn$ of order~$2$.
\end{Proposition}

Mirroring \autoref{D:AltHecke}, we define the second algebra
appearing in \autoref{T:Main}.

\begin{Definition}\label{D:sgndef}
  The \textbf{alternating quiver Hecke algebra of~$\An$} is the fixed-point
  subalgebra $\RAn=\mathscr{R}_e^\Zcal(\An)=\set{a\in\RSn|a^\sgn=a}$ of~$\RSn$
  under the involution~$\sgn$.
\end{Definition}

Since $\sgn$ is a homogeneous involution of~$\RSn$, an immediate and important
consequence of \autoref{D:sgndef} is the following.

\begin{Corollary}\label{P:GradedSubalgebra}
  The alternating quiver Hecke algebra $\RAn$ is a $\Z$-graded subalgebra
  of~$\RSn$.
\end{Corollary}

We finish this section with an example of how our main result,
\autoref{T:Main}, works when $n=3=e$.

\begin{Example}\label{Ex:A3Basis}
  Suppose that $e=3$ and $n=3$. Then $\An[3]\cong\Z/3\Z$ is the cyclic group of
  order~$3$. By \autoref{Ex:OS3}, $\mathscr{R}_3(\An[3])$ is spanned
  by the three elements
  \[
  1=e(012)+e(021),\quad \Psi=\psi_2(e(012)-e(021)),\quad Y=y_3(e(012)-e(021)).
  \]
  These three elements are homogeneous, with $\deg\Psi=1$ and $\deg
  Y=2$, and \autoref{T:BKiso} implies that they are non-zero. Therefore,
  $\{1,\Psi, Y\}$ is a basis of~$\mathscr{R}_3(\An[3])$.  Using the
  relations in \autoref{D:klr}, $\Psi^2=-Y$ and $\Psi^3=-Y\Psi=0$.
  Therefore, the map $\Psi\mapsto x$ determines an isomorphism of graded
  algebras
  \[\RAn[3]=\<\Psi\mid \Psi^3=0\>\cong\Z[x]/(x^3),\]
  where we put $\deg x=1$. It is now easy to see that
  $\mathscr{R}_3(\An[3])\otimes \F_3\cong\F_3\An[3]$, where
  $\F_3=\Z/3\Z$. For example, an isomorphism is determined
  by~$\Psi\mapsto 1-s_1s_2$.

  The isomorphism above is not unique. In the special case when $n=3$
  and $\xi=1\in\F_3$, the proof of \autoref{T:Main} in \autoref{S:Main}
  constructs a different isomorphism
  $\F_3\An[3]\cong\HAn[3]\bijection\RAn[3]$ that is determined by
  $\Psi\mapsto s_2s_1-s_1s_2$.
\end{Example}

The algebra $\RSn$ is defined in terms of the subalgebras $\Ralpha$. To
give a presentation for~$\RAn$ we need to work with the blocks of these
algebras. As in the introduction, set $\Qpm=Q^+_n/{\sim}$. Using the
notation introduced before \autoref{P:SgnHomogeneous}, if $\alpha\in
Q^+_n$ then $\set{\alpha,\alpha'}$ is its $\sim$-equivalence class. If
$\gamma\in\Qpm$ set~$\Ralpha[\gamma]=\bigoplus_{\alpha\in\gamma}\Ralpha$.
By \autoref{P:SgnHomogeneous}, $\sgn$ induces a homogeneous automorphism
of~$\Ralpha[\gamma]$. Define
\begin{equation}\label{E:RAngamma}
  \RAn_\gamma=(\Ralpha[\gamma])^\sgn=\set{a\in\Ralpha[\gamma]|a=a^\sgn}
\end{equation}
to be the fixed-point subalgebra of $\Ralpha[\gamma]$ under the
$\sgn$ automorphism. Since $\RSn=\bigoplus_\alpha\Ralpha$ we have the
following decomposition of~$\RAn$ as a direct sum of two-sided
$\Z$-graded ideals.

\begin{Corollary}\label{C:RAnBlocks}
  As a graded algebra,
  $\displaystyle\RAn=\bigoplus_{\gamma\in\Qpm}\RAn_\gamma$.
\end{Corollary}

\subsection{A presentation for alternating quiver Hecke algebras of $\An$}
\label{S:GeneratorsRelations}

In this section we prove \autoref{T:MainRelations}. To do this we  first give
a ``super'' presentation for the quiver Hecke algebra~$\RSn$.
For convenience, we identify the group $\Z_2=\Z/2\Z$ with $\{0,1\}$ in
the obvious way.

We start by defining a new algebra $\RpSn$ that, it turns out, is
isomorphic to~$\RSn$. The advantage of $\RpSn$ is that it is
$(\Z_2\times\Z)$-graded where the $\Z_2$-grading encodes the effects of~$\sgn$.
Abusing notation, we use similar notation for the generators of~$\RSn$
and~$\RpSn$. This is justified by \autoref{P:RpSnIsomorphism} below.

The sequence $\bi=(0,\dots,0)\in I^n$, which is the unique sequence such that
$\bi=-\bi$, is potentially problematic for us. The next result
resolves this.

\begin{Lemma}\label{L:ZeroSequence}
  Suppose that $\bi\in I^n$ and $e(\bi)\ne0$. Then $i_1=0$ and
  $i_2=\pm1$. In particular, $e(0,\dots,0)=0$.
\end{Lemma}

\begin{proof}By \cite[Lemma 4.1c]{HuMathas:GradedCellular},
  $e(\bi)\ne0$ if and only $\bi$ is the residue sequence of some
  standard tableau (see \autoref{S:tableaux}), which readily implies
  the result. As we need this argument later, we give a direct proof
  following \cite[Proposition~2.4.6]{Mathas:Singapore}.  First, $y_1=0$ and
  $e(\bi)\ne0$ only if $i_1=0$ by the cyclotomic relation
  $y_1^{(\Lambda_0,\alpha_{i_1})}e(\bi)=0$. If $i_1=i_2=0$ then
  $e(\bi)=(y_2\psi_1-\psi_1y_1)e(\bi)=y_2\psi_1e(\bi)=y_2e(\bi)\psi_1$,
  so that $e(\bi)=y_2^2e(\bi)\psi_1^2=0$. Hence,
  $e(0,0,i_3,\dots,i_n)=0$. Finally, suppose that $i_2\ne\pm1,0$. Then
  $e(\bi)=\psi_1^2e(\bi)=\psi_1e(i_2,i_1,i_3\dots,i_n)\psi_1=0$ since
  $e(\bj)=0$ whenever $j_1\ne0$.
\end{proof}

Set $I^n_+=\set{\bi\in I^n|i_1=0\text{ and }i_2=+1}$ and
$I^n_-=\set{\bi\in I^n|i_1=0\text{ and }i_2=-1}$. Then
\autoref{L:ZeroSequence}
shows that $e(\bi)\ne0$ only if $\bi\in I^n_+\cup I^n_-$.

Recall that if $\gamma\in\Qpm$ then $I^\gamma=\bigcup_{\alpha\in\gamma}I^\alpha$.

\begin{Definition}\label{D:RpSn}
  Suppose that $e\ne2$, $2$ is invertible in $\Zcal$ and that
  $\gamma\in\Qpm$. The algebra
  $\RpSn[\gamma]=\RpSn[\gamma](\Gamma_e,\Lambda_0)$ is the unital
  associative $\Zcal$-algebra with generators
  \[\set{\psi_r, y_s, \eps_a(\bi)|1\le r<n, 1\le s\le n,
  \bi\in I^\gamma\text{ and } a\in\Z_2}
  \]
  subject to the relations
  {\setlength{\abovedisplayskip}{2pt}
  \setlength{\belowdisplayskip}{1pt}
  \begin{align*}
    y_1^{(\Lambda_0,\alpha_{i_1})}\eps_0(\bi)&= 0,
    & \textstyle\sum_{\bi\in I^\gamma}\frac12\eps_0(\bi)&= 1,
    & \eps_0(\bi)\eps_0(\bj)&= \delta_{\bi\sim\bj}\eps_0(\bi),\\
    \eps_a(\bi)\eps_b(\bi)&= \eps_{a+b}(\bi),
    & \eps_a(\bi)&=(-1)^a\eps_a(-\bi),
    & \psi_r\eps_a(\bi) &=\eps_a(s_r\cdot\bi)\psi_r,
  \end{align*}
  \begin{align*}
    y_r\eps_a(\bi)&= \eps_a(\bi)y_r, &
    y_ry_s\eps_1(\bi)&= y_sy_r\eps_1(\bi),\\
    \psi_ry_{r+1}\eps_1(\bi)&= (y_r\psi_r+\delta_{i_r i_{r+1}})\eps_1(\bi), &
    y_{r+1}\psi_r\eps_1(\bi)&= (\psi_ry_r+\delta_{i_r i_{r+1}})\eps_1(\bi),
  \end{align*}
  \begin{align*}
    \psi_ry_s\eps_1(\bi)&= y_s\psi_r\eps_1(\bi),&&\text{if }s\neq r,r+1,\\
    \psi_r\psi_s\eps_1(\bi)&= \psi_s\psi_r\eps_1(\bi),&&\text{if }|r-s|>1,
  \end{align*}
  \begin{align*}
    \psi_r^2\eps_1(\bi)&= \begin{cases}
      0,&\text{if }i_r=i_{r+1},\\
      (y_r-y_{r+1})\eps_0(\bi),&\text{if }i_r\to i_{r+1},\\
      (y_{r+1}-y_r)\eps_0(\bi),&\text{if }i_r\leftarrow i_{r+1},\\
      \eps_1(\bi),&\text{otherwise,}
    \end{cases}\\
    \psi_r\psi_{r+1}\psi_r\eps_0(\bi)&=\begin{cases}
      \psi_{r+1}\psi_r\psi_{r+1}\eps_0(\bi)-\eps_1(\bi),
      &\text{if }i_r=i_{r+2}\to i_{r+1},\\
      \psi_{r+1}\psi_r\psi_{r+1}\eps_0(\bi)+\eps_1(\bi),
      &\text{if }i_r=i_{r+2}\leftarrow i_{r+1},\\
      \psi_{r+1}\psi_r\psi_{r+1}\eps_0(\bi),&\text{otherwise,}
    \end{cases}
  \end{align*}
  }%
  for $\bi,\bj\in I^\gamma$, $a,b\in\Z_2$ and all
  admissible~$r$ and~$s$. Let
  $\RpSn=\bigoplus_{\gamma\in \Qpm}\RpSn[\gamma]$.
\end{Definition}

It is routine to check that the relations in \autoref{D:RpSn}
are homogeneous with respect to the degree function
$\Deg\map\RpSn\Z_2\times\Z$ that is determined by
\[\Deg\psi_r\eps_0(\bi)=(1,-c_{i_r,i_{r+1}}),\quad
\Deg y_s=(1,2)\quad \text{and}\quad
\Deg\eps_a(\bi)=(a,0),
\]
for $1\le r<n$, $1\le s\le n$ and $\bi\in I^n$. Hence, $\RpSn$ is a
$(\Z_2\times\Z)$-graded algebra.

By \autoref{D:RpSn}, if $\bi\in I^\gamma$ then
$\eps_a(\bi)=\pm\eps_a(-\bi)$ so $\RpSn$ is generated by the elements
$\set{\psi_1,\dots,\psi_{n-1}}\cup\set{y_1,\dots,y_n}
\cup\set{\eps_a(\bi)|a\in\Z_2\text{ and }\bi\in I^\gamma_+}$,
where $I^\gamma_+=I^\gamma\cap I^n_+$. Similarly, set
$I^\gamma_-=I^\gamma\cap I^n_-$. We use~$I^\gamma$ in
\autoref{D:RpSn} because it compactly encodes a sign change
in the relation $\psi_r\eps_a(\bi)=\eps_a(s_r\cdot\bi)\psi$
when $r=2$.

We need a partial analogue of \autoref{L:ZeroSequence} for $\RpSn[\gamma]$.

\begin{Lemma}\label{L:ZeroSequence2}
  Suppose that $\bi\in I^\gamma$ and that $\bi=-\bi$. Then
  $\eps_a(\bi)=0$, for $a\in\Z_2$.
\end{Lemma}

\begin{proof}
  If $\bi=-\bi$ then $\eps_1(\bi)=-\eps_1(\bi)=0$. Hence,
  $\eps_0(\bi)=\eps_1(\bi)\eps_1(\bi)=0$.
\end{proof}

By forgetting the $\Z_2$-grading on $\RpSn[\gamma]$ we obtain a $\Z$-graded
algebra. Given the similarity of the relations in \autoref{D:klr} and
\autoref{D:RpSn} the next result should not surprise the reader.

\begin{Proposition}\label{P:RpSnIsomorphism}
  Suppose that $\gamma\in\Qpm$, $n\ge0$, $e\ne2$ and that~$2$ is
  invertible in~$\Zcal$. Then, as $\Z$-graded algebras,
  $\RpSn[\gamma]\cong\Ralpha[\gamma]$.
\end{Proposition}

\begin{proof}
  Define a map $\theta$ from the generators of $\RpSn[\gamma]$ to
  $\Ralpha[\gamma]$ by
  \begin{equation*}
    \theta(\psi_r) = \psi_r,\quad
    \theta(y_s) = y_s,\quad\text{ and }\quad
    \theta\(\eps_a(\bi)\) = e(\bi) + (-1)^ae(-\bi),
  \end{equation*}
  for $1\le r<n$, $1\le s\le n$, $a\in Z_2$ and $\bi\in I^\gamma$.  The
  relations in \autoref{D:RpSn} are very similar to those of
  \autoref{D:klr}, so it is straightforward to check that~$\theta$
  extends to an algebra homomorphism $\RpSn[\gamma]\to\Ralpha[\gamma]$.
  By definition, if $\bi\in I^\gamma$ then
  $e(\bi)=\frac12\theta\big(\eps_0(\bi)+\eps_1(\bi)\big)$.
  Therefore, the image of~$\theta$ contains
  all of the generators of $\Ralpha[\gamma]$. Hence, $\theta$ is
  surjective.

  Rather than proving directly that $\theta$ is an isomorphism we define
  an inverse map. Define $\vartheta$ to be the map from
  the set of non-zero generators of~$\Ralpha[\gamma]$ into~$\RpSn[\gamma]$ given by
  \[\vartheta(\psi_r)=\psi_r,\quad
  \vartheta(y_s)=y_s, \quad\text{and}\quad
  \vartheta\(e(\bi)\)=\tfrac12\(\eps_0(\bi)+\eps_1(\bi)\),
  \]
  for $1\le r<n$, $1\le s\le n$ and $\bi\in I^\gamma$. If $\bi\in I^\gamma$ and
  $a\in\Z_2$ then
  \[
  \vartheta\big(e(\bi)+(-1)^ae(-\bi)\big)
  = \frac12\big( \eps_0(\bi)+\eps_1(\bi)
  +(-1)^a \eps_0(-\bi)+(-1)^a\eps_1(-\bi)\big)
  =\eps_a(\bi).
  \]
  Hence, by \autoref{L:ZeroSequence2}, the image of $\vartheta$ contains
  all of the generators of~$\RpSn[\gamma]$ so, if it is a homomorphism, it is
  surjective.

  Now, since $\eps_a(\bi)\eps_b(\bi)=\eps_{a+b}(\bi)$, for
  all $\bi\in I^\gamma$
  and $a,b\in\Z_2$, by multiplying the relations in \autoref{D:RpSn} on
  the right by $\eps_1(\bi)$ the following additional relations
  hold in $\RpSn[\gamma]$:
  {\setlength{\abovedisplayskip}{2pt}
  \setlength{\belowdisplayskip}{1pt}
  \begin{align*}
    y_ry_s\eps_0(\bi)&= y_sy_r\eps_0(\bi),
  \end{align*}
  \begin{align*}
    \psi_ry_{r+1}\eps_0(\bi)&= (y_r\psi_r+\delta_{i_r i_{r+1}})\eps_0(\bi), &
    y_{r+1}\psi_r\eps_0(\bi)&= (\psi_ry_r+\delta_{i_r i_{r+1}})\eps_0(\bi),
  \end{align*}
  \begin{align*}
    \psi_ry_s\eps_0(\bi)&= y_s\psi_r\eps_0(\bi),&&\text{if }s\neq r,r+1,\\
    \psi_r\psi_s\eps_0(\bi)&= \psi_s\psi_r\eps_0(\bi),&&\text{if }|r-s|>1,
  \end{align*}
  \begin{align*}
    \psi_r^2\eps_0(\bi)&= \begin{cases}
      0,&\text{if }i_r=i_{r+1},\\
      (y_r-y_{r+1})\eps_1(\bi),&\text{if }i_r\to i_{r+1},\\
      (y_{r+1}-y_r)\eps_1(\bi),&\text{if }i_r\leftarrow i_{r+1},\\
      \eps_0(\bi),&\text{otherwise,}
    \end{cases}\\
    \psi_r\psi_{r+1}\psi_r\eps_1(\bi)&=\begin{cases}
      \psi_{r+1}\psi_r\psi_{r+1}\eps_1(\bi)-\eps_0(\bi),
      &\text{if }i_r=i_{r+2}\to i_{r+1},\\
      \psi_{r+1}\psi_r\psi_{r+1}\eps_1(\bi)+\eps_0(\bi),
      &\text{if }i_r=i_{r+2}\leftarrow i_{r+1},\\
      \psi_{r+1}\psi_r\psi_{r+1}\eps_1(\bi),&\text{otherwise,}
    \end{cases}
  \end{align*}}%
  for all admissible $r$ and $s$ and $\bi\in I^\gamma$. As it was
  for~$\theta$, it is now straightforward to verify that
  $\vartheta$ respects all of the relations of~$\Ralpha[\gamma]$. Consequently,
  $\vartheta$ extends to an algebra homomorphism $\Ralpha[\gamma]\rightarrow\RpSn[\gamma]$.

  In view of \autoref{L:ZeroSequence} and \autoref{L:ZeroSequence2}, the
  automorphisms $\theta\circ\vartheta$ and $\vartheta\circ\theta$ act as
  the identity on the non-zero generators of $\Ralpha[\gamma]$ and
  $\RpSn[\gamma]$, respectively. Therefore,~$\theta$ and $\vartheta$
  are mutually inverse isomorphisms and
  $\RpSn[\gamma]\bijection\Ralpha[\gamma]$ as (ungraded) algebras.

  It remains to observe that $\theta$ and
  $\vartheta$ respect the $\Z$-gradings on both algebras, but this is
  immediate from the definitions of~$\theta$ and~$\vartheta$. Hence,
  $\Ralpha[\gamma]\cong\RpSn[\gamma]$ as $\Z$-graded algebras, completing the
  proof.
\end{proof}

Define $h\in\RpSn$ to be \textbf{even} if $\Deg h=(0,d)$
and~$h$ is \textbf{odd} if $\Deg h=(1,d)$, for some $d\in\Z$. Let~$\RpSnp$ and
$\RpSnp[-]$ be the sets of even and odd elements in~$\RpSn[\gamma]$,
respectively. Then~$\RpSnp$ is a subalgebra of~$\RpSn[\gamma]$ and
\begin{equation}\label{clifforddecomp}
  \RpSn[\gamma]=\RpSnp\oplus\RpSnp[-]
\end{equation}
as $\Zcal$-modules.  Moreover, as we next show, $\RpSnp$ is isomorphic
to $\RAn$ under the isomorphism of \autoref{P:RpSnIsomorphism}.

\begin{Corollary}\label{C:EvenBit}
  Suppose that $n\ge0$ and that $2$ is invertible in $\Zcal$. Let
  $\gamma\in \Qpm$. Then $\RpSnp\cong\RAn_\gamma$
  as $\Z$-graded algebras.
\end{Corollary}

\begin{proof}
  Under the isomorphism $\RpSn[\gamma]\cong\Ralpha[\gamma]$ of
  \autoref{P:RpSnIsomorphism}, the images of the
  even generators of~$\RpSn[\gamma]$ are $\sgn$-invariant and $\sgn$ multiplies
  the images of the odd generators by~$-1$. Hence, $\theta$ restricts to
  an isomorphism $\RpSnp\bijection\RAn_\gamma$.
\end{proof}

We can now prove \autoref{T:MainRelations} from the introduction.

\begin{proof}[Proof of \autoref{T:MainRelations}]
  Let $A_\gamma$ be the abstract algebra with the presentation given in
  \autoref{T:MainRelations}. By \autoref{C:EvenBit}, to prove
  \autoref{T:MainRelations} it is enough to show that $A_\gamma\cong \RpSnp$.
  Define a map $\Theta\map{A_\gamma}\RpSnp$ by
  \[\Psi_r(\bi)\mapsto\psi_r\eps_1(\bi),\quad
  Y_s(\bi)\mapsto y_s\eps_1(\bi)\quad\text{and}\quad
  \eps(\bi)\mapsto \eps_0(\bi),
  \]
  for all $\bi\in I^\gamma$, $1\le r\le n$ and $1\le s<n$. Using
  \autoref{D:RpSn}, and the relations in the proof of
  \autoref{P:RpSnIsomorphism}, it is straightforward to check that all
  of the relations in~$A_\gamma$ are satisfied in $\RpSnp$, so $\Theta$
  extends to an algebra homomorphism from~$A_\gamma$ to~$\RpSnp$.

  By \autoref{D:RpSn}, the algebra $\RpSnp$ is generated by
  arbitrary products of the generators of~$\RpSn$ such that the
  resulting element is even. However, the only even generators
  of~$\RpSn$ are the idempotents $\eps_0(\bi)$, for $\bi\in I^n$, so
  $\RpSnp$ is generated by these idempotents together with
  all words of even length in the odd generators of~$\RpSn$. As
  $\psi_r\eps_1(\bi)=\eps_1(s_r\cdot\bj)\psi_r$ and
  $y_s\eps_1(\bi)=\eps_1(\bi)y_s$, for admissible $r$, $s$ and $\bi\in I^n$,
  it follows that $\RpSnp$ is generated by the images of
  $A_\gamma$ under $\Theta$. Hence, $\Theta$ is surjective.

  The algebra $\RpSnp$ is defined by generators and relations, so
  $\RpSnp$ is the subalgebra of $\RpSn$ generated by the
  words of even length in the generators of~$\RpSn$ modulo the even part of the
  relational ideal that defines~$\RpSn$. The only even relations
  for~$\RpSn$ are the idempotent relations and the commutation
  relations (these are the relations appearing in the first three lines of
  the relations in \autoref{D:RpSn}). Hence, up to multiplication by an
  idempotent~$\eps_0(\bi)$, all of the even relations are
  essentially trivial. Therefore, the even component of
  the relational ideal for~$\RpSn$ is generated by words of even
  length in the odd relations for~$\RpSn$. In turn, all of the even
  products of the odd relations are products of the even relations given
  in the proof of \autoref{P:RpSnIsomorphism}, together with the even
  idempotent and commutation relations. All of these relations are the
  images under $\Theta$ of the relations of~$A_\gamma$. Hence, $\Theta$ is an
  isomorphism and the result follows.
\end{proof}

\section{The seminormal form}\label{S:SeminormalForm}
  To prove \autoref{T:Main} we will work mainly in the setting of the
  semisimple representation theory of~$\HSn$. The idea is to show that
  the fixed-point subalgebras of $\HAn=\HSn^\#$ and $\RSn=\RSn^\sgn$
  coincide under the Brundan-Kleshchev isomorphism of \autoref{T:BKiso}.
  Unfortunately, as shown by \autoref{Ex:BadBKRestriction}, this is not true. To
  get around this we use the machinery developed in
  \cite{HuMathas:SeminormalQuiver} to construct a new isomorphism
  $\RSn\bijection\HSn$ that does restrict to an isomorphism
  $\RAn\bijection\HAn$.

  \subsection{Tableau combinatorics}\label{S:tableaux}
  This section recalls the partition and tableau
  combinatorics  that are needed in this paper.

  A \textbf{ partition} $\la=(\lambda_1,\lambda_2,\dots)$ is a weakly decreasing
  sequence of non-negative integers. The integers $\lambda_r$ are the
  \textbf{parts} of $\lambda$, for $r\ge1$, and $\lambda$ is a partition of~$n$ if
  $\left|\lambda\right|=n$, where $\left|\lambda\right|=\lambda_1+\lambda_2+\cdots$.
  Let $\Parts$ be the set of partitions of~$n$.

  The \textbf{Young diagram} of a partition $\lambda$ is the set
  $\set{(r,c)|1\le c\le\lambda_r\text{ for }r\ge1}$, which we represent as
  a collection of left-justified boxes in the plane, with $\lambda_r$
  boxes in row~$r$ and with rows ordered from top to bottom by increasing
  row index. We identify a partition with its diagram. The partition
  $\la'$ with $\la'_r=\#\set{c\geq 1| \lambda_r\geq c}$ is the partition
  \textbf{conjugate} to $\la$.

  Suppose that $\lambda\in\Parts$. A \textbf{$\la$-tableau} is a bijective filling of
  the boxes of $\la$ with the numbers $1,2,\ldots,n$. If $\t$ is a
  $\lambda$-tableau then it has \textbf{shape} $\lambda$ and we write
  $\sh(\t)=\lambda$. For $m\ge1$ let $\t_{\downarrow m}$ be the subtableau
  of~$\t$ that contains the numbers $1,2,\dots.m$.

  A tableau is \textbf{standard} if its entries increase from left to
  right along each row and from top to bottom down each column.  Hence,
  $\t$ is standard if and only if~$\t_{\downarrow m}$ is standard for
  $1\le m\le n$.  Let~$\Std(\lambda)$ be the set of standard $\la$-tableaux
  and let
  \[\Std(\Parts)=\bigcup_{\lambda\in\Parts}\Std(\lambda)\qquad\text{and}\qquad
  \Std^2(\Parts)=\bigcup_{\lambda\in\Parts}\Std(\lambda)\times\Std(\lambda).
  \]

  If $\t$ is a standard $\lambda$-tableau then the \textbf{conjugate}
  tableau $\t'$ is the standard $\la'$-tableau obtained by swapping the rows and
  columns of~$\t$.

  The \textbf{initial $\lambda$-tableau} $\t^\lambda$ is the
  $\lambda$-tableau obtained by inserting the numbers $1, 2, \ldots,n$ in
  order along the rows of~$\lambda$, from left to right and then top to
  bottom. The \textbf{co-initial} tableau $\t_\lambda$ is the conjugate
  of~$\t^{\lambda'}$. Then $\t_\lambda$ is the unique $\lambda$-tableau that has
  the numbers $1,2,\dots,n$ entered in order down the columns
  of~$\lambda$, from left to right.

  Recall that $I=\Z/e\Z$. Let $\t$ be a standard tableau and
  suppose that $m$ appears in row~$r$ and column~$c$ of ~$\t$, where $1\le
  m\le n$. The
  \textbf{content} and \textbf{$e$-residue} of~$m$ in~$\t$ are given by
  \[
  c_m(\t)=c-r\in\Z\qquad\text{and}\qquad \res_m(\t)= c-r+e\Z\in I.
  \]
  respectively. The \textbf{$e$-residue sequence} of the tableau~$\t$ is
  the $n$-tuple
  \[
  \res(\t)=\(\res_1(\t),\dots,\res_n(\t)\)\in I^n.
  \]

  Given a sequence $\bi\in I^n$ let
  $\Std(\bi)=\set{\t\in\Std(\Parts)\mid\res(\t)=\bi}$ be the set of
  standard tableaux with residue sequence~$\bi$.

  \subsection{Seminormal forms}\label{Sect:Seminormal}
  To prove \autoref{T:MainTheorem} we make extensive use of the
  semisimple representation theory of $\HSn$ using seminormal forms.
  This section introduces Jucys-Murphy elements and seminormal forms and
  proves some basic facts relating seminormal forms and the
  $\#$-involution.

  Throughout this section we fix a field $\K$ and a non-zero scalar
  $t\in\K$. Let~$\KHSn$ be the Iwahori-Hecke algebra of~$\Sn$ over~$\K$ with
  parameter~$t$.

  If $k\in\Z$  define the \textbf{quantum integer} to be the scalar
  \[[k]_t=\begin{cases}
    \phantom{-(}(1+t+\dots+t^{k-1}),& \text{if }k\ge0,\\
    -(t^{-1}+t^{-2}+\dots+t^k),&\text{if }k<0.
  \end{cases}
  \]
  When~$t$ is understood we write $[k]=[k]_t$.

  We need a well-known result, which is easily proved by induction
  on~$n$. To state this, define the \textbf{Poincar\'e polynomial} of~$\KHSn$ to
  be $P_\H(t)=[1][2]\dots[n]\in\K$.

  \begin{Lemma}[\protect{See \cite[Lemma 3.34]{M:ULect}}]
    \label{L:separation}
    Suppose that $P_\H(t)\ne0$ and that
    $\s,\t\in\Std(\Parts)$. Then $\s=\t$ if and only if
    $[c_r(\s)]=[c_r(\t)]$, for $1\le r\le n$.
  \end{Lemma}

  We assume for the rest of this section that $P_\H(t)\ne0$. In fact,
  the results that follow imply that~$\KHSn$ is semisimple if and only if
  $P_\H(t)\ne0$ and, in turn, this is equivalent to the condition in
  \autoref{L:separation}.

  If $\t$ is a tableau and $1\le r\le n$ then the
  \textbf{axial distance} from $r+1$ to $r$ in~$\t$ is
  \begin{equation}\label{E:rhodef}
    \rho_r(\t)=c_r(\t)-c_{r+1}(\t)\in\Z.
  \end{equation}
  By definition, $-n<\rho_r(\t)<n$ so $[\rho_r(\t)]\ne0$ if $P_\H(t)\ne0$.

  The next definition will provide us with the framework to prove
  \autoref{T:MainTheorem}.

  \begin{Definition}[\protect{%
    Hu-Mathas~\cite[Definition~3.5]{HuMathas:SeminormalQuiver}}]
    \label{D:alphaSNCS}
    Suppose that $P_\H(t)\ne0$.  A \textbf{$*$-seminormal coefficient
    system} is a set of scalars
    \[\balpha=\set{\alpha_r(\t)\in\K|1\leq r<n\text{ and }
    \t\in\Std(\Parts)}
    \]
    such that if~$\t\in\Std(\Parts)$ and $1\le r<n$ then:
    \begin{enumerate}
      \item $\alpha_r(\t)=0$ whenever $s_r\t$ is not standard.
      \item $\alpha_r(\t)\alpha_k(s_r\t)=\alpha_k(\t)\alpha_r(s_k\t)$
      whenever $1\le k<n$ and $|r-k|>1$.
      \item
      $ \alpha_r(s_{r+1}s_r\t)\alpha_{r+1}(s_r\t)\alpha_r(\t)
      =\alpha_{r+1}(s_rs_{r+1}\t)\alpha_r(s_{r+1}\t)\alpha_{r+1}(\t)
      $
      if $r\ne n-1$,
      \item if $\v=s_r\t\in\Std(\Parts)$ then
      \begin{equation*}\label{E:sncf}
        \alpha_r(\t)\alpha_r(\v)
        =\frac{[1+\rho_r(\t)][1+\rho_r(\v)]}{[\rho_r(\t)][\rho_r(\v)]}.
      \end{equation*}
    \end{enumerate}
  \end{Definition}

  Many examples of seminormal coefficient systems are given in
  \cite[\S3]{HuMathas:SeminormalQuiver}. For example, $\set{\alpha_r(\t)}$
  is seminormal coefficient system, where
  $\alpha_r(\t)=\tfrac{[1+\rho_r(\t)]}{[\rho_r(\t)]}$ whenever $\t,
  s_r\t\in\Std(\Parts)$. In \autoref{S:AlternaatingCoefficients} we fix
  a particular choice of seminormal coefficient system but until then
  we will work with an arbitrary coefficient system.

  For $k=1,2,\ldots,n$ the \textbf{Jucys-Murphy} element $L_k\in \HO$ is
  defined by
  \[
  L_k=\sum_{j=1}^{k-1}t^{j-k} T_{(k-j,k)}.
  \]

  A basis $\set{f_{\s\t}|(\s,\t)\in\Std^2(\Parts)}$ of $\KHSn$ is a \textbf{seminormal
  basis} if
  \[L_kf_{\s\t}=[c_k(\s)]f_{\s\t}\qquad\text{and}\qquad
  f_{\s\t}L_k=[c_k(\t)]f_{\s\t},
  \]
  for all $(\s,\t)\in \Std^2(\Parts)$ and $1\le k\le n$.
  The basis $\{f_{\s\t}\}$ is a \textbf{$*$-seminormal basis} if, in addition,
  $f_{\s\t}^*=f_{\t\s}$, for all $(\s,\t)\in\Std^2(\Parts)$, where $*$ is
  the unique anti-isomorphism of $\KHSn$ that fixes $T_1,\dots,T_{n-1}$.

  Recall that $P_\H(t)=[1][2]\dots[n]$.

  \begin{Theorem}[\protect{%
    The seminormal form~\cite[Theorem 3.9]{HuMathas:SeminormalQuiver}}]
    \label{T:SeminormalForm}
    Suppose that $P_\H(t)\ne0$ and that $\balpha$ is a seminormal
    coefficient system for $\KHSn$. Then:
    \begin{enumerate}
      \item The algebra $\KHSn$ has a unique $*$-seminormal basis
      $\set{f_{\s\t}|(\s,\t)\in\Std^2(\Parts)}$ such that
      \[
      f_{\s\t}^*=f_{\t\s},\quad L_kf_{\s\t}=[c_k(\s)]f_{\s\t}\quad\text{and}\quad
      T_rf_{\s\t}=\alpha_r(\s)f_{\u\t}-\frac{1}{[\rho_r(\s)]}f_{\s\t},
      \]
      where $\u=(r,r+1)\s$. $($Set $f_{\u\t}=0$ if $\u$ is not standard.$)$
      \item For $\t\in\Std(\Parts)$ there exist non-zero scalars
      $\gamma_\t\in\K$ such that
      $f_{\s\t}f_{\u\v}=\delta_{\t\u}\gamma_\t f_{\s\v}$ and
      $\set{\tfrac1{\gamma_\t}f_{\t\t}| \t\in\Std(\Parts)}$ is a complete
      set of pairwise orthogonal primitive idempotents.
      \item The $*$-seminormal basis
      $\set{f_{\s\t}|(\s,\t)\in\Std^2(\Parts)}$ is uniquely determined by
      the $*$-seminormal coefficient system $\balpha$ and the scalars
      $\set{\gamma_{\t^\lambda}|\lambda\in\Parts}$.
    \end{enumerate}
  \end{Theorem}

%

  By \autoref{T:SeminormalForm}(b), if $\t\in\Std(\Parts)$ then
  $F_\t=\tfrac1{\gamma_\t}f_{\t\t}$ is a primitive idempotent in $\KHSn$. As is
  well-known (see, for example, \cite[(3.2)]{HuMathas:SeminormalQuiver}),
  \begin{equation*}\label{E:Ft}
    F_\t= \prod_{k=1}^n\prod_{\substack{\s\in\Std(\Parts)\\c_k(\s)\ne c_k(\t)}}
    \frac{L_k-[c_k(\s)]}{[c_k(\t)]-[c_k(\s)]}.
  \end{equation*}
  In particular, the idempotent $F_\t$ is independent of the choice of
  seminormal basis.

  Let $\Lscr$ be the commutative subalgebra generated by the
  Jucys-Murphy elements. \autoref{T:SeminormalForm}(a) implies that, as
  an $(\Lscr,\Lscr)$-bimodule, $\KHSn$ decomposes as
  \begin{equation}\label{E:LLBimodule}
    \KHSn=\bigoplus_{(\s,\t)\in\Std^2(\Parts)}H_{\s\t},
  \end{equation}
  where $H_{\s\t}=\K f_{\s\t}$. Equivalently,
 \[ H_{\s\t}=\set{h\in\KHSn|L_k h=[c_k(\s)]h\text{ and }
    hL_k=[c_k(\t)]h\text{ for }1\le k\le n},
  \]
  for $(\s,\t)\in\Std^2(\Parts)$.

  Let $\#$ be the hash involution from \autoref{E:hash} on~$\KHSn$. Then
  $T_r^\#=-T_r+t-1$, for $1\le r<n$.

  \begin{Lemma}\label{L:LFHash}
    Suppose that $1\le k\le n$ and $\s\in\Std(\Parts)$. Then
    \[
    L_k^\# f_{\s\s}=[c_k(\s')]f_{\s\s}.
    \]
  \end{Lemma}
  \begin{proof}
    For $1\le k\le n$ set
    $\hat{L}_k=t^{1-k}T_{k-1}T_{k-2}\cdots T_2T_1^2T_2\cdots T_{k-2}T_{k-1}$.
    It is well-known and easy to prove that
    $\hat{L}_k=(t-1)L_k+1$; see, for example,
    \cite[Exercise~3.6]{M:ULect}.
    By \autoref{T:SeminormalForm}(a), $L_k f_{\s\s}=[c_k(\s)]f_{\s\s}$, so
    $\hat L_kf_{\s\s}=t^{c_k(\s)}f_{\s\s}$. Now
    $(\hat{L}_k)^\#=\hat{L}_k^{-1}$ since $T_r^\#=-tT_r^{-1}$ by
    \autoref{E:hash}, for $1\le r<k\le n$. Therefore,
    \[\hat{L}_k^\# f_{\s\s}=\hat{L}_k^{-1}f_{\s\s}
    =t^{-c_k(\s)}f_{\s\s}
    =t^{c_k(\s')}f_{\s\s},
    \]
    where the last equality follows because $c_k(\s')=-c_k(\s)$ for $1\le k\le n$.
    Hence, $L_k^\#f_{\s\s}=[c_k(\s')]f_{\s\s}$ as claimed.
  \end{proof}

  \begin{Lemma}\label{L:Fhash}
    Suppose that $\s\in \Std(\Parts)$. Then $F_\s^\#=F_{\s'}$.
  \end{Lemma}
  \begin{proof}
    Since $F_\s=\frac{1}{\gamma_\s}f_{\s\s}$, applying \autoref{L:LFHash} gives
    \[
    L_kF_\s^\#=(L_k^\#F_\s)^\#=([c_k(\s')]F_\s)^\#=[c_k(\s')]F_\s^\#.
    \]
    Similarly $F_\s^\#L_k=[c_k(\s')]F_\s^\#$. Therefore, $F_\s^\#\in H_{\s'\s'}$ in the
    decomposition of \autoref{E:LLBimodule}. As $F_\s$ is an idempotent, and $\#$ is
    an algebra automorphism, it follows that $F_\s^\#=F_{\s'}$ since this is the
    unique idempotent in $H_{\s'\s'}=\K F_{\s'}=\K f_{\s'\s'}$.
  \end{proof}

  \begin{Corollary}\label{C:fttHash}
    Suppose that $\s\in\Std(\Parts)$. Then
    $f_{\s\s}^\#=\dfrac{\gamma_\s}{\gamma_{\s'}}f_{\s'\s'}.$
  \end{Corollary}

  \begin{proof} Using \autoref{T:SeminormalForm}(b) and \autoref{L:Fhash},
    $f_{\s\s}^\#=\frac1{\gamma_\s}F_\s^\#=\frac1{\gamma_\s}F_{\s'}
    =\frac{\gamma_{\s'}}{\gamma_\s}f_{\s'\s'}$.
  \end{proof}

  \begin{Lemma}\label{L:futhash}
    Let $\s,\u\in\Std(\Parts)$ be standard tableaux such that
    $\u=(r,r+1)\s$, for some integer~$r$ with $1\le r<n$. Then
    \[
    f_{\u\s}^\#=-\frac{\alpha_r(\s')\gamma_\s}{\alpha_r(\s)\gamma_{\s'}}f_{\u'\s'}.
    \]
  \end{Lemma}

  \begin{proof}
    By \autoref{T:SeminormalForm}(a),
    $f_{\u\s}=\frac{1}{\alpha_r(\s)}\Big(T_r+\frac1{[\rho_r(\s)]}\Big)f_{\s\s}$.
    Recall that $T_r^\#=-tT_r^{-1}=-T_r+t-1$. Therefore, using
    \autoref{E:hash} and \autoref{C:fttHash} for the second equality,
    \begin{align*}
      f_{\u\s}^\#&=\frac{1}{\alpha_r(\s)}\Big(T_r
       +\frac1{[\rho_r(\s)]}\Big)^\#f_{\s\s}^\#\\
      &=\frac{\gamma_\s}{\alpha_r(\s)\gamma_{\s'}}\Big(-T_r+t-1
       +\frac1{[\rho_r(\s)]}\Big)f_{\s'\s'}\\
      &=-\frac{\gamma_\s}{\alpha_r(\s)\gamma_{\s'}}\Big(T_r
       -\frac{t^{\rho_r(\s)}}{[\rho_r(\s)]}\Big)f_{\s'\s'}\\
      &=-\frac{\gamma_\s}{\alpha_r(\s)\gamma_{\s'}}\Big(T_r
       +\frac1{[\rho_r(\s')]}\Big)f_{\s'\s'},
    \end{align*}
    since
    $[\rho_r(\s)]=-t^{\rho_r(\s)}[-\rho_r(\s)]=-t^{\rho_r(\s)}[\rho_r(\s')]$.
    Hence, the result follows by another application of \autoref{T:SeminormalForm}(a).
  \end{proof}

  By \autoref{T:SeminormalForm}(c) any $*$-seminormal basis is uniquely
  determined by a seminormal coefficient system and a choice of scalars
  $\set{\gamma_{\t^\lambda}|\lambda\in\Parts}$. For completeness we
  determine these scalars for the seminormal basis
  $\set{f_{\s\t}^\#|(\s,\t)\in\Std^2(\Parts)}$.  Recall from
  \autoref{S:tableaux} that $\t_{\lambda}=(\t^{\lambda'})'$ is the
  co-initial $\lambda$-tableau.

  \begin{Proposition}\label{P:AltCS}
    The seminormal basis $\set{f_{\s\t}^\#|(\s,\t)\in\Std(\Parts)}$ of
    $\KHSn$ is the seminormal basis determined by the seminormal
    coefficient system
    \[
    \set{-\alpha_r(\s')|\s\in\Std(\Parts)\text{ and }1\le r<n}
    \]
    together with the $\gamma$-coefficients
    $\set{\gamma_{\t_\lambda}|\lambda\in\Parts}$.  That is, if
    $(\s,\t)\in\Std^2(\Parts)$ then
    \[
    L_kf_{\s\t}^\#=[c_r(\s')]f_{\s\t}^\#,\quad
    f_{\s\t}^\#L_k=[c_r(\t')]f_{\s\t}^\#\quad\text{and}\quad
    T_r f_{\s\t}^\# = -\alpha_r(\s)f_{\u\t}^\#-\frac1{[\rho_r(\s')]}f_{\s\t}^\#,
    \]
    where $\u=(r,r+1)\s$, $1\le k\le n$ and $1\le r<n$. Moreover,
    $f_{\s\t}^\#f_{\u\v}^\#=\delta_{\t\u}\gamma_\t f_{\s\v}^\#$,
    for $(\s,\t),(\u,\v)\in\Std^2(\Parts)$.
  \end{Proposition}

  \begin{proof}
    If $1\le k\le n$ then
    $L_kf_{\s\t}^\#=[c_r(\s')]f_{\s\t}^\#$ and
    $f_{\s\t}^\#L_k=[c_r(\t')]f_{\s\t}^\#$ by \autoref{L:LFHash}.
    Using \autoref{T:SeminormalForm}(a),
    \begin{align*}
      T_rf_{\s\t}^\# &= \Big(T_r^\# f_{\s\t}\Big)^\#
      = \Big((-T_r+t-1)f_{\s\t}\Big)^\#\\
      &=\Big(-\alpha_r(\s)f_{\u\t} +(t-1+\frac1{[\rho_r(\s)]})f_{\s\t}\Big)^\#\\
      &=\Big(-\alpha_r(\s)f_{\u\t}-\frac1{[\rho_r(\s')]}f_{\s\t}\Big)^\#\\
      &=-\alpha_r(\s)f_{\u\t}^\# -\frac1{[\rho_r(\s')]}f_{\s\t}^\#.
    \end{align*}
    Similarly,
    $f_{\s\t}^\#f_{\u\v}^\#=(f_{\s\t}f_{\u\v})^\#=\delta_{\t\u}\gamma_\t f_{\s\v}^\#$.
    By \autoref{L:futhash}, $f_{\s\t}^\#\in H_{\s'\t'}$, so the $\alpha$-coefficient
    corresponding to~$f_{\s\t}^\#$ is naturally indexed by~$\s'$ (and not by~$\s$).
    Similarly, the labelling for the $\gamma$-coefficients involves conjugation
    because $F_\t=\frac1{\gamma_{\t'}}f_{\t'\t'}^\#$ by \autoref{C:fttHash}. Hence,
    the result follows by \autoref{T:SeminormalForm}.
  \end{proof}

  \subsection{Idempotent subrings and KLR generators}\label{S:KLRGenerators}
  We are almost ready to introduce the generators of $\HO$ that we need to
  prove \autoref{T:MainTheorem}. This section defines roughly half of
  these generators. The definition of these elements involves lifting idempotents from the
  non-semisimple case to the semisimple case and to do this we need to
  place additional constraints upon the rings that we work over.

  If $\O$ is a ring let $\J=\J(\O)$ be the Jacobson radical of~$\O$.

  \begin{Definition}[\protect{\cite[Definition~4.1]{HuMathas:SeminormalQuiver}}]
    \label{D:idempotentSub}
    Suppose that $\O$ is a subring of a field~$\K$ and that $t\in\O^\times$.
    The pair $(\O,t)$ is an \textbf{$e$-idempotent subring} of $\K$ if:
    \begin{enumerate}
      \item The Poincar\'e polynomial $P_\H(t)$ is a non-zero element of~$\O$,
      \item If $k\notin e\Z$ then $[k]_t$ is invertible in $\O$,
      \item If $k\in e\Z$ then $[k]_t\in \mathcal{J}(\O)$.
    \end{enumerate}
  \end{Definition}

  Condition~(a) ensures that \autoref{L:separation} and
  \autoref{T:SeminormalForm} apply and, in particular, that $\KHSn$ has a
  seminormal basis $\set{f_{\s\t}}$. As discussed in
  \cite[Example~4.2]{HuMathas:GradedCellular}, when considering the Hecke
  algebra $\HSn$ defined over the field~$F$ with parameter $\xi\in
  F^\times$, one natural choice of idempotent subring is to let $x$ be an
  indeterminate over~$F$ and set $\K=F(x)$, $t=x+\xi$ and
  $\O=F[x,x^{-1}]_{(x)}$. Note that $\m=x\O$ is the unique maximal ideal of~$\O$
  and that $(\K,\O,F)$ is a modular system with $\HSn\cong\HO\otimes_\O
  F$, where $F$ is considered as an $\O$-module by letting~$x$ act as
  multiplication by~$0$.

  The hash involution $\#$ from \autoref{E:hash} is well-defined on $\HO$.
  Let $\HO(\An)$ the $\O$-subalgebra of $\#$-fixed points in~$\HO$.
  Then $\HAn\cong\HO(\An)\otimes_\O F$.

  Recall that if $\bi\in I^n$ then
  $\Std(\bi)=\set{\t\in\Std(\Parts)|\res(\t)=\bi}$. Using an idea that
  goes back to Murphy~\cite{M:Nak}, the \textbf{$\bi$-residue idempotent}
  is defined to be the element
  \begin{equation*}\label{E:ResidueIdempotent}
    \fo = \sum_{\t\in\Std(\bi)}\frac1{\gamma_\t}f_{\t\t}
    = \sum_{\t\in\Std(\bi)}F_\t.
  \end{equation*}
  By \autoref{T:SeminormalForm}(b), if $\s\in\Std(\bj)$ and $\t\in\Std(\bk)$
  are tableaux of the same shape then $\fo f_{\s\t}=\delta_{\bi\bj}f_{\s\t}$ and
  $f_{\s\t}\fo=\delta_{\bi\bk}f_{\s\t}$, for $\bi,\bj,\bk\in I^n$.

  By definition, $\fo\in\H_t^\K(\Sn)$ but, in fact,  $\fo\in\HO$.

  \begin{Lemma}\label{L:foIdempotents}
    Suppose that $\bi\in I^n$. Then $\fo\in\HO$ and $(\fo)^\#=\fo[-\bi]$.
  \end{Lemma}

  \begin{proof}
    Since $(\O,t)$ is an idempotent subring, $\fo\in\HO$ by
    \cite[Lemma~4.5]{HuMathas:SeminormalQuiver}.  To prove that
    $(\fo)^\#=\fo[-\bi]$ first observe that $\s\in\Std(\bi)$ if and only if
    $\s'\in\Std(-\bi)$. Therefore, by \autoref{L:Fhash},
    \[
    (\fo)^\#=\sum_{\s\in\Std(\bi)}F_\s^\#=\sum_{\s\in\Std(\bi)}F_{\s'}=\fo[-\bi]
    \]
    as claimed.
  \end{proof}

  Following \cite{BK:GradedKL,HuMathas:SeminormalQuiver}, define
  $M_r=1-L_r+tL_{r+1}$, for $1\le r\le n$.  If $(\s,\t)\in\Std^2(\Parts)$ then
  it follows easily using \autoref{T:SeminormalForm}(a) that
  \begin{equation}\label{E:Mr}
    M_r f_{\s\t}=t^{c_r(\s)}[1-\rho_r(\s)]f_{\s\t}.
  \end{equation}
  The next result says that these elements are invertible when projected
  onto certain residue idempotents $\fo$, for $\bi\in I^n$.

  \begin{Corollary}[\protect{\cite[Corollary~4.8]{HuMathas:SeminormalQuiver}}]
    \label{C:MrInvert}
    Suppose that $\bi\in I^n$ and $i_r\ne i_{r+1}+1$, for $1\le r<n$ and . Then
    \[\Sum_{\s\in\Std(\bi)}\frac{t^{-c_r(\t)}}{[1-\rho_r(\s)]}F_\s\in\HO.\]
  \end{Corollary}

  In view of \autoref{C:MrInvert}, if $i_r\ne i_{r+1}+1$ define the formal symbol
  \begin{align*}
    \frac{1}{M_r}\fo
    &=\Sum_{\s\in\Std(\bi)}\frac{t^{-c_r(\t)}}{[1-\rho_r(\s)]}F_\s.
    \end{align*}
    This abuse of notation is justified because
    $\frac1{M_r}\fo M_r=\fo$ 
    by \autoref{E:Mr}. We will use these elements to define the
    KLR-generators of~$\HO$ that we use to prove \autoref{T:Main}.

    The results of~\cite{HuMathas:SeminormalQuiver} depend upon choosing an
    arbitrary  section of the natural quotient map $\Z\twoheadrightarrow\Z/e\Z$. In
    this paper we are far less flexible and need to use a particular
    section of this map. If
    $i\in I$ let $\i\ge0$ be the smallest non-negative integer such that
    $i=\i+e\Z$. (If $e=\infty$ set $\hat i=i$.) This defines an
    embedding $I\hookrightarrow\Z; i\mapsto\i$.  For $\bi\in I^n$ set
    $\rho_r(\bi)=\i_r-\i_{r+1}$, for $1\le r<n$.

    By \autoref{T:SeminormalForm}, the identity element of~$\HO$ can
    be written as $1=\sum_{\bi\in I^n}\fo$. So if $h\in\HO$ then
    $h=\sum_{\bi\in I^n} h\fo$ is uniquely determined by its
    projection onto the idempotents $\fo$.

    \begin{Definition}[\protect{\cite[Definition~4.14]{HuMathas:SeminormalQuiver}}]
      \label{D:psip}
      Fix an integer $1\le r<n$ and define the element
      $\psip_r=\sum_{\bi\in I^n}\psip_r\fo$ by
      \[
      \psip_r\fo=\begin{cases}
        (1+T_r)\frac{t^{\i_r}}{M_r}\fo,&\text{if }i_r=i_{r+1},\\
        (T_rL_r-L_rT_r)t^{-\i_r}\fo,&\text{if }i_r\to i_{r+1},\\
        (T_rL_r-L_rT_r)\frac{1}{M_r}\fo,&\text{otherwise}.
      \end{cases}
      \]
      For $1\le s\le n$ define $\yp_s=\sum_{\bi\in I^n}t^{-\i_s}(L_s-[\i_s])\fo$.
    \end{Definition}

    Recall from \autoref{S:GeneratorsRelations} that $\Qpm=Q^+_n/{\sim}$.
    For $\alpha\in Q^+$ define
    \[
    \HO[\alpha]=\HO\fo[\alpha],\qquad\text{where }\quad
    \fo[\alpha]=\sum_{\ei\in I^\alpha}\fo.
    \]
    For $\gamma\in\Qpm$ set
    $\fo[\gamma]=\sum_{\alpha\in\gamma}\fo[\alpha]$ and set
    $\HO[\gamma]=\bigoplus_{\alpha\in\gamma}\HO[\alpha]=\HO\fo[\gamma]$.

    \begin{Proposition}\label{P:HOdecomp}
      Suppose that $(\O,t)$ is an idempotent subring. Then $\fo[\gamma]$ is
      a central idempotent in $\HO$ and
      $\HO=\displaystyle\bigoplus_{\gamma\in\Qpm}\HO[\gamma]$.
    \end{Proposition}

      \begin{proof}
        By \autoref{L:foIdempotents}, $\fo[\gamma]\in\HO$ and it follows from
        \autoref{T:SeminormalForm} that $\fo[\gamma]$ is a central
        idempotent and that $1=\sum_{\gamma\in \Qpm}\fo[\gamma]$. Hence,
        $\HO[\gamma]=\fo[\gamma]\HO\fo[\gamma]$ is a subalgebra of~$\HO$ and
        $\HO=\bigoplus_{\gamma\in \Qpm}\HO[\gamma]$.
      \end{proof}

      By \cite{JM:Schaper}, for $\alpha\in Q^+_n$ the algebras
    $\HO[\alpha]\otimes_\O F$ are indecomposable two-sided
    ideals of~$\HSn$. Later we need the counterpart of \autoref{C:RAnBlocks}
    for $\HAn$. If $\gamma\in\Qpm$ then $(\fo[\gamma])^\#=\fo[\gamma]$ by
    \autoref{L:foIdempotents}. Therefore, $\#$ restricts to an automorphism
    of~$\HO[\gamma]$. Define
    \begin{equation}\label{E:HAnBlocks}
      \HO(\An)_\gamma=\bigl(\HO[\gamma]\bigr)^\#
                     =\set{h\in\HO[\gamma]|h^\#=h}
                     =\HO(\An)\fo[\gamma].
    \end{equation}
    As $1=\sum_\gamma\fo[\gamma]$, \autoref{P:HOdecomp} immediately implies
    the following.

    \begin{Corollary}\label{C:HAnDecomp}
      Suppose that $(\O,t)$ is an idempotent subring. Then
      \[\HO(\An)=\bigoplus_{\gamma\in\Qpm}\HO(\An)_\gamma.\]
    \end{Corollary}

    The subalgebra $\HO(\An)_\gamma$ is a block of~$\HO(\An)$ in the
    sense that it is a two-sided ideal and a direct summand
    of~$\HO(\An)$.  Let $F$ be a field that is an $\O$-algebra and set
    $\FHAn_\gamma=\HO(\An)_\gamma\otimes_\O F$. Then $\FHAn_\gamma$ is
    almost always indecomposable. See \autoref{T:RAnBlocks} for
    the precise statement.

      \begin{Theorem}[\protect{%
        Hu-Mathas~\cite[Theorem~A]{HuMathas:SeminormalQuiver}}]
        \label{T:HMGradedIso}
        Suppose that  $\mathcal{K}$ is a field, $\gamma\in Q^+_n$ and that
        $(\O,t)$ an $e$-idempotent subring of $\mathcal{K}$, where $e\ne2$.
        As an $\O$-algebra, the Iwahori-Hecke algebra $\HO[\gamma]$ is generated
        by the elements
        \[
        \set{\fo|\ei\in I^\gamma}\cup \set{\psip_r|1\leq r<n}\cup\set{\yp_s|1\leq s\leq n}
        \]
        subject to the relations
        {\setlength{\abovedisplayskip}{2pt}
        \setlength{\belowdisplayskip}{1pt}
        \begin{alignat*}{3}
          (\yp_1)^{(\Lambda_0,\alpha_{i_1})}\fo&=0,
          &\qquad \fo\fo[\bj]&= \delta_{\ei\ej}\fo,
          &\qquad\sum_{\ei\in I^\gamma}\fo&= 1\\
          \yp_r \fo&= \fo \yp_r,&\qquad \psip_r \fo&= f_{s_r\cdot\ei}^\O\psip_r,&\qquad
          \yp_r \yp_s&= \yp_s \yp_r
        \end{alignat*}
        \begin{alignat*}{2}
          \psip_r \yp_{r+1}\fo&= (\yp_r \psip_r+\delta_{i_r i_{r+1}})\fo,
          &\qquad\yp_{r+1}\psip_r \fo&= (\psip_r \yp_r+\delta_{i_r i_{r+1}})\fo
        \end{alignat*}
        \begin{align*}
          \psip_r \yp_s&= \yp_s\psip_r,   &\text{if $s\neq r,r+1$,}\\
          \psip_r\psip_s&= \psip_s\psip_r,&\text{if $\left|r-s\right|>1$,}
        \end{align*}
        \begin{align*}
          (\psip_r)^2\fo&= \begin{cases}
            (\dyo{1+\rho_r(\bi)}_r-\yp_{r+1})\fo,&\text{if }i_r\rightarrow i_{r+1},\\
            (\dyo{1-\rho_r(\bi)}_{r+1}-\yp_r)\fo,&\text{if }i_r\leftarrow i_{r+1},\\
            0,&\text{if }i_r=i_{r+1},\\
            \fo,&\text{otherwise,}
          \end{cases}\\
          \psip_r\psip_{r+1}\psip_r \fo&= \begin{cases}
            (\psip_{r+1}\psip_r\psip_{r+1}-t^{1+\rho_r(\bi)})\fo,
            &\text{if }i_r=i_{r+2}\rightarrow i_{r+1},\\
            (\psip_{r+1}\psip_r\psip_{r+1}+1)\fo,&\text{if }i_r=i_{r+2}\leftarrow i_{r+1},\\
            \psip_{r+1}\psip_r\psip_{r+1}\fo,&\text{otherwise,}
          \end{cases}
        \end{align*}
        }\noindent%
        where $\dyo{d}_r\fo=(t^d\yp_r-[d])\fo$ for $d\in \Z$, and
        $\rho_r(\bi)=\i_r-\i_{r+1}$, for $\ei,\ej\in I^\gamma$ and all
        admissible $r,s$.
      \end{Theorem}

      \begin{proof}
        By \cite[Theorem~A]{HuMathas:SeminormalQuiver}, this result holds for the
        algebras $\HO[\alpha]$, for $\alpha\in Q^+_n$.
        As~$\HO[\gamma]=\bigoplus_{\alpha\in\gamma}\HO[\alpha]$ this gives the
        result for~$\HO[\gamma]$. The proof of
        \cite[Theorem~A]{HuMathas:SeminormalQuiver} assumes that $e<\infty$ (or
        $e=0$ in the notation of \cite{HuMathas:SeminormalQuiver}), however,
        this assumption is only needed for
        \cite[(3.2)]{HuMathas:SeminormalQuiver} which is automatic in
        level~$1$ when $\Lambda=\Lambda_0$. Alternatively, as in
        \cite[Corollary~2.15]{HuMathas:SeminormalQuiver}, it is enough to
        consider the case when~$n<e<\infty$.
      \end{proof}

      As above, let $\m$ be a maximal ideal of $\O$ and set $\F=\O/\m$ and
    $\xi=t+\m\in\F$ and let~$\HSn$ be the Iwahori-Hecke algebra over~$\F$
    with parameter~$\xi$.  Then $\FHSn[\F]\cong\HO\otimes_O\F$. The definition of
    an $e$-idempotent subring ensures that~$\xi$ has quantum
    characteristic~$e$. Comparing the relations in \autoref{D:klr} with
    those in \autoref{T:HMGradedIso}, modulo~$\m$, there is an  algebra
    isomorphism $\theta\map{\FRSn[\F]}\FHSn[\F]$ determined by
    \[
        \psi_r\mapsto\psip_r\otimes1_\F,\quad
        y_r\mapsto\yp_r\otimes1_\F \quad\text{and}\quad
        e(\bi)\mapsto \fo\otimes1_\F,
    \]
    for all admissible $r$ and $\bi\in I^n$. Unfortunately, as the next
    example shows, we cannot use~$\theta$ to prove \autoref{T:Main}
    because $\theta\circ\sgn\ne\#\circ\theta$.

    \begin{Example}\label{Ex:BadBKRestriction}
      Suppose that $n=3$, $\Lambda=\Lambda_0$ and work over $\F_3$, the field with
      three elements. By \autoref{Ex:OS3},
      \[
      \set[\big]{e(012),  e(021),  y_3 e(012),  y_3e(021),  \psi_2e(012),
      \psi_2e(021)}
      \]
      is a basis of $\RSn[3]\cong\F_3\Sn[3]$ and, by \autoref{Ex:A3Basis},
      \[
      \set[\big]{e(012)+e(021),\psi_2(e(012)-e(021)), y_3(e(012)-e(021))}
      \]
      is a basis of $\RAn[3]$. Let $\theta\map{\RSn[3]}\F_3\Sn[3]$ be the
      Brundan-Kleshchev isomorphism induced by \autoref{T:HMGradedIso}. With
      some work it is possible to show that:
      \begin{align*}
        \theta\big(e(012)+e(021)\big)&= 1\\
        \theta\big(y_3(e(012)-e(021))\big)&= 1+s_1s_2+s_2s_1\\
        \theta\big(\psi_2(e(012)-e(021))\big)&= s_2+2s_1s_2s_1.
      \end{align*}
      In particular, $\theta$ does not restrict to an isomorphism
      between $\RAn[3]$ and $\F_3\An[3]$.
    \end{Example}

    To obtain an isomorphism $\RSn\to\HSn$, which restricts to an
    isomorphism $\RAn\to\HAn$, we modify the generators of~$\HO$
    given in \autoref{T:HMGradedIso}.

    \begin{Definition}\label{D:psim}
      Let $\psim_r=(\psip_r)^\#$ and $\ym_s=(\yp_s)^\#$, for
      $1\le r<n$ and $1\le s\le n$.
    \end{Definition}

    Notice that $(\fo)^\#=\fo[-\ei]$ by \autoref{L:foIdempotents}.
    Therefore, since $\#$ is an automorphism, the elements
    $\set{\psim_r,\ym_s,\fo}$ generate $\HO$, subject to essentially the
    same relations as those given in \autoref{T:HMGradedIso} except that
    $\bi$ should be replaced with $-\bi$. As we need this result below we
    state it in full for easy reference.

    \begin{Corollary}\label{C:HMGradedIso}
      Suppose that $e>2$, $\gamma\in Q^+_n$ and  $n\ge0$. Let $(\O,t)$ be an $e$-idempotent subring of~$\K$. Then $\HO[\gamma]$ is generated as an $\O$-algebra by the elements
      \[
      \set{\psim_r|1\leq r<n}\cup\set{\ym_s|1\leq s\leq n}\cup \set{\fo|\ei\in I^\gamma}
      \]
      subject to the relations
      {\setlength{\abovedisplayskip}{2pt}
      \setlength{\belowdisplayskip}{1pt}
      \begin{alignat*}{3}
        (\ym_1)^{(\Lambda_0,\alpha_{i_1})}\fo&=0,
        &\qquad \fo\fo[\bj]&= \delta_{\ei\ej}\fo,
        &\qquad\sum_{\ei\in I^\gamma}\fo&= 1\\
        \ym_r \fo&= \fo \ym_r,&\qquad \psim_r \fo&= f_{s_r\cdot\ei}^\O\psim_r,&\qquad
        \ym_r \ym_s&= \ym_s \ym_r
      \end{alignat*}
      \begin{alignat*}{2}
        \psim_r \ym_{r+1}\fo&= (\ym_r \psim_r+\delta_{i_r i_{r+1}})\fo,
        &\qquad\ym_{r+1}\psim_r \fo&= (\psim_r \ym_r+\delta_{i_r i_{r+1}})\fo
      \end{alignat*}
      \begin{align*}
        \psim_r \ym_s&= \ym_s\psim_r,   &\text{if $s\neq r,r+1$,}\\
        \psim_r\psim_s&= \psim_s\psim_r,&\text{if $\left|r-s\right|>1$,}
      \end{align*}
      \begin{align*}
        (\psim_r)^2\fo&= \begin{cases}
          (\dyo{1-\rho_r(\bi)}_r-\ym_{r+1})\fo,&\text{if }i_r\leftarrow i_{r+1},\\
          (\dyo{1+\rho_r(\bi)}_{r+1}-\ym_r)\fo,&\text{if }i_r\rightarrow i_{r+1},\\
          0,&\text{if }i_r=i_{r+1},\\
          \fo,&\text{otherwise,}
        \end{cases}\\
        \psim_r\psim_{r+1}\psim_r \fo&= \begin{cases}
          (\psim_{r+1}\psim_r\psim_{r+1}-t^{1-\rho_r(\bi)})\fo,
          &\text{if }i_r=i_{r+2}\leftarrow i_{r+1},\\
          (\psim_{r+1}\psim_r\psim_{r+1}+1)\fo,&\text{if }i_r=i_{r+2}\rightarrow i_{r+1},\\
          \psim_{r+1}\psim_r\psim_{r+1}\fo,&\text{otherwise,}
        \end{cases}
      \end{align*}
      }\noindent%
      where $\dyo{d}_r\fo=(t^d\ym_r-[d])\fo$, for all $d\in \Z$, for
      $\ei,\ej\in I^\gamma$ and all admissible~$r$ and~$s$.
    \end{Corollary}

    The elements $\dyo{d}_r$ appearing in \autoref{T:HMGradedIso} and
    \autoref{C:HMGradedIso} are different. In the next section we
    introduce a third variation of this notation. The meaning
    of~$\dyo{d}_r$ will always be clear from context.

    \subsection{Signed KLR generators}\label{S:AlternaatingCoefficients}

    This section sets up the machinery that will be used to construct the
    isomorphism $\RAn\to\HAn$. The idea is to
    use the results of the last two sections to give a new presentation
    of~$\HO$, which induces an isomorphism $\RSn\to\HSn$ that restricts
    to an isomorphism $\RAn\to\HAn$. To do this we use the generators of~$\HO$
    given in \autoref{T:HMGradedIso} and \autoref{C:HMGradedIso} together
    with a particular seminormal coefficient system and idempotent
    subring.

    The seminormal coefficient system that we use to prove
    \autoref{T:Main} forces us to work over a ring that contains
    ``enough'' square roots. We start by defining this ring, following
    \cite[Definition~3.1]{MathasRatliff}. Suppose that the
    Iwahori-Hecke algebra $\H$ is defined over the field~$F$ with
    parameter $\xi\in F$ of quantum characteristic~$e$. Recall that
    in this paper we are assuming that characteristic of~$F$
    is not~$2$ and that $e>2$.

    \begin{Definition}\label{D:squareroots}
      Let $x$ be an indeterminate over~$F$ and set $t=x+\xi$. In the
      algebraic closure $\overline{F(x)}$ of $F(x)$ fix square roots
      $\sqrt{-1}$, $\sqrt{t}$ and $\sqrt{[h]}$, for $1<h\le n$. Let
      \[\O=
      F\bigr[\sqrt{t},\sqrt{[h]}\ \big|\ 1<h\le n\bigr]_{(x)}
      \]
      be the localization of $F[\sqrt{t},\sqrt{[h]}\mid1< h\le n]$ at
      the maximal ideal generated by~$x$.  Let $\K$ be the field
      of fractions of $\O$.
    \end{Definition}

    Note that $t$ is invertible in~$\O$ so that we can consider the
    Iwahori-Hecke algebra~$\HO$ with parameter~$t$. By
    \cite[Corollary~5.12]{MathasRatliff}, the field of fractions~$\K$
    of~$\O$ is a splitting field for the semisimple algebra $\KHAn$.

    Let $\m=x\O$ be the maximal ideal of~$\O$ and set $\F=\O/\m$.  Then
    $F$ is (isomorphic to) a subfield of~$\F$.
    Moreover,~$\xi$ is identified with the image of~$t$ under the natural
    map
    $\O\twoheadrightarrow\F$. Hence,
    $\HSn\otimes_F\F\cong\HO\otimes_\O\F$. (Note that working over~$\F$ does not
    change the representation theory of~$\HSn$ because any field is a
    splitting field for~$\HSn$ since it is a cellular
    algebra~\cite{GL,M:ULect}.) By construction, $\F$ contains square
    roots $\sqrt{-1}$, $\sqrt{\xi}$ and $\sqrt{[h]_\xi}$, for $-n\le h\le n$.
    In general,~$\F$ is a non-trivial extension of~$F$.

    For $0<h\le n$ fix a choice of ``negative'' square roots in~$\O$ by setting
    \begin{equation}\label{E:Squareroots}
      \sqrt{[-h]}=\sqrt{-1}(\sqrt{t})^{-h}\sqrt{[h]}.
    \end{equation}
    Then $\sqrt{[-h]}\in\O$ for $-n\le h\le n$.
    If $h>0$ then $[-h]=-t^{-h}[h]$ so the effect of~\autoref{E:Squareroots} is
    to fix the sign of the square root of~$[-h]$.

    In order to apply the results of \autoref{T:HMGradedIso}
    and \autoref{C:HMGradedIso} we need to check that $(\O,t)$ is an
    idempotent subring in the sense of \autoref{D:idempotentSub}. Part~(a)
    of \autoref{D:idempotentSub} is automatic whereas parts~(b) and~(c)
    follow from the observation that if~$k\in\Z$ then the
    polynomial $[k]=[k]_t\in\O$ has zero constant term, as a polynomial
    in~$x$, if and only if $k\in e\Z$. Hence, we have the following.

    \begin{Lemma}\label{L:IdempotentSubring}
      The pair $(\O,t)$ is an idempotent subring.
    \end{Lemma}

    Now that we have fixed an idempotent subring we turn to the proof of
    \autoref{T:Main}. The idea is to use the generators of $\HO[\gamma]$
    from \autoref{T:HMGradedIso} for ``half'' of~$\HO[\gamma]$ and to use
    the generators from \autoref{C:HMGradedIso} the rest of the time.  To
    make this more
    precise, recall from after \autoref{L:ZeroSequence} that
    \[
    I^n_+=\set{\bi\in I_n|i_1=0\text{ and } i_2=1}
    \quad\text{ and }\quad
    I^n_-=\set{\bi\in I_n|i_1=0\text{ and } i_2=-1}.
    \]
    These sets are disjoint because $e\neq 2$. Set
    $I^\gamma_+=I^\gamma\cap I^n_+$ and $I^\gamma_-=I^\gamma\cap I^n_-$,
    for $\gamma\in\Qpm$. Then the map
    $\bi\mapsto-\bi$ is a bijection of sets $I^\gamma_+\bijection I^\gamma_-$.

    \begin{Lemma}\label{L:res}
      Suppose that $\bi\in I^\gamma$ and $\fo\ne0$. Then $\bi\in I^\gamma_+$ or
      $\bi\in I^\gamma_-$.
    \end{Lemma}

    \begin{proof}
      By definition, $\fo\ne0$ only if $\Std(\bi)\ne\emptyset$ or,
      equivalently, $\bi=\res(\t)$ for some standard tableau $\t\in\Std(\Parts)$.
      If $\t\in\Std(\bi)$ then $i_1=\res_1(\t)=0$ and $i_2=\res_2(\t)=\pm1$, so
      $\bi\in I^\gamma_+\cup I^\gamma_-$.
    \end{proof}

    In particular, if $h\in\HO[\gamma]$ then $h=\sum_{\bi\in
    I^\gamma_+}(h\fo+h\fo[-\bi])$. In what follows we apply
    \autoref{L:res}, and this observation, without further mention.

    Motivated in part by \autoref{P:AltCS} we make the following
    definition.

    \begin{Definition}\label{D:AltCS}
      An \textbf{alternating coefficient system} is a $*$-seminormal
      coefficient system $\balpha=\set{\alpha_r(\t)}$ such that
      $\alpha_r(\t)=-\alpha_r(\t')$, for $1\leq r<n$ and
      $\t\in\Std(\Parts)$.
    \end{Definition}

    Consider the case when $n=3$ and $\balpha$ is an alternating coefficient system.
    Let $\t=\Tableau{{1,2},{3}}$ and $\s=\Tableau{{1,3},{2}}$, so that
    $\res(\t)\in I^3_+$. By \autoref{D:AltCS} and \autoref{D:alphaSNCS},
    \[
    \alpha_2(\t)^2 = - \alpha_2(\s)\alpha_2(\t)=-\frac{t[3]}{\hspace*{1em}[2]^2}.
    \]
    So $\alpha_2(\t)=\pm\sqrt{-1}\sqrt{t}\sqrt{[3]}/[2]$. In the
    argument that follows we need  $\alpha_2(\t)$ to have such
    values for all $\t\in\Std(\Parts)$. Following
    \cite[\S3]{MathasRatliff}, for $\bi\in I^n$ and $1\le r<n$ define
    \begin{equation}\label{E:AlternatingCS}
      \alpha_r(\t)=\begin{cases}
        \frac{t^{\rho_r(\t)/2}\sqrt{[1+\rho_r(\t)]}\sqrt{[1-\rho_r(\t)]}}{[\rho_r(\t)]},
        &\text{if }\bi\in I^n_+,\\[2mm]
        -\alpha_r(\t'), &\text{if }\bi\in I^n_-,\\
        0,&\text{otherwise}.
      \end{cases}
    \end{equation}
    By \autoref{D:squareroots}, $\alpha_r(\t)\in\O$
    for all $\t\in\Std(\Parts)$ and $1\le r<n$.
    Moreover, if~$\bi\in I^n_\pm$ and $\t\in\Std(\bi)$ then
    $\alpha_2(\t)=\pm\sqrt{-1}\sqrt{t}\sqrt{[3]}/[2]$ by
    \autoref{E:Squareroots}.

    It is straightforward to check that  $\{\alpha_r(\t)\}$ is an
    alternating coefficient system. In particular, if $\t$ is standard,
    $1\le r<n$ and $s_r\t$ is not standard then
    $\rho_r(\t)=\pm1$ so that $\alpha_r(\t)=0$.

    Using \autoref{T:SeminormalForm}, we fix an arbitrary seminormal basis
    $\set{f_{\s\t}}$ for $\KHSn$ that is compatible with the
    seminormal coefficient system defined by \autoref{E:AlternatingCS}.
    Note that \autoref{D:psio}, and hence the results that follow, do
    not depend on this choice of seminormal basis.

    \begin{Definition}\label{D:psio}
      Suppose that $1\le r<n$ and  $1\le s\le n$. If $r\ne2$ define
      \[
      \psio_r=\sum_{\bi\in I^\gamma_+}(\psip_r \fo-\psim_r\fo[-\bi])
      \quad\text{and}\quad
      \yo_s=\sum_{\bi\in I^\gamma_+}(\yp_s\fo-\ym_s\fo[-\bi]),
      \]
      and when $r=2$ set
      $\psio_2=\Sum_{\bi\in I^\gamma_+}\kappa_\bi(\psip_2\fo-\psim_2\fo[-\bi])$, where
      \[ \kappa_\bi = \begin{cases}
        t^{-1},&\text{if }e=3,\\
        \frac{\sqrt{t}}{\sqrt{[3]}},&\text{if }e>3.
      \end{cases}
      \]
    \end{Definition}

    For convenience, set $\kappa_{-\bi}=\kappa_\bi$, for $\bi\in
    I^\gamma_+$.  The scalars $\kappa_\bi$ are needed to ensure that
    $\psio_2$ satisfies analogues of the quadratic and braid relations in
    \autoref{D:klr}.

    By \autoref{D:AltCS}, $\kappa_\bi$ is invertible in~$\O$ for all
    $\bi\in I^\gamma$. The reason why $\kappa_\bi$ depends only on~$e$,
    and not on the quiver~$\Gamma_e$, goes back to \autoref{L:ZeroSequence}: if
    $\bi\in I^\gamma$ and $e(\bi)\ne0$ or, equivalently, $\fo\ne0$ then
    the possible values for $i_1$, $i_2$ and $i_3$ are tightly constrained.

    By \autoref{D:squareroots}, the elements in $\set{\psio_r|1\leq
    r<n}\cup\set{\yo_s|1\leq s\leq n}$ belong to~$\HO[\gamma]$.  The aim
    is now to show that these elements, together with the idempotents
    $\set{\fo|\bi\in I^n}$, generate~$\HO[\gamma]$ subject to relations
    that are similar to those in \autoref{T:HMGradedIso}.  This will
    imply that these elements induce an isomorphism
    $\RSn\otimes_\Z\F\cong\HSn\otimes_F\F$. Before we start the proof we
    note the following consequence of \autoref{D:psio} and
    \autoref{L:foIdempotents}. Ultimately, this observation will imply that
    $\HAn\otimes_F\F\cong\RAn\otimes_\Z\F$.

    \begin{Corollary}\label{C:PsioHash}
      Suppose that $1\le r<n$, $1\le s\le n$ and $\bi\in I^\gamma$. Then
      \[(\psio_r)^\#=-\psio_r,\qquad
      (\yo_s)^\#=-\yo_s\qquad\text{and}\qquad
      (\fo)^\# = \fo[-\bi]. \]
    \end{Corollary}

    The first step is to give a new generating set for $\HO[\gamma]$.

    \begin{Proposition}\label{P:NewGenerators}
      Suppose $\gamma\in \Qpm$. Then $\HO[\gamma]$ is generated by
      \[
      \set{\psio_r|1\le r<n}\cup\set{\yo_s|1\le s\le n}\cup\set{\fo|\bi\in I^\gamma}.
      \]
    \end{Proposition}

    \begin{proof}Let $H_\gamma$ be the $\O$-subalgebra of $\HO[\gamma]$
      generated by the elements in the statement of the proposition.
      It is enough to show that $T_r\fo\in H_\gamma$,
      for $1\le r<n$ and~$\bi\in I^\gamma$, since these elements generate
      $\HO[\gamma]$. Further, by \autoref{C:PsioHash}, $H_\gamma^\#=H_\gamma$
      so it is enough to show that $T_r\fo\in H_\gamma$, for $\bi\in I^\gamma_+$ and
      $1\le r<n$. Let $\fo[+]=\sum_{\bi\in I^\gamma_+}\fo$. As remarked
      above, $\kappa_\bi$ is an invertible scalar in~$\O$. Therefore, the
      $\O$-module $H_\gamma\fo[+]$ contains the elements
      $\set{\psip_r \fo, \yp_s\fo, \fo|\bi\in I^\gamma_+}$. Hence,
      $H_\gamma\fo[+]=\HO\fo[+]$ by \autoref{T:HMGradedIso}. This completes the proof.
    \end{proof}

    For the rest of this paper, for $d\in\Z$, $1\le r\le n$ and $\bi\in I^\gamma$ we
    set
    \begin{equation}\label{E:YShiftDef}
        \dyo{d}_r\fo=\begin{dcases*}
                       (t^d\yo_r-[d])\fo,& if $\bi\in I^\gamma_+$,\\
                       (t^d\yo_r+[d])\fo,& if $\bi\in I^\gamma_-$.
                    \end{dcases*}
    \end{equation}
    Since $\yo_r=\sum_{\bi\in I^\gamma_+}(\yp_r\fo-\ym_r\fo[-\bi])$
    this is compatible with the two definitions of~$\dyo{d}_r\fo$ used in the
    last section

    The rest of this section determines a set of defining relations
    for $\HO[\gamma]$ for the generators of~$\HO[\gamma]$ from
    \autoref{P:NewGenerators}.  Fortunately, much of the work has
    already been done because \autoref{T:HMGradedIso} and
    \autoref{C:HMGradedIso} give us a large number of relations.  More
    precisely, they give the following list of relations,
    \textit{not} involving~$\psi_2^\O$.

    \begin{Lemma}\label{L:AutomaticRelations}
      Suppose that  $\gamma\in\Qpm$. The following
      identities hold in $\HO[\gamma]$:
      {\setlength{\abovedisplayskip}{2pt}
      \setlength{\belowdisplayskip}{1pt}
      \begin{alignat*}{3}
        (\yo_1)^{(\Lambda_0,\alpha_{i_1})}\fo&=0,
        &\qquad \fo\fo[\bj]&= \delta_{\ei\ej}\fo,&\qquad\textstyle\sum_{\ei\in I^\gamma}\fo&= 1\\
        \yo_t \fo&= \fo \yo_t,&\qquad \psio_r \fo&= \fo[s_r\cdot\ei]\psio_r,&\qquad
        \yo_r \yo_t&= \yo_t \yo_r
      \end{alignat*}
      \begin{alignat*}{2}
        \psio_r \yo_{r+1}\fo&= (\yo_r \psio_r+\delta_{i_r i_{r+1}})\fo,
        &\qquad\yo_{r+1}\psio_r \fo&= (\psio_r \yo_r+\delta_{i_r i_{r+1}})\fo
      \end{alignat*}
      \begin{align*}
        \psio_r \yo_t&= \yo_t\psio_r,   &\text{if $t\neq r,r+1$,}\\
        \psio_r\psio_s&= \psio_s\psio_r,&\text{if $\left|r-s\right|>1$,}
      \end{align*}
      \begin{align*}
        (\psio_r)^2\fo&= \begin{cases}
          (\dyo{1+\rho_r(\bi)}_r-\yo_{r+1})\fo,
          &\text{if $i_r\to i_{r+1}$ and $\bi\in I^\gamma_+$},\\
          (\yo_{r+1}-\dyo{1-\rho_r(\bi)}_r)\fo,
          &\text{if $i_r\leftarrow i_{r+1}$ and $\bi\in I^\gamma_-$}\\
          (\dyo{1-\rho_r(\bi)}_{r+1}-\yo_r)\fo,
          &\text{if $i_r\leftarrow i_{r+1}$ and $\bi\in I^\gamma_+$}\\
          (\yo_r-\dyo{1+\rho_r(\bi)}_{r+1})\fo,
          &\text{if $i_r\to i_{r+1}$ and $\bi\in I^\gamma_-$},\\
          0,&\text{if }i_r=i_{r+1},\\
          \fo,&\text{otherwise,}
        \end{cases}\\
        \psio_r\psio_{r+1}\psio_r \fo&= \begin{cases}
          (\psio_{r+1}\psio_r\psio_{r+1}-t^{1+\rho_r(\bi)})\fo,
          &\text{if $i_r=i_{r+2}\to i_{r+1},$ and $\bi\in I^\gamma_+$}\\
          (\psio_{r+1}\psio_r\psio_{r+1}+t^{1-\rho_r(\bi)})\fo,
          &\text{if $i_r=i_{r+2}\leftarrow i_{r+1}$ and $\bi\in I^\gamma_-$},\\
          (\psio_{r+1}\psio_r\psio_{r+1}+1)\fo,
          &\text{if $i_r=i_{r+2}\leftarrow i_{r+1}$ and $\bi\in I^\gamma_+$},\\
          (\psio_{r+1}\psio_r\psio_{r+1}-1)\fo,
          &\text{if $i_r=i_{r+2}\to i_{r+1},$ and $\bi\in I^\gamma_-$}\\
          \psio_{r+1}\psio_r\psio_{r+1}\fo,&\text{otherwise,}
        \end{cases}
      \end{align*}
      }\noindent%
      for all admissible $\bi,\bj\in I^\gamma$ and $r,s,t$ satisfying
      $2<r,s<n$ and $1\le t\le n$.
    \end{Lemma}

    \begin{proof}
      First notice that if $\bi\notin I^\gamma_+\cup I^\gamma_-$ then $\fo=0$ by
      \autoref{L:res}, so all of the relations above are trivially true.
      We may assume then that $\bi\in I^\gamma_+\cup I^\gamma_-$.

      The first three identities follow directly from \autoref{T:HMGradedIso} and
      \autoref{C:HMGradedIso}. For the remaining formulas, observe that if
      $2<r<n$ then $\bi\in I^\gamma_+$ if and only if $s_r\cdot\bi\in I^\gamma_+$ and,
      similarly, $\bi\in I^\gamma_-$ if and only if $s_r\cdot\bi\in I^\gamma_-$.
      Therefore, if $\bi\in I^\gamma_+$ the relations hold by virtue of
      \autoref{T:HMGradedIso} and if $\bi\in I^\gamma_-$ then they hold by
      \autoref{C:HMGradedIso}. Note that if~$\bi\in I^\gamma_-$ then there is a sign
      change in the last two relations, in comparison with \autoref{C:HMGradedIso},
      because $\psio_r\fo=-\psim_r\fo$ and $\yo_t\fo=-\ym_s\fo$.
    \end{proof}

    Next we need analogues of the relations in
    \autoref{L:AutomaticRelations} for~$\psio_1$. We could replace the
    next result with the single relation~$\psio_1=0$, however, this
    is not sufficient for our later arguments because the proof of
    \autoref{T:MainTheorem} relies on the fact that the generators
    of~\autoref{P:NewGenerators} satisfy relations that are compatible
    with \autoref{D:klr}.

    \begin{Lemma}\label{L:PsiOne}
      Suppose that $\gamma\in\Qpm$. The following
      identities hold in $\HO[\gamma]$:
      {\setlength{\abovedisplayskip}{2pt}
      \setlength{\belowdisplayskip}{1pt}
      \begin{alignat*}{3}
        \psio_1\fo&=\fo[s_1\cdot\bi]\psio_1,
        &\qquad \psio_1 \yo_s&= \yo_s\psio_1,
        &\qquad \psio_1\psio_r&= \psio_r\psio_1,
      \end{alignat*}
      \begin{alignat*}{2}
        \psio_1 \yo_2\fo&= (\yo_1 \psio_1+\delta_{i_1 i_2})\fo,
        &\qquad\yo_2\psio_1 \fo&= (\psio_1 \yo_1+\delta_{i_1 i_2})\fo,
      \end{alignat*}
      \begin{align*}
        (\psio_1)^2\fo&= \begin{cases}
          (\yo_1-\yo_2)\fo,
          &\text{if }i_1\to i_2,\\
          (\yo_2-\yo_1)\fo,
          &\text{if }i_1\leftarrow i_2,\\
          0,&\text{if }i_1=i_2,\\
          \fo,&\text{otherwise,}
        \end{cases}\\
        \psio_1\psio_2\psio_1 \fo&= \begin{cases}
          (\psio_2\psio_1\psio_2-1)\fo,
          &\text{if }i_1=i_3\to i_2,\\
          (\psio_2\psio_1\psio_2+1)\fo,
          &\text{if }i_1=i_3\leftarrow i_2,\\
          \psio_2\psio_1\psio_2\fo,&\text{otherwise,}
        \end{cases}
      \end{align*}
      }\noindent%
      for all admissible $\bi\in I^\gamma$ and $r,s$ satisfying
      $2<r<n$ and $3\le s\le n$.
    \end{Lemma}

    \begin{proof}
      By definition, $\fo\ne0$ if and only if~$\bi=\res(\s)$ for some
      standard tableau~$\s$. In particular, if $\bi=(i_1,\dots,i_n)$ then
      $i_1=0$, $i_2\in\set{-1,1}$ and $i_3\in\set{-2,-1,1,2}$. Hence, it
      follows from \autoref{T:HMGradedIso} and \autoref{C:HMGradedIso}
      that $\psip_1=0=\psim_1$. Therefore, $\psio_1=0$ and the first three
      relations are trivially true. The next two relations hold
      because $\delta_{i_1i_2}=0$ whenever $\fo\ne0$ and the quadratic
      relation for $(\psio_1)^2$ holds in view of \autoref{T:HMGradedIso}
      and \autoref{C:HMGradedIso}. For the final ``braid'' relation,
      if~$i_1=i_3\rightarrow i_2$ or $i_1=i_3\leftarrow i_2$ then
      $\fo=0$ by the remarks at the start of the proof, so the
      braid relation is trivially true in these cases. In the remaining cases
      $\psio_1\psio_2\psio_1=0=\psio_2\psio_1\psio_2$ since
      $\psio_1=0$. This completes the proof.
    \end{proof}

    It remains to determine the relations involving $\psio_2$. The first step is
    easy.

    \begin{Lemma}\label{L:SimplePsi2Relations}
      Suppose that $\bi\in I^\gamma$, $2<r<n$ and $1\le s\le n$ with
      $t\ne 2,3$. Then
      \[
            \psio_2\psio_r= \psio_r\psio_2
                \qquad\text{and}\qquad
            \psio_2 \yo_s= \yo_s\psio_2.
      \]
    \end{Lemma}

    \begin{proof}
      Since $\psio_2\fo=\pm\kappa_\bi\psi^\pm_2\fo$, where $\kappa_\bi\in\O$ is
      invertible for $\bi\in I^\gamma$, the result follows directly from
      \autoref{T:HMGradedIso} and \autoref{C:HMGradedIso}.
    \end{proof}

    For the remaining relations we need a more precise description of how the
    generators of~\autoref{D:psio} act on the seminormal basis. Suppose that
    $\s\in\Std(\bi)$ and that $\u=(r,r+1)\s$, where $1\le r<n$.
    Following \cite[(4.21)]{HuMathas:SeminormalQuiver}, define
    \begin{equation}\label{E:beta}
      \beta_r(\s)=\begin{cases}
        -\beta_r(\s'),&\text{if }\bi\in I^\gamma_-,\\[1pt]
        \dfrac{t^{\i_{r}-c_{r}(\s)}\alpha_r(\s)}{[1-\rho_r(\s)]},
        &\text{if $\bi\in I^\gamma_+$ and }i_r=i_{r+1},\\[1pt]
        t^{c_{r+1}(\s)-\i_r}\alpha_r(\s)[\rho_r(\s)],
        &\text{if $\bi\in I^\gamma_+$ and }i_r=i_{r+1}+1,\\
        \dfrac{t^{-\rho_r(\s)}\alpha_r(\s)[\rho_r(\s)]}{[1-\rho_r(\s)]},
        &\text{if $\bi\in I^\gamma_+$ and }i_r\notin\set{i_{r+1}, i_{r+1}+1}.
      \end{cases}
    \end{equation}
    Note that $\beta_r(\s)=0$ if $\u$ is not standard because
    $\alpha_r(\u)=0$ whenever $\u\notin\Std(\Parts)$. More explicit formulas
    for $\beta_r(\s)$ can be obtained using
    \autoref{E:AlternatingCS}, however, we will only need these in one special case;
    see \autoref{E:BetaTwo} below.

    Note that if $\s$ is standard and $\bi=\res(\s)$ then $\res(\s')=-\bi$. Therefore,
    the four cases in \autoref{E:beta} are mutually exclusive.
    We need to be slightly careful, however, because
    if $\bi=\res(\s)$ and $\bj=\res(\s')$ then it is
    not usually true that $\j_r=-\i_r$, for $1\le r\le n$.

    Following \cite[Lemma~4.23]{HuMathas:SeminormalQuiver} we can now
    describe the action of $\psio_r$ and $\yo_r$ on the seminormal basis
    $\set{f_{\s\t}}$. This result is the only place where we explicitly use
    the assumption of \autoref{D:AltCS}.

    \begin{Proposition}\label{P:betas}
      Suppose that $\s,\t\in\Std(\lambda)$ for $\lambda\in\Parts$ and let
      $\bi=\res(\s)$, $\bj=\res(\s')$ for $\bi,\bj\in I^\gamma$. Fix $1\leq r<n$
      and let $\u=(r,r+1)\s$. Then
      \[
      \psio_r f_{\s\t} =\begin{cases}
        \kappa_\bi\beta_2(\s)f_{\u\t},&\text{if }r=2,\\
        \beta_r(\s)f_{\u\t}-\delta_{i_ri_{r+1}}
        \dfrac{t^{\i_{r+1}-c_{r+1}(\s)}}{[\rho_r(\s)]}f_{\s\t},
        &\text{if $r\ne2$ and }\bi\in I^\gamma_+,\\[3mm]
        \beta_r(\s)f_{\u\t}-\delta_{i_ri_{r+1}}
        \dfrac{t^{\j_{r+1}-c_{r+1}(\s)}}{[\rho_r(\s)]}f_{\s\t},
        &\text{if $r\ne2$ and }\bi\in I^\gamma_-.\\
      \end{cases}
      \]
      Moreover, if $1\le k\le n$ then
      \[\yo_k f_{\s\t}=\begin{cases}
        \phantom{-}[c_k(\s)-\i_k]f_{\s\t},&\text{if }\bi\in I^\gamma_+,\\
        -[c_k(\s')-\j_k]f_{\s\t},
        &\text{if }\bi\in I^\gamma_-.\\
      \end{cases}\]
    \end{Proposition}

    \begin{proof}
      Without loss of generality, we can
      assume that $\t=\s$ by \autoref{T:SeminormalForm}(a), so we need to
      compute $\psio_rf_{\s\s}$ and $\yo_r f_{\s\s}$.

      First consider $\psio_rf_{\s\s}$ when $r\ne2$. If $\bi\in I^\gamma_+$ then
      $\psio_r f_{\s\s}=\psip_rf_{\s\s}$ and the lemma is a restatement of
      \cite[Lemma~4.23]{HuMathas:SeminormalQuiver}. Suppose then that
      $\bi\in I^\gamma_-$,
      so that $\bj\in I^\gamma_+$ and $\psio_r f_{\s'\s'}$ is given
      by the formulas above. As $\#$ is an involution,
      using \autoref{C:fttHash} for the third equality and \autoref{L:futhash}
      for the last equality,
      \begin{align*}
        \psio_r f_{\s\s} &=-(\psip_r)^\#f_{\s\s}
        =-\big(\psip_r f_{\s\s}^\#\big)^\#
        =-\frac{\gamma_\s}{\gamma_{\s'}}\big(\psip_r f_{\s'\s'}\big)^\#,\\
        &=-\frac{\gamma_\s}{\gamma_{\s'}}\Big(\beta_r(\s')f_{\u'\s'}
        -\delta_{j_rj_{r+1}}
        \frac{t^{\j_{r+1}-c_{r+1}(\s')}}{[\rho_r(\s')]}f_{\s'\s'}\Big)^\#\\
        &=\frac{\alpha_r(\s)\beta_r(\s')}{\alpha_r(\s')}f_{\u\s}
        -\delta_{i_ri_{r+1}}\frac{t^{\j_{r+1}-c_{r+1}(\s)}}{[\rho_r(\s)]}f_{\s\s},
      \end{align*}
      since $[\rho_r(\s')]=[-\rho_r(\s)]=-t^{-\rho_r(\s)}[\rho_r(\s)]$. By
      \autoref{D:AltCS}, $\alpha_r(\s')=-\alpha_r(\s)$ and
      $\beta_r(\s')=-\beta_r(\s)$, so this establishes the formula for
      $\psio_r f_{\s\t}$ when $r\ne2$.

      Now consider $\psio_2f_{\s\s}$. If $\fo\ne0$ then $i_2\ne i_3$
      because $\fo\ne0$ only if $\bi$ is
      the residue sequence of some standard tableau. Therefore,
      if $\bi\in I^\gamma_+$ and $\fo\ne0$ then
      the argument of the last paragraph shows that
      $\psip_2f_{\s\t}=\beta_2(\s)f_{\u\t}$ and if $\bi\in I^\gamma_-$ then
      $-\psim_2f_{\s\t}=\beta_2(\s)f_{\u\t}$. As
      $\psio_2=\sum_{\bi\in I^\gamma_+}\kappa_\bi(\psip_2\fo-\psim_2\fo[-\bi])$,
      it follows that $\psio_2 f_{\s\t}=\kappa_\bi\beta_2(\s)f_{\u\t}$ as
      claimed.

      For the action of $\yo_k$, if $\bi\in I^\gamma_+$ then
      $\yo_k f_{\s\s}=\yp_k f_{\s\s}=[c_k(\s)-\i_k]f_{\s\s}$ by
      \cite[Lemma~4.23]{HuMathas:SeminormalQuiver}. On the other hand, if
      $\bi\in I^\gamma_-$ then, using \autoref{L:foIdempotents} twice,
      \[\yo_k f_{\s\s}
      = -\frac{\gamma_\s}{\gamma_{\s'}}(\yp_k f_{\s'\s'})^\#
      = -\frac{\gamma_\s}{\gamma_{\s'}}\big([c_k(\s')-\j_k] f_{\s'\s'}\big)^\#
      =-[c_k(\s')-\j_k]f_{\s\s}.
      \]
      as required.
    \end{proof}

    We can now determine the remaining ``KLR-like'' relations satisfied by
    $\psio_2$.

    \begin{Lemma}\label{L:psi2Intertwiner}
      Suppose that $\bi\in I^\gamma$ and let $\bj=s_2\cdot\bi$. Then
      $\psio_2\fo=\fo[\bj]\psio_2$.
    \end{Lemma}

    \begin{proof}
      If $\bi\notin I^\gamma_+\cup I^\gamma_-$ then $\fo=0=\fo[\bj]$ and
      there is nothing to prove. Therefore, we may assume that $\bi\in
      I^\gamma_+\cup I^\gamma_-$. Recall that
      $\fo[\gamma]=\sum_{\bk\in I^\gamma}\fo[\bk]$ is the identity element
      of~$\HO[\gamma]$.  By \autoref{P:betas},
      \[ \psio_2=\psio_2\fo[\gamma]
              =\sum_{\bk\in I^\gamma}\psio_2\fo[\bk]
              =\sum_{\substack{\bk\in I^\gamma\\\t\in\Std(\bk)}}
                 \frac1{\gamma_\t}\psio_2f_{\t\t}
              =\sum_{\substack{\bk\in I^\gamma\\\t\in\Std(\bk)\\\s=s_2\t}}
              \frac{\kappa_\bk\beta_2(\t)}{\gamma_\t}f_{\s\t}.
      \]
      If $f_{\s\t}$ is a term in the right-hand sum, with
      $\t\in\Std(\bk)$ and $\s=s_2\t$, then
      \[
           f_{\s\t}\fo[\bi]=\delta_{\bi\bk}f_{\s\t}
                           =\delta_{s_2\cdot\bi,s_2\cdot\bk}f_{\s\t}
                           =\fo[\bj]f_{\s\t},
      \]
      by \autoref{T:SeminormalForm}(b). Hence,
      $\psio_2\fo=\sum_{\t\in\Std(\bi)}\frac1{\gamma_\t}\kappa_\bi\beta_2(\t)f_{\s\t}
                 =\fo[\bj]\psio_2$,
      as required.
    \end{proof}

    Recall from \autoref{E:YShiftDef} that if $d\in\Z$ and
    $\bi\in I^\gamma_\pm$
    then $\dyo{d}_r\fo=(t^d\yo_r\mp[d])\fo$.

    \begin{Lemma}\label{L:Mixed}
      Suppose that $\gamma\in\Qpm$ and $\bi\in I^\gamma$. Then:
      \[
        \yo_3 \psio_2\fo=\big(\psio_2\dyo{-e}_2+\delta_{i_2i_3}\big)\fo\\
          \quad\text{and}\quad
        \psio_2 \yo_3\fo=\big(\dyo{-e}_2\psio_2+\delta_{i_2i_3}\big)\fo.
      \]
    \end{Lemma}

    \begin{proof}
      Both identities are proved similarly so we consider only the first
      one.  If~$\fo\ne0$ then $\bi=\res(\s)$, for some standard
      tableau~$\s$, in which case $i_2\ne i_3$. Hence, if~$\fo\ne0$ then
      $\delta_{i_2i_3}=0$ so we can assume that $\delta_{i_2i_3}=0$ in
      what follows.  (We include the term for $\delta_{i_2i_3}$
      because to prove \autoref{T:Main} we need to compare the
      identity in the lemma with the relations in \autoref{D:klr}.)
      Without loss of generality, we may assume that $\bi\in I^\gamma_+$.

      By \autoref{T:SeminormalForm}(b),
      $\fo=\sum_{\s\in\Std(\bi)}\frac1{\gamma_\s}f_{\s\s}$. Therefore,
      to prove the lemma it is enough to verify that
      $\yo_3 \psio_2f_{\s\s}=\psio_2\dyo{-e}_2f_{\s\s}$, for all
      $\s\in\Std(\bi)$.  Fix $\s\in\Std(\bi)$ and set $\u=(2,3)\s$ and
      $\bj=s_2\cdot\bi\in I^\gamma_-$ so that $\bj=\res(\u)$ if~$\u$
      is standard. If~$\u$ is not standard then $\beta_2(\s)=0$ so
      $\yo_3 \psio_2f_{\s\s}=0=\psio_2\dyo{-e}_2f_{\s\s}$ by \autoref{P:betas}.
      Suppose then that $\u$ is standard so that
      \[  \s_{\downarrow3}=\Tableau{{1,2},{3}} \quad\text{and}\quad
          \u_{\downarrow3}=\Tableau{{1,3},{2}}.
      \]
      Using \autoref{P:betas} again,
      $\psio_2\dyo{-e}_2f_{\s\s}=-\kappa_\bi\beta_2(\s)[-e]f_{\u\s}$
      since $\yo_2f_{\s\s}=0$. Similarly,
      $\yo_3\psio_2f_{\s\s}
               =\kappa_\bi\beta_2(\s)\yo_3f_{\u\s}
               =-\kappa_\bi\beta_2(\s)[-e]f_{\u\s}.
      $
      Hence, $\yo_3\psio_2f_{\s\s} =\psio_2\dyo{-e}_2f_{\s\s}$ in all cases,
      completing the proof.
    \end{proof}

    The proof of the next result explains why $\kappa_\bi$ is needed in the
    definition of~$\psio_2$. Fix $\s\in\Std(\bi)$ such that  $\bi\in I^\gamma_+$
    and $(2,3)\s$ is standard. Then $\rho_2(\s)=2$ and either
    $e=3$ and $i_2\rightarrow i_3$, or $e>3$ and $i_2\noedge i_3$.  Hence,
    $i_2\not\in\set{i_3,i_3+1}$, so by \autoref{E:beta} and \autoref{D:psio}
    \begin{equation}\label{E:BetaTwo}
      \beta_2(\s)
      =\frac{t^{-\rho_2(\s)}\sqrt{-1}\sqrt{t}\sqrt{[3]}[\rho_2(\s)]}%
      {[1-\rho_2(\s)][2]}
      =-\frac{\sqrt{-1}\,\sqrt{[3]}}{\sqrt{t}}.
    \end{equation}
    We can now determine the quadratic relation for $\psio_2$.

    \begin{Lemma}\label{L:Quadratic}
      Suppose that $\bi\in I^\gamma$. Then
      \[(\psio_2)^2\fo = \begin{dcases*}
        (\yo_2-\yo_3)\fo, &if $i_2\to i_3$,\\
        (\yo_3-\yo_2)\fo, &if $i_2\leftarrow i_3$,\\
        0,&if $i_2=i_3$,\\
        \fo,&otherwise.
      \end{dcases*}
      \]
    \end{Lemma}

    \begin{proof}
      It is enough to consider the case when $\bi\in I^\gamma_+$. Since
      $\fo=\sum_{\s\in\Std(\bi)}\frac1{\gamma_\s}f_{\s\s}$ we are
      reduced to computing $(\psio_2)^2f_{\s\s}$, for $\s\in\Std(\bi)$ and
      $\bi\in I^\gamma_+$. Fix~$\s\in\Std(\bi)$ and let
      $\u=(2,3)\s\in\Std(\bj)$. By \autoref{P:betas},
      \[
        (\psio_2)^2f_{\s\s} = \kappa_\bi\beta_2(\s)\psio_2f_{\u\s}
                            = -\kappa_\bi^2\beta_2(\s)^2f_{\s\s}.
      \]
      If $3$ is in the first row of~$\s$ then $\u$ is not standard so
      $\beta_2(\s)=0$ and $(\psio_2)^2f_{\s\s}=0$.
      In this case, $i_2\to i_3$ and $\yo_2 f_{\s\s}=0=\yo_3f_{\s\s}$, so the
      lemma holds. The only other possibility is that~$3$ is in the first
      column of~$\s$, so that $\rho_2(\s)=2$. Then $i_2\to i_3$ if~$e=3$ and
      $i_2\noedge i_3$ if~$e>3$. Hence, using \autoref{D:psio} and
      \autoref{E:BetaTwo},
      \[   (\psio_2)^2f_{\s\s}=
      -\kappa_\bi^2\beta_2(\s)^2f_{\s\s}=\begin{dcases*}
        t^{-3}[3]f_{\s\s},&if $i_2\to i_3$ (and $e=3$),\\
        f_{\s\s},&if $i_2\noedge i_3$ (and $e>3$).
      \end{dcases*}
      \]
      Hence, if $i_2\noedge i_3$ then $(\psio_2)^2\fo=\fo$ as claimed.
      Finally, if $i_2\to i_3$ then
      \[  (\yo_2-\yo_3)f_{\s\s}
                =(0-[-e])f_{\s\s}
                =t^{-3}[3]f_{\s\s}
                =(\psio_2)^2f_{\s\s},
      \]
      where the middle equality holds only because $e=3$. This completes the proof.
    \end{proof}

    \begin{Remark*}
      The proof of \autoref{L:Quadratic} suggests that $\kappa_\bi$ is
      uniquely determined, for $\bi\in I^\gamma$. In fact, this is not
      quite true. What the proof shows is that the value of $\kappa_\bi$
      is uniquely determined by the quadratic relation satisfied
      by~$\psio_2$. For the proof of our main results we only need
      $\psio_2$ to satisfy a ``deformed'' version of the quadratic
      relation for $\psi_2$ in \autoref{L:Quadratic}. For example, we
      can obtain slightly different relations by replacing~$\yo_r$
      with~$\dyo{ke}_r$, for some $k\in\Z$. Such relations would require
      a different value for~$\kappa_\bi$.
    \end{Remark*}

    Finally, it remains to check the braid relation for $\psio_2$ and~$\psio_3$.

    \begin{Lemma}
      Suppose that $\bi\in I^\gamma$. Then
      \begin{align*}
        \psio_2\psio_3\psio_2 \fo&= \begin{cases}
          (\psio_3\psio_2\psio_3-1)\fo, &\text{if }i_2=i_4\to i_3,\\
          (\psio_3\psio_2\psio_3+1)\fo, &\text{if }i_2=i_4\leftarrow i_3,\\
          \psio_3\psio_2\psio_3\fo,&\text{otherwise,}
        \end{cases}
      \end{align*}
    \end{Lemma}

    \begin{proof}Again, it is enough to consider the case when
      $\bi\in I^\gamma_+$. We fix $\s\in\Std(\bi)$ and show that the two
      sides of the identity in the lemma act in the same way on~$f_{\s\s}$.
      Let $\bj=(2,4)\cdot\bi\in I^\gamma$.
      By \autoref{L:psi2Intertwiner} and \autoref{L:SimplePsi2Relations},
      $\psio_2\psio_3\psio_2 \fo=\fo[\bj]\psio_2\psio_3\psio_2$,
      so  $\psio_2\psio_3\psio_2 \fo=0$ unless $\bj$ is the
      residue sequence of a standard tableau. Similarly
      $\psio_3\psio_2\psio_3 \fo=0$ unless $\bj$ is the
      residue sequence of a standard tableau. Let $\s_{\downarrow4}$ be the
      subtableau of~$\s$ containing the numbers $1,2,3,4$. Then
      \[
      \s_{\downarrow4}\in\set[\Bigg]{ \Tableau[-4]{{1,2,3,4}},  \quad
      \Tableau[-4]{{1,2,3},{4}},\quad
      \Tableau[-4]{{1,2,4},{3}},\quad
      \Tableau[-4]{{1,2},{3,4}},\quad
      \Tableau[-4]{{1,2},{3},{4}} }
      \]
      since $\bi\in I^\gamma_+$. We consider two cases.

      \smallskip\noindent\textbf{Case 1: $e>3$:}
      Inspecting the list of possibilities for $\s_{\downarrow4}$, in all cases
      $i_2\ne i_4$ and  $\bj\ne\res(\t)$ for any standard tableau~$\t$. Therefore,
      \[ \psio_2\psio_3\psio_2 \fo=0=\psio_3\psio_2\psio_3\fo,\]
      in agreement with the statement of the lemma.

      \smallskip\noindent\textbf{Case 2: $e=3$:}
      Except for the last tableau in the set above, $i_2\ne i_4$ and
      $\bj$ is not a residue sequence for a standard tableau. Hence, as
      in Case~1, the lemma holds when $i_2\ne i_4$ as both sides are
      zero. Moreover, if $\fo\ne0$ then the case $i_2=i_4\leftarrow i_3$
      does not arise, so the lemma is vacuously true in this case.  It
      remains to consider the case when $i_2=i_4\to i_3$, which occurs
      only if $s_{\downarrow4}$ is the last tableau in the set above.
      Noting that $\psio_3f_{\s\s}=0$, \autoref{P:betas} and
      \autoref{E:BetaTwo} quickly imply that
      \[ \big(\psio_2\psio_3\psio_2-\psio_3\psio_2\psio_3\big)f_{\s\s}
      =-\kappa_\bi^2\beta_2(\s)^2\frac{t^3}{[3]}f_{\s\s}
      =-f_{\s\s},
      \]
      where the last equality follows using \autoref{E:BetaTwo} exactly as
      in the proof of \autoref{L:Quadratic}. This completes the proof.
    \end{proof}

    \subsection{The isomorphism $\RAn\cong\HAn$}\label{S:Main}
    We now have almost everything in place that we need to prove
    \autoref{T:Main}. We first prove a stronger version of
    \autoref{T:BKiso} over~$\O$. For this we need the following
    definition, which should be viewed as an $\O$-deformation of~$\RSn$.
    The reader should compare this result with \autoref{T:HMGradedIso}.

    \begin{Definition}\label{D:RO}
      Suppose that $\gamma\in\Qpm$. Let $\RO$ be the unital
      associative $\O$-algebra generated by the elements
      \[
      \set{\wpsio_1,\wpsio_2,\ldots,\wpsio_{n-1}}
      \cup\set{\wyo_1,\wyo_2,\ldots,\wyo_n}
      \cup\set{\wfo| \bi\in I^\gamma}
      \]
      subject to the relations
      {\setlength{\abovedisplayskip}{2pt}
      \setlength{\belowdisplayskip}{1pt}
      \begin{alignat*}{3}
        (\wyo_1)^{(\Lambda_0,\alpha_{i_1})}\wfo&=0,
        &\wfo\wfo[\bj]&= \delta_{\ei\ej}\wfo,
        &\textstyle\sum_{\ei\in I^\gamma}\wfo&= 1,\\
        \wyo_t \wfo&= \wfo \wyo_t,&\qquad
        \wpsio_r \wfo&= \wfo[s_r\cdot\ei]\wpsio_r,&
        \wyo_r \wyo_t&= \wyo_t \wyo_r,
      \end{alignat*}
      \begin{align*}
        \wpsio_1 \wyo_2\wfo=\big(\wyo_1\wpsio_1+\delta_{i_1i_2}\big)\wfo,
        \qquad
        \wyo_2 \wpsio_1\wfo=\big(\wpsio_1\wyo_1+\delta_{i_1i_2}\big)\wfo,\\
        \wpsio_2 \wyo_3\wfo=\big(\ddyo{-e}_2\wpsio_2+\delta_{i_2i_3}\big)\wfo,
        \qquad
        \wyo_3 \wpsio_2\wfo=\big(\wpsio_2\ddyo{-e}_2+\delta_{i_2i_3}\big)\wfo,
      \end{align*}
      \begin{align*}
        \wpsio_r \wyo_t&= \wyo_t\wpsio_r,   \quad\text{if }t\neq r,r+1,\\
        \wpsio_r\wpsio_s&= \wpsio_s\wpsio_r,\quad\text{if }\left|r-s\right|>1,
      \end{align*}
      if $r=1$ or $r=2$ then
      \begin{align*}
        (\wpsio_r)^2\wfo&= \begin{cases}
          (\wyo_r-\wyo_{r+1})\wfo,
          &\text{if }i_r\rightarrow i_{r+1},\\
          (\wyo_{r+1}-\wyo_r)\wfo,
          &\text{if }i_r\leftarrow i_{r+1},\\
          0,&\text{if }i_r=i_{r+1},\\
          \wfo,&\text{otherwise,}
        \end{cases}
      \end{align*}
      \begin{align*}
        \wpsio_r\wpsio_{r+1}\wpsio_r \wfo&= \begin{cases}
          (\wpsio_{r+1}\wpsio_r\wpsio_{r+1}-1)\wfo,
          &\text{if }i_r=i_{r+2}\rightarrow i_{r+1},\\
          (\wpsio_{r+1}\wpsio_r\wpsio_{r+1}+1)\wfo,
          &\text{if }i_r=i_{r+2}\leftarrow i_{r+1},\\
          \wpsio_{r+1}\wpsio_r\wpsio_{r+1}\wfo,&\text{otherwise},
        \end{cases}
      \end{align*}
      and if $2<r<n$ then
      \begin{align*}
        \wpsio_r \wyo_{r+1}\wfo= (\wyo_r \wpsio_r+\delta_{i_r i_{r+1}})\wfo,
        \qquad
        \wyo_{r+1}\wpsio_r \wfo= (\wpsio_r \wyo_r+\delta_{i_r i_{r+1}})\wfo,
      \end{align*}
      \begin{align*}
        (\wpsio_r)^2\wfo&= \begin{cases}
          (\ddyo{1+\rho_r(\bi)}_r-\wyo_{r+1})\wfo,
          &\text{if $i_r\rightarrow i_{r+1}$ and $\bi\in I^\gamma_+$}\\
          (\wyo_{r+1}-\ddyo{1-\rho_r(\bi)}_r)\wfo,
          &\text{if $i_r\leftarrow i_{r+1}$ and $\bi\in I^\gamma_-$}\\
          (\ddyo{1-\rho_r(\bi)}_{r+1}-\wyo_{r})\wfo,
          &\text{if $i_r\leftarrow i_{r+1}$ and $\bi\in I^\gamma_+$}\\
          (\wyo_{r}-\ddyo{1+\rho_r(\bi)}_{r+1})\wfo,
          &\text{if $i_r\rightarrow i_{r+1}$ and $\bi\in I^\gamma_-$}\\
          0,&\text{if }i_r=i_{r+1},\\
          \wfo,&\text{otherwise,}
        \end{cases}
      \end{align*}
      \begin{align*}
        \wpsio_r\wpsio_{r+1}\wpsio_r \wfo&= \begin{cases}
          (\wpsio_{r+1}\wpsio_r\wpsio_{r+1}-t^{1+\rho_r(\ei)})\wfo,
            &\text{if $i_r=i_{r+2}\rightarrow i_{r+1}$ and $\bi\in I^\gamma_+$,}\\
          (\wpsio_{r+1}\wpsio_r\wpsio_{r+1}+t^{1-\rho_r(\ei)})\wfo,
            &\text{if $i_r=i_{r+2}\leftarrow i_{r+1}$ and $\bi\in I^\gamma_-$,}\\
          (\wpsio_{r+1}\wpsio_r\wpsio_{r+1}-1)\wfo,
            &\text{if $i_r=i_{r+2}\rightarrow i_{r+1}$ and $\bi\in I^\gamma_-$,}\\
          (\wpsio_{r+1}\wpsio_r\wpsio_{r+1}+1)\wfo,
            &\text{if $i_r=i_{r+2}\leftarrow i_{r+1}$ and $\bi\in I^\gamma_+$,}\\
          \wpsio_{r+1}\wpsio_r\wpsio_{r+1}\wfo,&\text{otherwise,}
        \end{cases}
      \end{align*}
      }\noindent%
      for all admissible $\bi,\bj\in I^\gamma$ and all admissible $r,s$ and $t$
      and where for $d\in\Z$
      \[   \ddyo{d}_r\wfo=\begin{dcases*}
                    (t^d\wyo_r-[d])\wfo,& if $\bi\in I^\gamma_+$,\\
                    (t^d\wyo_r+[d])\wfo,& if $\bi\in I^\gamma_-$.
           \end{dcases*}
      \]
      If $\F$ is an $\O$-module let $\RO[\F]=\RO\otimes_\O\F$.
    \end{Definition}

    To show that $\RO$ is finitely generated as an $\O$-module we need
    the following technical lemma, which is an analogue of
    \cite[Lemma~4.31]{HuMathas:SeminormalQuiver}.

    \begin{Lemma}\label{L:yorder}
      Suppose that $1\leq r\leq n$ and $\ei\in I^\gamma$. If $\bi\notin
      I^\gamma_+\cup I^\gamma_-$ then $\wfo=0$ and if $\bi\in
      I^\gamma_\pm$ then there exists a set $X_r(\ei)\subseteq
      e\Z\times\N$ such that
      \[
      \prod_{(c,m)\in X_r(\ei)}(\wyo_r\mp[c])^m\fo=0
      \]
      in $\RO$.
    \end{Lemma}

    \begin{proof}
      Arguing exactly as in proof of \autoref{L:ZeroSequence}, if
      $\bi\in I^\gamma$ then $\wfo\ne0$ only if $i_1=0$
      and $i_2=\pm1$.  That is, $\wfo\ne0$ only if $\bi\in
      I^\gamma_+\cup I^\gamma_-$. Hence, we may assume that
      $\bi\in I^\gamma_+\cup I^\gamma_-$.

      Checking the relations in \autoref{D:RO}, $\RO$ has an
      automorphism~$\dhash$ such that
      \[    (\wpsio_r)^{\dhash}=-\wpsio_r, \quad
            (\wyo_s)^{\dhash}=-\wyo_s\quad\text{and}\quad
            (\wfo)^{\dhash}=\wfo[-\bi],
      \]
      for all $1\le r<n$, $1\le s\le n$ and $\bi\in I^\gamma$.
      Therefore, it is enough to consider the case when $\bi\in I^\gamma_+$.

      By \autoref{D:RO}, $\wyo_1\wfo=0$, so we may take
      $X_1(\bi)=\set{(0,1)}$.  As
      $\wpsio_1\wfo=\wfo[s_1\cdot\bi]\wpsio_1$, it follows that
      $\wpsio_1=0$.  Therefore, if $\wfo\ne0$ then
      $0=(\wpsio_1)^2\wfo=(\wyo_1-\wyo_2)\wfo=-\wyo_2\wfo$, so
      $\wyo_2=0$.  Hence, we can set $X_2(\bi)=\set{(0,1)}$, for all
      $\bi\in I^\gamma$.

      Now consider $\wyo_3\wfo$, for $\bi\in I^\gamma_+$. If $i_2=i_3$ then the
      commutation relations for~$\wpsio_2$ and~$\wyo_3$ give
      $\wfo=(\wyo_3\wpsio_2-\wpsio_2\ddyo{-e}_2)\wfo
            =(\wyo_3+[-e])\wpsio_2\wfo$ since $\wyo_2=0$.
      Similarly, $\wfo=\wpsio_2(\wyo_3+[-e])\wfo$. Therefore,
      $\wfo=(\wyo_3+[e])(\wpsio_2)^2(\wyo_3+[-e])\wfo=0$.
      Hence, we can assume that $i_2\ne i_3$. If $i_2\noedge i_3$ and
      $\bi\in I^\gamma_+$ then
      \begin{align*}
      (\wyo_3-[-e])\fo&=(\wyo_3-[-e])(\wpsio_2)^2\wfo
                       =(\wyo_3-[-e])\wpsio_2\wfo[s_2\cdot\bi]\wpsio_2\\
                      &=\wpsio_2(\wyo_2+[-e]-[-e])\wfo[s_2\cdot\bi]\wpsio_2
                       =0.
      \end{align*}
      Hence, if $i_2\noedge i_3$ set $X_3(\bi)=\set{(-e,1)}$.
      Similarly, if $i_2\to i_3$ then
      \begin{align*}
        (\wyo_3-[-e])\wyo_3\wfo &=(\wyo_3-[-e])(\wyo_2-\wpsio_2)^2\wfo
                  =-(\wyo_3-[-e])(\wpsio_2)^2\wfo=0.
      \end{align*}
      Consequently, we can set $X_3(\bi)=\set{(-e,1),(0,1)}$.
      The case when $i_2\leftarrow i_3$ is similar and easier with
      $X_3(\bi)=\set{(0,1)}$.

      The last two paragraphs show that if $1\le r\le 3$ and $\bi\in I^\gamma_+$
      then there exists a set $X_r(\bi)\subseteq e\Z\times\N$ such that
      $\prod_{(c,m)\in X_r(\bi)}(\wyo_r-[c])^m\wfo=0$. If $3\le r<n$ and $\bi\in
      I^\gamma_+$ then $s_r\cdot\bi\in I^\gamma_+$. Moreover, the
      elements $\wpsio_3,\dots,\wpsio_{n-1}$ and $\wyo_4,\dots,\wyo_n$
      satisfy the same defining relations as $\psip_3,\dots,\psip_{n-1}$,
      $\yp_4,\dots,\yp_n$. Therefore, the inductive argument in
      \cite[Lemma~4.31]{HuMathas:SeminormalQuiver} shows that
      there exists a set $X_r(\bi)\subset e\Z\times\N$ such that
      \[\prod_{(c,m)\in X_r(\bi)}(\wyo_r-[c])^m\wfo=0,
          \qquad\text{ for }\bi\in I^\gamma_+\text{ and }1\le r\le n.
      \]
      (Note that if $1\le r\le 3$, or if $\bi\in I^\gamma_-$, then the
      argument from \cite{HuMathas:SeminormalQuiver} does not apply
      because $\wpsio_1$ and~$\wpsio_2$ satisfy slightly different
      relations to the corresponding elements considered in that paper.)
    \end{proof}

    Finally, we are able to prove the enhanced version of \autoref{T:BKiso}
    that we use to prove our main result. If $A$ is an $\O$-algebra let
    $\FHSn[A]=\HO\otimes_\O A$.

    \begin{Theorem}\label{T:MainTheorem}
      Suppose that $\gamma\in\Qpm$ and that $(\O,t)$ is the idempotent
      subring defined in $\autoref{E:AlternatingCS}$. Then
      $\RO\cong\HO[\gamma]$ as $\O$-algebras.
    \end{Theorem}

  \begin{proof}
    By the results in \autoref{S:AlternaatingCoefficients}, from
    \autoref{P:NewGenerators} onwards,
    there is a unique surjective algebra homomorphism
    $\RO\twoheadrightarrow\HO[\gamma]$ such that
    \[
    \wpsio_r\mapsto\psio_r,\qquad
    \wyo_s\mapsto\yo_s\quad\text{and}\quad
    \wfo\mapsto\fo,
    \]
    for $1\le r<n$, $1\le s\le n$ and $\bi\in I^\gamma$. To prove that this
    map is an isomorphism we use the argument from
    \cite[Theorem~4.32]{HuMathas:SeminormalQuiver} to show that $\RO$ is
    free as an $\O$-module with the same rank as~$\HO[\gamma]$.

    First, using the relations in \autoref{D:RO} it is straightforward to
    show that $\RO$ is spanned by elements of the form $f_w(\dot y)\wpsio_w\wfo$,
    where $f_w(\dot y)$ is a polynomial in $\O[\wyo_1,\dots,\wyo_n]$,
    $\bi\in I^\gamma$ and for each $w\in\Sn$ we fix a reduced expression
    $w=s_{r_1}\dots s_{r_k}$ and set
    $\wpsio_w=\wpsio_{r_1}\dots\psio_{r_k}$. Hence, $\RO$ is finitely
    generated as an $\O$-module by \autoref{L:yorder}.

    Next, let $\m=x\O$ be the maximal ideal of~$\O$ and set $\F=\O/\m$ and
    $\xi=t+\m\in\F$. Then~$\xi$ has quantum characteristic~$e$ because if
    $k\in\Z$ then $[k]_t\in\J(\O)=\m$ if and only if $k\in e\Z$ by
    \autoref{D:idempotentSub}. By \autoref{D:RO}, the relations in
    $\RO[\F]$ collapse and become the KLR relations for $\FRSn[\F]_\gamma$
    given in \autoref{D:klr}. That is, $\RO[\F]\cong\FRSn[\F]_\gamma$ as
    $\F$-algebras. Consequently,
    \[\dim_\F\RO[\F]=\dim_\F\FRSn[\F]_\gamma=\rank_\O\HO[\gamma],\]
    where the last equality follows by \cite[Theorem
    4.20]{BK:GradedDecomp} (alternatively, use \autoref{T:BKiso}).
    Since $\m$ is the unique maximal ideal~$\m$ of $\O$, and $\RO$ is
    finitely generated as an $\O$-module, Nakayama's Lemma implies
    that~$\RO$ is free as an $\O$-module of rank
    $\dim_\F\FHSn[\F]_\gamma=\rank_\O\HO[\gamma]$. Hence, as an
    $\O$-module, $\RO$ is free of the same rank as~$\HO[\gamma]$. Since
    $\HO[\gamma]$ is also free over~$\O$, it follows that the surjective
    algebra homomorphism $\RO\twoheadrightarrow\HO[\gamma]$ given in the
    first paragraph of the proof is actually an isomorphism and the
    theorem is proved.
  \end{proof}

  Recalling \autoref{E:HAnBlocks}, for $\gamma\in\Qpm$ define
  $\FHAn[\F]_\gamma=\HO(\An)_\gamma\otimes_\O\F$. Then
  $\FHAn[\F]_\gamma$ is a direct summand of $\FHAn[\F]$ by
  \autoref{C:HAnDecomp}. By construction, $F$ is (isomorphic to) a subfield of~$\F$
  and the algebra $\FHAn$ is the $F$-subalgebra of $\FHAn[\F]$
  generated by the elements $T_1,\dots,T_{n-1}$. By \autoref{L:foIdempotents} and
  \autoref{P:HOdecomp}, $e_\gamma=\fo[\gamma]\otimes1_\F$ is central
  idempotent in~$\FHAn$. Define
  \[ \FHAn_\gamma=\FHAn e_\gamma.  \]
  Then $\FHAn_\gamma$ is the $F$-subalgebra of~$\FHAn[\F]_\gamma$
  generated by $T_1e_\gamma,\dots,T_{n-1}e_\gamma$.

  We are assuming that $F$ is a field and that $\xi\in F$ an element of quantum
  characteristic~$e$.  Recall from before \autoref{T:Main} that a field~$F$ is
  \textbf{large enough} for~$\xi$ if~$F$ contains squareroots $\sqrt{\xi}$
  and $\sqrt{1+\xi+\xi^2}$ whenever $e>3$. (In particular, if $e=3$ then
  any field is large enough for~$\xi$.)

\begin{Theorem}\label{altiso}
  Suppose that $\gamma\in\Qpm$, $e>2$ and that $\xi\in F$ an element of
  quantum characteristic~$e$. Let $F$ be a large enough field for~$\xi$ of
  characteristic different from~$2$. Then $\FHAn_\gamma\cong\FRAn_\gamma$.
\end{Theorem}

\begin{proof}
  Let $(\O,t)$ be the idempotent subring given in
  \autoref{L:IdempotentSubring}, starting from~$F$ and~$\xi$, and let
  $\F=\O/\m$, where $\m=x\O$ is the maximal ideal of~$\O$. Now
  $\RO[\F]\cong\FHSn[\F]_\gamma$ by \autoref{T:MainTheorem}, so there is an
  isomorphism of $\F$-algebras $\Theta\map{\FRSn[\F]_\gamma}{\FHSn[\F]_\gamma}$ given by
  \[
      \psi_r\otimes1_\F\mapsto\psio_r\otimes1_\F,\qquad
      y_s\otimes1_\F\mapsto\yo_s\otimes1_\F\quad\text{and}\quad
      e(\bi)\otimes1_\F\mapsto\fo\otimes1_\F,
  \]
  for $1\le r<n$, $1\le s\le n$ and $\bi\in I^\gamma$. By
  \autoref{C:PsioHash} the following diagram commutes:
  \begin{center}
    \begin{tikzpicture}[>=stealth,->,shorten >=2pt,looseness=.5,auto]
      \matrix (M)[matrix of math nodes,row sep=1cm,column sep=16mm]{
           \FRSn[\F]_\gamma & \FHSn[\F]_\gamma \\
           \FRSn[\F]_\gamma & \FHSn[\F]_\gamma \\
       };
       \draw(M-1-1)--node[above]{$\Theta$}(M-1-2);
       \draw(M-2-1)--node[above]{$\Theta$}(M-2-2);
       \draw(M-1-1)--node[left]{$\sgn$}(M-2-1);
       \draw(M-1-2)--node[right]{$\#$}(M-2-2);
    \end{tikzpicture}
  \end{center}
  Therefore, $\Theta$ restricts to an isomorphism
  $\Theta\map{\FRAn[\F]_\gamma}{\FHAn[\F]_\gamma}$.

  We have now shown that $\FRAn[\F]_\gamma$ and $\FHAn[\F]_\gamma$ are
  isomorphic over~$\F$ but, of course, we want the isomorphism over~$F$,
  which is a subfield of~$\F$. Since~$F$ is large enough for~$\xi$, by
  \autoref{D:psio} the generators of~$\FHSn[\F]_\gamma$ listed in
  \autoref{P:NewGenerators} all belong
  to~$\FHSn[F]_\gamma$, which we consider as a subalgebra
  of~$\FRSn[\F]_\gamma$. The coefficients in the relations of
  \autoref{D:RO} also belong to $\FHSn[F]_\gamma$. Hence, there is a
  surjective algebra homomorphism
  $\FRSn[F]_\gamma\twoheadrightarrow\FHSn[F]_\gamma$. Counting
  dimensions, this map is an isomorphism
  so~$\FRSn[F]_\gamma\cong\FHSn[F]_\gamma$ as $F$-algebras.  Applying
  \autoref{C:PsioHash}, as above, it follows that
  $\FRAn[F]_\gamma\cong\FHAn[F]_\gamma$ as required.
\end{proof}

In view of \autoref{C:HAnDecomp}, we obtain \autoref{T:Main} from the
introduction.

\begin{Corollary}\label{C:AlternatingHecke}
  Let $F$ be a field of characteristic different from~$2$ and let
  $\xi\in F$ be an element of quantum characteristic~$e$.  Suppose that
  $e>2$ and that $F$ is a large enough field for~$\xi$. Then
  $\FHAn\cong\FRAn$.
\end{Corollary}

Hence, as noted in \autoref{C:HAnZGraded}, the alternating Hecke algebra
$\FHAn[F]$ is a $\Z$-graded algebra. In particular, $F\An$ is a
$\Z$-graded algebra when $F$ is large enough for~$\xi=1$.

\section{A homogeneous basis for $\FHAn[F]$}\label{S:Basis}
  We have now proved \autoref{T:Main} and \autoref{T:MainRelations} from
  the introduction. It remains to prove \autoref{T:GDim}, which
  gives the graded dimension of~$\RAn$. To do this we give a homogeneous
  basis for~$\RAn$ by combining the two graded cellular bases of~$\RSn$
  defined by Hu and the second-named author
  \cite{HuMathas:GradedCellular}. In order to define these bases we need some definitions.

Fix a partition $\lambda\in\Parts$. If $A=(r,c)$ and $B=(s,d)$ are nodes
of $\lambda$ then $A$ is \textbf{strictly above} $B$, or $B$ is
\textbf{strictly below} $A$, if $r<s$. Following
Brundan, Kleshchev and Wang~\cite[\S1]{BKW:GradedSpecht}, define
integers
\begin{align*}
    d_A(\lambda)&=\#\SetBox{addable $i$-nodes of $\lambda$\\[-2mm]
                            strictly below $A$}
                 -\#\SetBox[36]{removable $i$-nodes of $\lambda$\\[-2mm]
                            strictly below $A$},\\
    d^A(\lambda)&=\#\SetBox{addable $i$-nodes of $\lambda$\\[-2mm]
                            strictly above $A$}
                 -\#\SetBox[36]{removable $i$-nodes of $\lambda$\\[-2mm]
                            strictly above $A$}.
\end{align*}

\begin{Definition}[\protect{%
  Brundan, Kleshchev and Wang~\cite[\S1]{BKW:GradedSpecht}}]
  Let $\t$ be a standard $\lambda$-tableau, for $\lambda\in\Parts$, and let
  $A=\t^{-1}(n)$. Then the \textbf{degree} $\deg\t$ and $\codeg\t$
  of~~$\t$ are defined inductively by
  \begin{align*}
    \deg\t&=\begin{cases}
      \deg\t_{\downarrow{(n-1)}}+d_A(\lambda),&\text{if },n>0,\\
              1,&\text{if }n=0,\\
    \end{cases}\\
    \codeg\t&=\begin{cases}
      \codeg\t_{\downarrow{(n-1)}}+d^A(\lambda),&\text{if }n>0,\\
        1,&\text{if }n=0,
    \end{cases}
  \end{align*}
\end{Definition}

Fix $\lambda\in\Parts$. If $1\le m\le n$ and $\t\in\Std(\lambda)$
let $\col_m(\t)=c$ if $m$ appears in column~$c$ of~$\t$ and let
$\row_m(\t)=r$ if $m$ appears in row~$r$ of~$\t$.
Following~\cite{HuMathas:GradedCellular} set
\[
  y_\lambda=\prod_{\substack{1\le m\le n\\\col_m(\t^\lambda)\equiv0\pmod e}}y_m
  \qquad\text{and}\qquad
  y'_\lambda=\prod_{\substack{1\le m\le n\\\row_m(\t_\lambda)\equiv0\pmod e}}y_m.
\]
As we are considering the special case when $\Lambda=\Lambda_0$ the
definitions of~$y_\lambda$ and~$y'_\lambda$ from
\cite{HuMathas:GradedCellular} simplify and are equivalent to the formulas above.

If $\t\in\Std(\lambda)$ define
permutations $d(\t)$ and $d'(\t)$ in~$\Sn$ by
\[
  \t=d(\t)\t^\lambda\qquad\text{and}\qquad\t=d'(\t)\t_\lambda,
\]
where $\Sn$ acts on $\t$ by permuting its entries. For each $w\in\Sn$
fix a reduced expression $w=s_{r_1}\dots s_{r_k}$ and define
$\psi_w=\psi_{r_1}\dots\psi_{r_k}\in\RSn$. In general, $\psi_w$ depends
upon the choice of reduced expression for~$w$.

\begin{Definition}[\protect{%
   Hu and Mathas~\cite[Definitions~5.1 and~6.9]{HuMathas:GradedCellular}}]
Suppose that $\s,\t\in\Std(\lambda)$, for $\lambda\in\Parts$. Set
$\bi^\lambda=\res(\t^\lambda)$ and $\bi_\lambda=\res(\t_\lambda)$
and define
\[
\psi_{\s\t}=\psi_{d(\s)}y_\lambda e(\bi^\lambda) \psi_{d(\t)}^*
\quad\text{and}\quad
\psi_{\s\t}'= \psi_{d'(\s)}y_\lambda'e(\bi_\lambda)\psi_{d'(\t)}^*.
\]
\end{Definition}

By construction, $\psi_{\s\t}$ and $\psi'_{\s\t}$ are homogeneous
elements of~$\RSn$.

\begin{Remark}
For the reasons explained in \cite[Remark~3.12]{HuMathas:QuiverSchurI},
we are following the conventions of
\cite{HuMathas:QuiverSchurI,KMR:UniversalSpecht} here rather than those
of~\cite{HuMathas:GradedCellular}. In particular, the
element~$\psi'_{\s\t}$ defined above is equal to $\psi'_{\s'\t'}$ in the
notation of \cite{HuMathas:GradedCellular}.
\end{Remark}

\begin{Theorem}[Hu and Mathas~\cite{HuMathas:GradedCellular}
                and Li~\cite{GeLi:IntegralKLR}]\label{T:SnCellularBasisThm}
Let $\Zcal$ be a commutative ring and suppose that $e>2$ and $n\geq 0$.
Then $\RSn$ is free as a $\Zcal$-module with homogeneous bases
$\set{\psi_{\s\t}| \s,\t\in\Std^2(\Parts)}$ and
$\set{\psi'_{\s\t}| \s,\t\in\Std^2(\Parts)}$.
Moreover, if~$(\s,\t)\in\Std^2(\Parts)$ then
$\deg\psi_{\s\t}=\deg\s+\deg\t$ and $\deg\psi'_{\s\t}=\codeg\s+\codeg\t$.
\end{Theorem}

This result was first proved in~\cite{HuMathas:GradedCellular} with some
restrictions on the ring~$\Zcal$. Li's~\cite{GeLi:IntegralKLR} extension
of this result to arbitrary rings is a difficult theorem. A second proof
of this result, using geometry, is given in \cite{StroppelWebster:QuiverSchur}.

Although we will not need this, by \cite[Theorems~5.8
and~6.11]{HuMathas:GradedCellular} both of the bases of
\autoref{T:SnCellularBasisThm} are graded \textit{cellular} bases of
$\RSn$ in the sense of Graham and
Lehrer~\cite{GL,HuMathas:GradedCellular}.

\begin{Example}\label{psibasiscalc}
Set $n=e=3$ and set
\[
        \s=\Tableau[-1]{{1,2,3}},\quad
        \t=\Tableau{{1,2},{3}},\quad
        \u=\Tableau{{1,3},{2}}\quad\text{and}\quad
        \v=\Tableau[-4]{{1},{2},{3}}.
\]
Using the definitions,
\begin{xalignat*}{2}
  \psi_{\s\s}&=y_3e(012),         & \psi_{\v\v}'&=y_3e(021),\\
  \psi_{\t\t}&=e(012),            & \psi_{\u\u}'&=e(021),\\
  \psi_{\t\u}&=e(012)\psi_2,      & \psi_{\u\t}'&=e(021)\psi_2,\\
  \psi_{\u\t}&=\psi_2e(012),      & \psi_{\t\u}'&=\psi_2e(021),\\
  \psi_{\u\u}&=\psi_2e(012)\psi_2,& \psi_{\t\t}'&=\psi_2e(021)\psi_2,\\
  \psi_{\v\v}&=e(021),            & \psi_{\s\s}'&=e(012).
\end{xalignat*}
Notice that $\psi_{\t\u}=\psi_2e(021)=\psi'_{\t\u}$,
$\psi_{\u\u}=\psi_2^2e(021)=-y_3e(021)=-\psi'_{\v\v}$ and, similarly,
$\psi'_{\t\t}=-y_3e(012)=-\psi_{\s\s}$. Therefore, up to sign, the
$\psi$ and $\psi'$ bases coincide with the basis of $\RSn[3]$ given
in \autoref{Ex:OS3}.
\end{Example}

We now use the two homogeneous bases for $\RSn$ from
\autoref{T:SnCellularBasisThm} to construct homogeneous bases for
$\RAn$.  The next result, which follows easily from the definitions,
shows that the $\psi$ and $\psi'$-bases are interchanged by the sign
automorphism.

\begin{Lemma}[\protect{%
  Hu and Mathas~\cite[Proposition 3.26]{HuMathas:QuiverSchurI}}]
  \label{L:psisgn}
  Suppose that $\s,\t\in\Std(\Parts)$. Then
  $\psi_{\s\t}^\sgn = (-1)^{\ell(d(\s))+\ell(d(\t))+\deg \t^\lambda}
         \psi_{\s'\t'}'$.
  In particular, $\deg\psi_{\s\t}=\deg\psi_{\s'\t'}'$.
\end{Lemma}

In fact, it follows easily from the definitions that $\deg\t=\codeg\t'$
for any standard tableau~$\t$. Hence,
$\deg\psi_{\s'\t'}'=\codeg\s'+\codeg\t'=\deg\s+\deg\t=\deg\psi_{\s\t}$
as claimed.

Recall from \autoref{P:RpSnIsomorphism} that
$\RSn_\gamma\cong\RpSn[\gamma]=\RpSnp\oplus\RpSnp[-]$ is a $(\Z_2\times\Z)$-graded
algebra and that $\RAn_\gamma\cong\RpSnp$ by \autoref{C:EvenBit}. To find a
basis for~$\RAn_\gamma$ we first give a basis of~$\RSn_\gamma$ that is homogeneous
with respect to the $(\Z_2\times\Z)$-grading.

\begin{Definition}
Fix $e>2$ and $\lambda\in\Parts$. For $\s,\t\in \Std(\lambda)$ define
elements
\[
    \Psi^+_{\s\t}=\psi_{\s\t}+\psi_{\s\t}^\sgn \quad\text{and}\quad
    \Psi^-_{\s\t}=\psi_{\s\t}-\psi_{\s\t}^\sgn.
\]
\end{Definition}

Fix $(\s,\t)\in\Std^2(\Parts)$. By \autoref{L:psisgn},
$\Psi^+_{\s\t}$ and $\Psi^-_{\s\t}$ are homogeneous with respect to the $\Z$-grading
on~$\RSn$. Furthermore, by \autoref{C:EvenBit},~$\Psi^+_{\s\t}$ is even
and~$\Psi^-_{\s\t}$ is odd with respect to the $\Z_2$-grading. Hence,
the elements $\Psi^\pm_{\s\t}$ are homogeneous with respect to the
$(\Z_2\times\Z)$-grading on~$\RSn$.

Fix $(\s,\t)\in\Std^2(\Parts)$ and set $\bi^\s=\res(\s)$ and
$\bi^\t=\res(\t)$.  By \cite[(3.13)]{HuMathas:QuiverSchurI}, if
$\bi,\bj\in I^n$
\begin{equation}\label{E:psiIdempotents}
  e(\bi)\psi_{\s\t}e(\bj)=\delta_{\bi^\s\bi}\delta_{\bi^\t\bj}\psi_{\s\t}
  \quad\text{and}\quad
  e(\bi)\psi'_{\s\t}e(\bj)=\delta_{\bi^\s\bi}\delta_{\bi^\t\bj}\psi'_{\s\t}
\end{equation}
Set $\eps_1(\bi)=e(\bi)-e(-\bi)$ and
$\eps_0(\bi)=\eps_1(\bi)^2=e(\bi)+e(-\bi)$.
Observe that \autoref{E:psiIdempotents} implies that if
$\bi\in I^\gamma_+$ and $(\s,\t)\in\Std(\Parts)$ with
$\bi^\s=\res(\s)\in I^\gamma$ and $\bi^\t=\res(\t)\in I^\gamma$ then
\begin{equation}\label{E:epsPsi}
\eps_1(\bi)\Psi^+_{\s\t}
        \begin{cases}
          \phantom{-}\Psi^-_{\s\t}&\text{if }\bi=\bi^\s,\\
          -\Psi^-_{\s\t}&\text{if }\bi=-\bi^\s,\\
          \phantom{-}0,&\text{otherwise},
        \end{cases}
\quad\text{and}\quad
\Psi^+_{\s\t}\eps_1(\bi)=
        \begin{cases}
          \phantom{-}\Psi^-_{\s\t}&\text{if }\bi=\bi^\t,\\
          -\Psi^-_{\s\t}&\text{if }\bi=-\bi^\t,\\
          \phantom{-}0,&\text{otherwise}.
        \end{cases}
\end{equation}
Since $\eps_1(-\bi)=-\eps_1(\bi)$ and $\eps_0(\bi)=\eps_1(\bi)^2$,
this readily implies the corresponding formulas for the action of
$\eps_a(\bi)$ on $\Psi^\pm_{\s\t}$, for all $a\in\Z_2$ and $\bi\in
I^\gamma$.

The \textbf{dominance} order on~$\Parts$ is the partial order~$\gedom$
given by $\lambda\gedom\mu$ if
\[
    \sum_{j=1}^k\lambda_j\ge\sum_{\j=1}^k\mu_j\qquad\text{ for all }k\ge0
\]
Write $\lambda\gdom\mu$ if $\lambda\gedom\mu$ and $\lambda\ne\mu$.

For $\lambda\in\Parts$ define
$\Std_+(\lambda)=\set{\s\in\Std(\lambda)|\res(\s)\in I^n_+}$.  We will
use this set to index a homogeneous basis for~$\RSn$, with respect to its
$(\Z_2\times\Z)$-grading.
The following simple combinatorial result is probably well-known.

\begin{Lemma}\label{L:CountingLemma}
  Suppose that $n\ge0$. Then
  \[
     \sum_{\substack{\lambda\in\Parts\\\lambda\gedom\lambda'}}
         |\Std_+(\lambda)|\cdot|\Std(\lambda)|
         =\frac{n!}{2}.
  \]
\end{Lemma}

\begin{proof} Implicit in \autoref{T:SnCellularBasisThm}, is the
  well-known
  fact that $n!=\sum_{\lambda\in\Parts}|\Std(\lambda)|^2$.
  Since $|\Std(\lambda)|=|\Std(\lambda')|$, via the map $\t\mapsto\t'$,
  it follows that
  \begin{align*}
    \frac{n!}2
      &=\sum_{\substack{\lambda\in\Parts\\\lambda\gdom\lambda'}}|\Std(\lambda)|^2
         +\frac12\sum_{\substack{\lambda\in\Parts\\\lambda=\lambda'}}
                |\Std(\lambda)|^2\\
      &=\sum_{\substack{\lambda\in\Parts\\\lambda\gdom\lambda'}}
          \Bigl(|\Std_+(\lambda)|+|\Std_+(\lambda')|\Bigr)\cdot|\Std(\lambda)|
         +\sum_{\substack{\lambda\in\Parts\\\lambda=\lambda'}}
               |\Std_+(\lambda)|\cdot |\Std(\lambda)|.
  \end{align*}
\end{proof}

Recall from \autoref{S:GeneratorsRelations} that
$\Deg\map{\RSn}\Z_2\times\Z$ is the degree function for the
$(\Z_2\times\Z)$-grading on~$\RSn$.

\begin{Theorem}\label{T:OddEvenBasis}
  Let $\Zcal$ be a commutative ring such that $2$ is invertible
  in~$\Zcal$ and suppose that $e>2$ and $n\ge0$. Then $\RSn$ is free as
  a $\Zcal$-module with basis
  \[
  \set{\Psi^+_{\s\t}, \Psi^-_{\s\t}|\s\in\Std_+(\lambda), \t\in\Std(\lambda)
           \text{ for } \lambda\in\Parts}.
  \]
  Moreover, this basis is homogeneous with respect to the
  $(\Z_2\times\Z)$-grading on~$\RSn$.
\end{Theorem}

\begin{proof}
  We have already noted $\Psi^\pm_{s\t}$ is homogeneous with respect to
  the $(\Z_2\times\Z)$-grading. More precisely, if
  $(\s,\t)\in\Std^2(\Parts)$ then
  \[\Deg\Psi^+_{\s\t}=(0,\deg\s+\deg\t)\quad\text{and}\quad
    \Deg\Psi^-_{\s\t}=(1,\deg\s+\deg\t).
  \]

  Let $R_n$ be the $\Zcal$-submodule of $\RSn$ spanned by the
  elements in the statement of the theorem.
  Fix $\s\in\Std_+(\lambda)$ and $\t\in\Std(\lambda)$, for some
  $\lambda\in\Parts$. Since $2$ is invertible in~$\Zcal$,
  \[
      \psi_{\s\t}=\frac12\big(\Psi^+_{\s\t}+\Psi^-_{\s\t}\big)
      \quad\text{and}\quad
      \psi^\sgn_{\s\t}=\frac12\big(\Psi^+_{\s\t}-\Psi^-_{\s\t}\big).
  \]
  Hence, $\psi_{\s\t}, \psi'_{\s'\t'}\in R_n$ since
  $\psi'_{\s'\t'}=\pm\psi_{\s\t}^\sgn$ by \autoref{L:psisgn}.  Recall
  from \autoref{L:ZeroSequence} that $e(\bi)\ne0$ only if
  $\bi\in I^n_+$ or $\bi\in I^n_-$. Set
  $e_+=\sum_{\bi\in I^n_+}e(\bi)$ and $e_-=\sum_{\bi\in I^n_-}e(\bi)$.
  By \autoref{T:SnCellularBasisThm} and \autoref{E:psiIdempotents},
  as $\Zcal$-modules,
  \begin{align*}
    e_+\RSn&=\<\psi_{\s\t}\mid (\s,\t)\in\Std^2(\Parts)\text{ and }
            \res(\s)\in I^n_+ \>_\Zcal\subseteq e_+ R_n\\
  \intertext{
  Similarly, since $\res(\s')=-\res(\s)$, \autoref{T:SnCellularBasisThm}
  also implies that}
      e_-\RSn&=\<\psi'_{\s\t}\mid (\s,\t)\in\Std^2(\Parts)\text{ and }
            \res(\s)\in I^n_- \>_\Zcal\subseteq e_- R_n.
  \end{align*}
  Hence, $e_+R_n\subseteq e_+\RSn\subseteq e_+R_n$ and
  $e_-R_n\subseteq e_-\RSn\subseteq e_-R_n$, so that
  \[  \RSn=e_+\RSn\oplus e_-\RSn=e_+R_n\oplus e_-R_n=R_n. \]
  We have now shown that the set of elements $\set{\Psi^\pm_{\s\t}}$ in
  the statement of the theorem span~$\RSn$. Let $F$ be the field of
  fractions of~$\Zcal$. Using \autoref{L:CountingLemma} to count
  dimensions, it follows that $\set{\Psi^\pm_{\s\t}\otimes1_F}$ is a
  basis of~$\FRSn\cong\RSn\otimes_\Zcal F$. Hence,
  $\set{\Psi^\pm_{\s\t}}$ is $\Zcal$-linearly independent, completing the proof.
\end{proof}

By \autoref{C:EvenBit}, $\RAn$ is the even component of~$\RSn$, with
respect to the $\Z_2$-grading. Hence, we have the following.

\begin{Corollary}\label{C:EvenBasis}
  Let $\Zcal$ be a commutative ring such that $2$ is invertible
  in~$\Zcal$ and suppose that $e>2$ and $n\ge0$.
  Then $\RAn$ is free as a $\Zcal$-module with basis
  \[
  \set{\Psi^+_{\s\t}|\s\in\Std_+(\lambda) \text{ and } \t\in\Std(\lambda)
           \text{ for } \lambda\in\Parts}.
  \]
  Moreover, this basis is homogeneous with respect to the
  $\Z$-grading on~$\RAn$.
\end{Corollary}

\begin{Example}\label{psibasisalt}
Continuing the notation of \autoref{psibasiscalc},
\[
\Psi^+_{\s\s}= \mathcal{Y}_3= -\Psi^+_{\u\u}\quad
\Psi^+_{\t\t}= 1= \Psi^+_{\v\v}\quad\text{and}\quad
\Psi^+_{\t\u}= \Psi^+_2= -\Psi^+_{\t\u}.
\]
Hence,
$\set{\Psi^+_{\a\b}| \a\in\Std_+(\lambda)\text{ and }\b\in\Std(\lambda)}
      =\set{\Psi^+_{\s\s}, \Psi^+_{\t\t}, \Psi^+_{\t\u}}$
is the basis of~$\RAn$ constructed in \autoref{Ex:A3Basis}.
\end{Example}

If $M=\bigoplus_{d\in\Z}M_d$ is a $\Z$-graded module then its
\textbf{graded dimension} is the Laurent polynomial
\[
   \qdim M=\sum_{d\in\Z} (\dim M_d)\,q^d\in\Acal=\N[q,q^{-1}],
\]
where $q$ is an indeterminate over~$\Z$.  Adding up the degrees of the
homogeneous basis elements in \autoref{C:EvenBit} gives \autoref{T:GDim}
from the introduction.

\begin{Corollary}\label{gdim}
  Let $F$ be a field of characteristic different from~$2$ and suppose that
  $e>2$ and $n\ge0$. Then the graded dimension of~$\RAn$ is
\[
\qdim \RAn=\sum_{\lambda\in\Parts}
  \sum_{\substack{\s\in\Std_+(\lambda)\\\t\in\Std(\lambda)}} q^{\deg\s+\deg\t}.
\]
\end{Corollary}

%

By \autoref{E:epsPsi}, if $\gamma\in\Qpm$ then the basis of $\RAn$ given in
\autoref{C:EvenBasis} restricts to give a basis of~$\RAn_\gamma$. Note
that if $(\s,\t)\in\Std^2(\Parts)$ and $\res(\s)\in I^\gamma$ then
$\res(\t)\in I^\gamma$.

\begin{Corollary}\label{C:BlockBasis}
  Fix $\gamma\in\Qpm$ and let $\Zcal$ be a commutative ring such that
  $2$ is invertible in~$\Zcal$ and suppose that $e>2$ and $n\ge0$.
  Then $\RAn_\gamma$ is free as a $\Zcal$-module with basis
  $\set{\Psi^+_{\s\t}|(\s,\t)\in\Std^2(\Parts)\text{ for }
             \res(\s)\in I^\gamma_+}$.
\end{Corollary}

We next show that $\RAn_\gamma$ is a graded symmetric algebra. Repeating
the arguments leading to \autoref{C:BlockBasis} we obtain a second
homogeneous basis for~$\RAn_\gamma$. We will use the two homogeneous
bases of~$\RAn_\gamma$ to prove that $\RAn_\gamma$ is graded symmetric.

\begin{Corollary}\label{C:EvenBasisII}
  Fix $\gamma\in\Qpm$ and let $\Zcal$ be a commutative ring such that
  $2$ is invertible in~$\Zcal$ and suppose that $e>2$ and $n\ge0$.
  Then $\RAn_\gamma$ is free as a $\Zcal$-module with basis
  $\set{\Psi^+_{\s\t}|(\s,\t)\in\Std^2(\Parts)\text{ for }
             \res(\s)\in I^\gamma_-}$.
\end{Corollary}

Before we show that the blocks of $\RAn$ are graded symmetric algebras
we recall some definitions. A \textbf{trace form} on an algebra $A$ is a
linear map $\tau\map AF$ such that $\tau(ab)=\tau(ba)$, for all $a,b\in
A$. The algebra~$A$ is \textbf{symmetric} if $A$ is equipped with a
non-degenerate symmetric bilinear form $\theta: A\times A\rightarrow F$
which is associative in the following sense:
\[ \theta(xy,z)=\theta(x,yz),\quad\text{for all }x,y,z\in A. \]
A graded algebra $A$ is a \textbf{graded symmetric} algebra if there
exists a non-degenerate homogeneous trace form $\tau\map AF$. Suppose
that $A$ is equipped with a homogeneous anti-isomorphism $\sigma$ of
order~$2$. In view of \cite[Lemma~6.13]{HuMathas:GradedCellular}, giving
a non-degenerate trace form~$\tau$ is equivalent to requiring that the bilinear
form
\[\<a,b\> = \tau(ab^\sigma)\qquad\text{ for all }a,b\in A\]
is non-degenerate. We work interchangeably with the trace form $\tau$
and its associated bilinear form.

The algebras $\RSn$ and $\RAn$ are both symmetric algebras but their
homogeneous trace forms are defined on the blocks $\RAn_\gamma$ of these
algebras, for~$\gamma\in\Qpm$.  Recall that $\Qpm=Q^+_n/{\sim}$, where
$\alpha\sim\alpha'$ if $(\Lambda_i,\alpha)=(\Lambda_{-i},\alpha')$, for
all $i\in I$. Fix $\alpha\in Q^+_n$.  The \textbf{defect} of~$\alpha$ is the
non-negative integer
\[\defect\alpha=(\Lambda_0,\alpha)-\frac12(\alpha,\alpha),\]
where $(\ ,\ )\map{P^+\times Q^+}\Z$ is the pairing defined in
\autoref{S:KLRAlgebras}. Hence, if $\alpha\sim\alpha'$ then
$\defect\alpha=\defect\alpha'$. Therefore, if~$\gamma\in\Qpm$ we can
define the \textbf{defect} of~$\gamma$ to be
$\defect\gamma=\defect\alpha$, for any $\alpha\in\gamma$.  Set
\[\Parts[\alpha]=\set{\lambda\in\Parts|\bi^\lambda\in I^\alpha}
   \quad\text{and}\quad
  \Parts[\gamma]=\bigcup_{\alpha\in\gamma}\Parts[\alpha],
\]
for $\alpha\in Q^+_n$ and $\gamma\in\Qpm$. In the usual language from
the representation theory of the symmetric groups, the partitions
in~$\Parts[\gamma]$ have \textbf{$e$-weight} $\defect\gamma$.

The following useful fact is straightforward to establish from the
definitions.

\begin{Lemma}[\protect{%
  Brundan and Kleshchev~\cite[Lemma 3.11, 3.12]{BKW:GradedSpecht}}]
  \label{L:degcodeg}
  Let $\alpha\in Q^+_n$ and suppose that $\s\in\Std(\Parts[\alpha])$.
  Then $\deg\s+\codeg\s=\defect\alpha$.
\end{Lemma}

  As is well-known and easy to prove (see, for example,
  \cite[Proposition~1.16]{M:ULect}), the Iwhahori-Hecke $\HSn$ is a
  symmetric algebra with trace form $\tau$. Explicitly, if
  $h=\sum_{w\in\Sn}a_w T_w\in\HSn$ then $\tau(h)=a_1$. For $\alpha\in Q^+$
  let $\tau_\alpha$ be the homogeneous component of~$\tau$ of
  degree~$-2\defect\alpha$ restricted to $\HSn_\alpha$. Let $\<\ ,\ \>_\alpha$
  be the homogeneous bilinear form associated with~$\tau$.

  Let $\star$ be the unique homogeneous anti-isomorphism of~$\RSn$ that
  fixes each of the generators of~$\RSn$. Then
  $\psi_{\s\t}^\star=\psi_{\t\s}$ and
  $(\psi'_{\s\t})^\star=\psi'_{\t\s}$, for all
  $(\s,\t)\in\Std^2(\Parts)$.

  The following result was first proved in~\cite{HuMathas:GradedCellular}.

  \begin{Theorem}[\protect{%
      Hu and Mathas~\cite[Theorem~6.7]{HuMathas:GradedCellular}}]
    \label{T:RSnSymmetric}
    Suppose that $\alpha\in Q^+_n$ and that~$F$ is a field. Then the KLR
    algebra~$\Ralpha$ is a graded symmetric algebra with homogeneous
    bilinear form $\<\ ,\ \>_\alpha$ of degree $-2\defect\alpha$.
    Moreover, if $\s,\t\in\Std(\lambda)$ and $\u,\v\in\Std(\mu)$ then
    \[   \<\psi_{\s\t},\psi'_{\u\v}\>_\alpha
             =\begin{dcases*}
               c_{\s\t},& if $(\u,\v)=(\s,\t)$,\\
                   0,& if $(\u,\v)\notgedom(\s,\t)$,
               \end{dcases*}
     \]
     where $c_{\s\t}$ is a non-zero element of~$F$ that depends only
     on~$\s$ and~$\t$.
  \end{Theorem}

  Recall from \autoref{E:RAngamma} that
  $\RAn_\gamma=\big(\bigoplus_{\alpha\in\gamma}\Ralpha\big)^\sgn$.
  We will use the bilinear forms on $\RSn_\alpha$, for
  $\alpha\in\gamma$, to define a bilinear form on~$\RAn_\gamma$. We do
  not take the obvious extension of these forms to~$\RSn_\gamma$,
  however, because the arguments below require a $\sgn$-invariant form.
  If $h\in\RSn_\gamma$ then $h=\sum_{\bi\in I^\gamma}e(\bi)h_i$. Hence,
  define the trace form $\tau_\gamma\map{\RSn_\gamma}F$ by
  \begin{equation}\label{E:SgnInvariance}
       \tau_\gamma\big(e(\bi)h\big)=\begin{dcases*}
          \tau_\alpha\big(e(\bi)h\big),&
              if $\bi\in I^\alpha_+\subseteq I^\gamma_+$,\\
          \tau_{\alpha'}\big(e(-\bi)h^\sgn\big),&
              if $\bi\in I^{\alpha}_-\subseteq I^\gamma_-$.
        \end{dcases*}
  \end{equation}
  Importantly, $\tau_\gamma(h)=\tau_\gamma(h^\sgn)$, for all
  $h\in\RSn_\gamma$.  Let $\<\ ,\ \>_\gamma$ be the corresponding
  bilinear form on~$\RSn_\gamma$.

  Extend the dominance ordering $\gedom$ to standard tableaux by defining
  \[
     \s\gedom\t\quad\text{if}\quad
           \sh(\s_{\downarrow m})\gedom\sh(\t_{\downarrow m})
           \quad\text{ for }1\le m\le n,
  \]
  for $\s,\t\in\Std(\Parts)$. We can now prove that $\RAn_\gamma$ is a
  graded symmetric algebra.

  \begin{Theorem}\label{T:RAnGradeSymmetric}
    Let $F$ be a field of characteristic different from~$2$ and suppose
    that $e>2$ and $\gamma\in\Qpm$. Then $\RAn_\gamma$ is a graded
    symmetric algebra with homogeneous bilinear form $\<\ ,\ \>_\gamma$
    of degree $-2\defect\gamma$.
  \end{Theorem}

  \begin{proof}
    By \autoref{T:RAnGradeSymmetric} and the definitions above,
    $\<\ ,\ \>_\gamma$ is a (not necessarily associative) homogeneous bilinear form
    on~$\RSn_\gamma$ of degree $-2\deg\gamma$. By restriction, we can consider
    $\<\ ,\ \>_\gamma$ as a bilinear form on~$\RAn_\gamma$. By
    construction, $\<\ ,\ \>_\gamma$ is an associative bilinear form on~$\RAn$.

    We need to show that $\<\ ,\ \>_\gamma$ is non-degenerate on~$\RAn_\gamma$.
    To do this we use the two bases of $\RAn_\gamma$ given by
    \autoref{C:BlockBasis} and \autoref{C:EvenBasisII}. Fix
    $\lambda,\mu\in\Parts$ and tableaux $\s,\t\in\Std(\lambda)$,
    $\u,\v\in\Std(\mu)$ with $\res(\s)\in I^\gamma_+$
    and $\res(\u)\in I^\gamma_-$. By assumption, $\res(\s)\ne\res(\u)$ so
    $\<\psi_{\s\t},\psi_{\u\v}\>_\gamma =\tau_\gamma(\psi_{\s\t}\psi_{\v\u})
        =\tau_\gamma(\psi_{\v\u}\psi_{\s\t})=0$
    by~\autoref{E:psiIdempotents}. Hence,
    $\<\psi^\sgn_{\s\t},\psi^\sgn_{\u\v}\>_\gamma
               =\tau_\gamma(\psi_{\s\t}\psi_{\v\u}) =0$.
    Therefore,
    \begin{align*}
        \<\Psi^+_{\s\t},\Psi^+_{\u\v}\>_\gamma
          &=\<\psi_{\s\t}+\psi_{\s\t}^\sgn,\psi_{\u\v}+\psi_{\u\v}^\sgn\>_\gamma
           =\<\psi_{\s\t},\psi_{\u\v}^\sgn\>_\gamma
                +\<\psi_{\s\t}^\sgn,\psi_{\u\v}\>_\gamma\\
          &=2\tau_\gamma(\psi_{\s\t}\psi^\sgn_{\u\v})
           =\begin{dcases*}
            \pm 2c_{\s\t},& if $(\u',\v')=(\s,\t)$,\\
            0,& if $(\u',\v')\notgedom(\s,\t)$.
            \end{dcases*}
    \end{align*}
    The second last equality follows because
    $\tau_\gamma(h)=\tau_\gamma(h^\sgn)$ for $h\in\RSn_\gamma$, by
    \autoref{E:SgnInvariance}, and the last equality follows from
    \autoref{T:RAnGradeSymmetric}, since
    $\psi_{\u\v}^\sgn=\pm\psi'_{\u'\v'}$ by \autoref{L:psisgn}. Hence,
    by ordering the two bases $\set{\Psi^+_{\s\t}}$ and
    $\set{\Psi^+_{\u\v}}$ in a way that is compatible with dominance and
    reverse dominance, respectively, it follows that the Gram matrix
    $\big(\<\Psi^+_{\s\t},\Psi^+_{\u\v}\>_\gamma\big)$ is triangular
    with non-zero entries on the diagonal. Therefore, $\<\ ,\ \>_\gamma$
    is a non-degenerate associative bilinear form on~$\FHAn_\gamma$ of
    degree~$-2\defect\gamma$, so the theorem is proved.
  \end{proof}

  Finally, we describe the blocks and irreducible modules
  of~$\FHAn$.

  Let $\Res=\Res^{\RSn}_{\RAn}$ be the restriction functor from
  the category of finitely generated graded $\RSn$-modules to the
  the category of finitely generated graded $\RAn$-modules.

  Let $\Klesh=\set{\mu\in\Parts|\mu_r-\mu_{r+1}<e\text{ for all }r\ge1}$
  be the set of \textbf{$e$-restricted} partitions of~$n$. By
  \cite[Theorem~4.11]{BK:GradedKL}, there is a unique self-dual
  irreducible graded $\RSn$-module $D^\mu$ for each $e$-restricted
  partition~$\mu$ and
  \[\set{D^\mu\<d\>|\mu\in\Klesh\text{ and }d\in\Z}\]
  is a complete set of pairwise non-isomorphic irreducible graded
  $\RSn$-modules. By \cite[Corollary~5.11]{HuMathas:GradedCellular}, the
  module $D^\mu$ arises as a quotient of the corresponding graded Specht
  module~\cite{BKW:GradedSpecht,HuMathas:GradedCellular}.

  If $M$ is an $\RSn$-module let $M^\sgn$ be the $\RSn$-module that is
  isomorphic to~$M$ as a vector space but where the $\RSn$-action is
  twisted by~$\sgn$.  By \cite[Theorem~3.6.6]{Mathas:Singapore},
  $(D^\mu)^\sgn\cong D^{\mull(\mu)}$ where $\mull\map\Klesh\Klesh$ is
  the \textbf{Mullineux map}. Therefore, a straightforward application of
  Clifford theory implies that if $\mu\ne\mull(\mu)$ then $\Res D^\mu\cong
  \Res D^{\mull(\mu)}$ is an irreducible graded $\RAn$-module and if
  $\mu=\mull(\mu)$ then, over an algebraically closed field, $\Res
  D^\mu= D^\mu_+\oplus D^\mu_-$, for non-isomorphic irreducible graded
  $\RAn$-modules $D^\mu_+$ and $D^\mu_-$. Set
  \[
     \Klesh^{\mull\gdom}=\set{\mu\in\Klesh|\mull(\mu)\gdom\mu}
        \quad\text{and}\quad
     \Klesh^\mull=\set{\mu\in\Klesh|\mu=\mull(\mu)}.
  \]
  Clifford theory implies that every irreducible graded $\RAn$-module arises in
  the manner described above, so we obtain the following.

  \begin{Theorem}\label{T:GRadedSimples}
    Suppose that $F$ is an algebraically closed field of characteristic
    different from~$2$ and that $e>2$. Then
    \[\set{D^\mu\<d\>|d\in\Z \text{ and } \mu\in\Klesh^{\mull\gdom}}\cup
    \set{D^\mu_+\<d\>, D^\mu_-\<d\>|d\in\Z \text{ and } \mu\in\Klesh^\mull}.\]
   is a complete set of pairwise non-isomorphic irreducible graded $\RAn$-modules.
   \end{Theorem}

   In the semisimple case, $\Klesh=\Parts$ and $\mull(\mu)=\mu'$. If
   $\mu\in\Parts$ and $\mu=\mu'$ then
   \cite[Proposition~3.9]{MathasRatliff} gives an explicit construction
   of the modules $D^\mu_+$ and $D^\mu_-$ over the field of
   fractions~$\K$ of the idempotent subring $\O$ from
   \autoref{D:squareroots}.

   Finally we turn to the blocks of $\RAn$. By \autoref{C:RAnBlocks}, \[  \RAn =
   \bigoplus_{\gamma\in\Qpm}\RAn_\gamma,\]
   where $\RAn_\gamma$ is a two-sided graded ideal of~$\RAn$. Our last result
   says that if $\gamma\in\Qpm$ then $\RAn_\gamma$ is a block, or
   indecomposable two-sided ideal, of~$\RAn$ except when
   $\defect\gamma=0$ and $|\gamma|=1$. In the traditional language of
   the symmetric groups, $\defect\gamma=0$ if and only if  the
   partitions in $\Parts[\gamma]$ are $e$-cores and, in this case,
   $|\gamma|=1$ if and only if $\Parts[\gamma]=\set{\mu}$, where
   $\mu=\mull(\mu)$ is a Mullineux self-conjugate partition. As~$\mu$ is
   an $e$-core, $D^\mu$ is an irreducible Specht module and $\Res
   D^\mu=D^\mu_+\oplus D^\mu_-$, similar to the situation considered in
   the last paragraph.

  \begin{Theorem}\label{T:RAnBlocks}
    Suppose that $F$ is a field of characteristic different from~$2$
    and that $e>2$. Let $\gamma\in\Qpm$. Then:
    \begin{enumerate}
      \item If $|\gamma|=2$ or $\defect\gamma>0$ then
      $\FRAn_\gamma$ is an indecomposable two-sided graded ideal of~$\FRAn$.
      \item If $F$ is algebraically closed, $|\gamma|=1$ and $\defect\gamma=0$
      then $\FRAn_\gamma$ is a direct sum of two conjugate matrix
      algebras.
    \end{enumerate}
  \end{Theorem}

  \begin{proof}
    This follows by the general theory of covering blocks for
    $\Z_2$-graded algebras as can be found, for example, in
    \cite{Witherspoon:Clifford}. In more detail a block~$A$ of $\HSn$
    \textit{covers} a block~$B$ of $\RAn$ if~$B$ is a direct summand of
    the restriction of~$A$ to~$\RAn$. Since $|\Z_2|=2$ the blocks
    of~$\RAn$ are covered by at most two blocks of~$\HSn$ and, in
    particular, $\RAn_\gamma$ is indecomposable if $|\gamma|=2$.
    If $|\gamma|=1$ and $\defect\gamma>0$ then there exists a partition
    $\mu\in\Klesh\cap\Parts[\gamma]$ such that
    $\mu\ne\mull(\mu)\in\Parts[\gamma]$. Therefore, $\RAn_\gamma$ is a
    self-conjugate block, so that it is indecomposable. Finally, if~$F$
    is algebraically closed, $\defect\gamma=0$ and $|\gamma|=1$ then
    $\RAn_\gamma$ is the direct sum of two matrix algebras, in view of
    the remarks in the paragraph before the theorem.
  \end{proof}

%
\bibliographystyle{andrew}

\end{document}